%% file: paper-inflatable-CaS-2025-R-arxiv.tex
\documentclass[preprint,12pt]{elsarticle_arxiv}

\usepackage{graphicx}
\usepackage{placeins}
\usepackage{subcaption}

\usepackage{amssymb}
\usepackage{amsmath}
\usepackage{amsthm}
\usepackage{amsfonts}

\usepackage{hyperref}
\usepackage{esint}
\usepackage{color}

\journal{Computers \& Structures}

\usepackage{changes}
\definechangesauthor[name={VL}, color=orange]{vl}
\definechangesauthor[name={ER}, color=orange]{er}

\usepackage{algpseudocode}
\usepackage{algorithm}

\algdef{SE}[DOWHILE]{Do}{doWhile}{\algorithmicdo}[1]{\algorithmicwhile\ #1}%
\def\sbb{{\bf{s}}}
\newcommand\revE[1]{{{#1}}}
\newcommand\newSecE[1]{#1}
\newcommand\chE[1]{#1}
\newcommand\trialE[1]{{}} 

\usepackage{appendix}

\newcommand\todosol[1]{{}}

\graphicspath{.}

\newtheorem{remark}{Remark}[section]

\def\eq#1{(\ref{#1})}

\def\bmi#1{\textbf{\textit{#1}}}
\def\ol#1{\overline{#1}}

\def\zerobf{{\bf 0}}
\def\onebm{{\bf{1}}}

\def\psibf{{\boldsymbol{\psi}}}
\def\sigmabf{{\boldsymbol{\sigma}}}

\def\omegabf{{\boldsymbol{\omega}}}

\def\Dlt{\Delta}
\def\dlt{\delta}

\def\veps{\varepsilon}

\def\vphi{\varphi}
\def\vtheta{\vartheta}
\def\pd{\partial}
\def\dd{\,{\rm{d}}}

\def\Om{\Omega}

\def\RR{{\mathbb{R}}}
\def\intY{\fint}

\def\dV{} 
\def\dVx{} 
\def\dVy{} 
\def\dVxy{} 
\def\dY{\mathrm{\,dV}_{\kern-0.1em y}}  %
\def\dS{\mathrm{\,dS}}  %
\def\dSx{\mathrm{\,dS}_{x}}  %
\def\dSy{\mathrm{\,dS}_{y}}  %

\def\fx{{[f]}}
\def\cx{{[c]}}
\def\px{{[P]}}
\def\qx{{[Q]}}

\def\lhs{{l.h.s.{~}}}
\def\ext{{\rm{ext}}}
\def\mic{{\rm{mic}}}
\def\inlet{{\rm{inl}}}
\def\outlet{{\rm{out}}}
\def\visc{{\rm{visc}}}

\def\sgn{{\rm{sgn}}}

\def\inout{{\rm{in-out}}}

\def\aYf#1#2{a_Y^f \left ({#1},\,{#2}\right )}
\def\ipYf#1#2{\left \langle{#1},\,{#2}\right \rangle_{Y_f}}
\def\ipGG#1#2{\left\langle{#1},\,{#2}\right\rangle_{\Gamma_f}}

\newcommand\wrt{{w.r.t.{~}}}
\newcommand\fsi{{\it{fs}}}
\newcommand\ie{{\it{i.e.~}}}
\newcommand\eg{{\it{e.g.~}}}
\newcommand\etc{{\it{etc.~}}}
\newcommand\cf{{{cf.~}}}
\newcommand\wtilde[1]{\widetilde{#1}}
\newcommand\ul[1]{\underline{#1}}
\newcommand\ull[1]{\underline{\underline{#1}}}

\newcommand\Rnil{{}} 

\def\what#1{\widehat{#1}}

\def\Vcal{\mathcal{V}}

\def\Ical{\mathcal{I}}

\def\Qcal{\mathcal{Q}}

\def\Scal{\mathcal{S}}

\def\Ucalbf{\boldsymbol{\mathcal{U}}}
\def\Wcalbf{\boldsymbol{\mathcal{W}}}
\def\Qcalbf{\boldsymbol{\mathcal{Q}}}
\def\Scalbf{\boldsymbol{\mathcal{S}}}

\def\IIcalbf{{\mbox{\boldmath$\mathcal{I}$\unboldmath \kern-0.5em\boldmath$\mathcal{I}$\unboldmath}}}
\def\DDcalbf{{\boldsymbol{\mathcal{I}} \kern-0.5em\boldsymbol{\mathcal{D}}}}

\newcommand\GrxyS[2]{\left(\eebx{#1} + \eeby{#2}\right)}

\newcommand{\jump}[1]{\left\lceil{#1}\right\rfloor}
\newcommand{\jumpG}[1]{{\left\lceil{#1}\right\rfloor}_{\Gamma_f}}
\newcommand{\posPart}[1]{\left[{#1}\right]_{+}}

\newcommand{\Tuf}[1]{{\mathcal{T}}_\veps{\left ({#1}\right )}}
\newcommand\dltsh{\dlt_\tau}

\newcommand\fbav{{\what{\bmi{f}}}}
\newcommand\bbav{{\what{\bmi{b}}}}

\newcommand\tp{{\triangleleft}}
\def\kkk{{k}}

\def\bb{{\bmi{b}}}
\def\ub{{\bmi{u}}}
\def\vb{{\bmi{v}}}
\def\wb{{\bmi{w}}}
\def\fb{{\bmi{f}}}
\def\nb{{\bmi{n}}}

\def\eb{{\bmi{e}}}
\def\ssb{{\bmi{s}}}
\def\qb{{\bmi{q}}}

\def\Bb{{\bmi{B}}}
\def\Db{{\bmi{D}}}

\def\Wb{{\bmi{W}}}
\def\Gb{{\bmi{G}}}
\def\Kb{{\bmi{K}}}

\def\Qb{{\bmi{Q}}}
\def\Ib{{\bmi{I}}}

\def\ubm{{\bf{u}}}

\def\pbm{{\bf{p}}}

\def\fbm{{\bf{f}}}

\def\gbm{{\bf{g}}}

\def\sbm{{\bf{s}}}

\def\Gbm{{\bf{G}}}

\def\Rbm{{\bf{R}}}
\def\Sbm{{\bf{S}}}

\def\Bbm{{\bf{B}}}

\def\Kbm{{\bf{K}}}
\def\Nbm{{\bf{N}}}
\def\Abm{{\bf{A}}}
\def\Mbm{{\bf{M}}}

\def\Jbm{{\bf{J}}}

\def\1bm{{\bf{1}}}

\def\Dltb{\boldsymbol{\Delta}_+^{fc}}

\def\Pibf{\boldsymbol{\Pi}}
\def\Xibf{\boldsymbol{\Xi}}
\def\Hdb{{\bf{H}}^1}
\def\Hpdb{{\bf{H}}_\#^1}
\def\HpdbO{{\bf{H}}_{\#0}^1}
\def\Hdb{{\bf{H}}^1}

\def\Hpdbav{\wtilde{\bf{H}}_\#^1}

\def\Lb{{\bf{L}}}

\def\Dop{{{\rm I} \kern-0.2em{\rm D}}}
\def\Hop{{{\rm I} \kern-0.2em{\rm H}}}
\def\Fop{{{\rm I} \kern-0.2em{\rm F}}}
\def\Aop{{{\rm A} \kern-0.6em{\rm A}}}%
\def\Cop{{{\rm C} \kern-0.6em{\rm C}}}%
\def\Iop{{{\rm I} \kern-0.2em{\rm I}}}%

\def\thetabf{\boldsymbol{\theta}}

\def\eeb#1{\eb({#1})}
\def\eeby#1{\eb_y({#1})}
\def\eebx#1{\eb_x({#1})}

\def\R{\hbox{\rm I\kern-0.2em R}}
\def\Z{\hbox{\rm Z\kern-0.3em Z}}

\def\sx{{[s]}}
\def\fx{{[f]}}
\def\kx{{[k]}}

\def\circparpP#1#2{{\langle{#1}\rangle_{p_{#2}}}}
\def\circpare#1{{\langle{#1}\rangle_\eb}}

\def\aYs#1#2{a_{Y}^{s} \left ({#1},\,{#2}\right )}

\newcommand\Appx[1]{Appendix~\ref{#1}}
\def\prevstep{^{\triangleleft}}
\def\upstep{\ub\prevstep}
\def\ppstep{p\prevstep}

\def\shp{{\rm{sh}}}

\def\bara{\bar} 

\begin{document}

\begin{frontmatter}

  \title{Homogenization of flow in inflatable periodic structures with nonlinear effects}

\author[NTIS]{E.~Rohan\corref{cor1}}
\ead{rohan@kme.zcu.cz}
\cortext[cor1]{Corresponding author}
\author[NTIS]{V.~Luke\v{s}}
\ead{vlukes@kme.zcu.cz}

\address[NTIS]{Department of Mechanics \&
  NTIS New Technologies for Information Society, Faculty of Applied Sciences, University of West Bohemia in Pilsen, \\
Univerzitn\'\i~22, 30100 Plze\v{n}, Czech Republic}

\begin{abstract}
  The paper presents a new type of weakly nonlinear two-scale model of
  inflatable periodic poroelastic structures saturated by Newtonian fluids. The
  periodic microstructures incorporate fluid inclusions connected to the fluid
  channels by admission and ejection valves respected by a 0D model. This
  induces a nonlinearity in the macroscopic Biot-type model, whereby the Darcy
  flow model governs the fluid transport due to the channels. \chE{Moreover, the
  fluid channels consist of compartments separated by  semipermeable membranes
  inducing the pressure discontinuity.} The homogenized model is derived under
  the small deformation assumption, however the equilibrium is considered in the
  Eulerian frame. Deformation-dependent homogenized coefficients of the
  incremental poroelasticity constitutive law and the permeability are
  approximated using the sensitivity analysis, to avoid coupled two-scale
  iterations. Numerical simulations illustrate the inflation process in time.
  Example of a bi-material cantilever demonstrates the inflation induced
  bending. \chE{The proposed two-scale model is intended to provide a
  computational tool for designing of porous metamaterials for fluid transport,
  or shape morphing with various potential applications.}
\end{abstract}

\trialE{
  Recent advances in 3D printing technologies open new possibilities in design of porous metamaterials intended for fluid transport, or shape morphing with various potential applications. The paper presents a new type of weakly nonlinear two-scale model of inflatable periodic poroelastic structures saturated by Newtonian fluids. The periodic microstructures incorporate fluid inclusions connected to the fluid channels by admission and ejection valves respected by a 0D model. This induces a nonlinearity in the macroscopic Biot-type model, whereby the Darcy flow model governs the fluid transport due to the channels. The homogenized model is derived under the small deformation assumption, however the equilibrium is considered in the Eulerian frame. Deformation-dependent homogenized coefficients of the incremental poroelasticity constitutive law and the permeability are approximated using the sensitivity analysis, to avoid coupled two-scale iterations. Numerical simulations illustrate the inflation process in time. Example of a bi-material cantilever demonstrates the inflation induced bending.
  }


\end{frontmatter}



\section{Introduction}

Inflatable structures (InfS) are usually understood as objects made of flat membrane elements designed to attain a desired shape (volume) after their inflation by a fluid (gas) under a given pressure. Another type if InfS is obtained when inflatable membranes are combined with a kind of mesh-like structures, or deformable lattices which determine the morphing process during the inflation -- increase of the volume and internal pressure. Such constructions present the basis for designs of soft robots, namely the robotic arms, ``fingers'', \etc. The design options provide a large variability of InfS which are broadly applicable also as temporary shelters, floating vessels, or safety and medical equipment (inflatable cushions, prosthetic parts). In this context, modelling of InfS has become the topic of numerous publications. These comprise analysis of  inflatable skeletons, or lattices consisting of beams \cite{Sinha-2024} to obtain desired effective medium properties. Such structures are applicable in robotics \cite{Wang-2021}. The nature of inflatable structures require a comprehensive understanding of their mechanical behavior under various dynamic loading regimes. Robustness of the design is important to enhance their functionality, \cf \cite{Sinha-2024}.  In general, a kind of homogenization is the alternative way to treat lattice structures as a continuum \cite{Braides2023,encyclopedia2022}. The homogenization approach was employed in \cite{Ren-2024} as the computational framework for analyzing and designing surface-based inflatable structures with arbitrary fusing patterns from a database of geometrical and corresponding deformation characteristics. 

In this paper, we mean by the InfS a kind of porous fluid-saturated periodic structure which enable to increase its volume and modify its shape depending in the pore pressure value. To this aim, the microarchitecture involves two types of pores, one intended for the fluid transport, \ie fluid redistribution within the whole structures, the other is constituted by ``inflatable inclusions'' connected to the first type porosity locally trough automatic valves operated passively by the pressure difference. This design presents  an extension of the one proposed in \cite{Rohan-Lukes-CaS2024}, where the homogenized model of the fluid redistribution due to peristaltic deformation induced by piezoelectric segments (actuators) was derived and extended to treat nonlinearity issues. In particular, the nonlinearity treated in the present paper stems from the equilibrium considered in the deformed configuration, although we consider only deformations well described using the linear kinematic ansatz. As pointed out above, this nonlinearity is important and even necessary to model functionality of the inflated structures in any case of their contact interaction with other objects -- the situation featuring the robotic applications.



The literature devoted to the use of homogenization in modelling porous media and the fluid structure interaction in periodic, or locally periodic structures is vast. Here we only point out some works of interest which in some sense are related to the particular InfS treated in this paper. Concerning different periodic micro-pores intended for fluid redistribution in a bulk, pores arranged in hexagonal patterns seem to provide optimal performance \cite{Roose-JB2012}. Strongly heterogeneous porous structures involving dual porosities has been studied for decades \cite{Auriault1992,Arbogast1990,rohan-etal-jmps2012-bone,Rohan-zamp2020}; in the context of the problem treated in this paper, mixed dimensional models of blood flow through interfaces between hierarchically arranged porosities coupling 0D-1D-3D vascular models were proposed in \cite{Fritz-Koeppl-Oden-Wohlmuth-2022}, a similar type of models is convenient to describe flows in multiscale fractured-porous media \cite{Panfilov-2021}. Another approach of multicompartment flows \cite{Rohan-etal-JMB2018} involves models of individual compartments which are analogous to the ones obtained by upscaling double porosity microstructures \cite{rohan-cimrman-perfusionIJMCE2010}. Discontinuity in pressure emerges on interfaces between primary and dual porosity (distinguished by scale-dependent contrast in the permeability) as the consequence of reduced interface permeability \cite{Li2018,Daly-Roose-2014}, \cf \cite{GrisoRohan2014} where overall stress discontinuity emerged as a consequence of nonlocal effects associated with force couples imposed at distant points of interconnected interfaces. Homogenization of species transport in problems with jumps in chemical potential was induced by semipermeable membranes in \cite{Ptashnyk-2013}. For studying the
fluid-structure interaction in porous  media under large
deformation in the microstructure, the homogenization approach seems
to be perspective \cite{Sandstrom_2016,Lukes_Rohan_2022}, \revE{\cf \cite{Brown2014,Collis2017,Miller2021}}, although it naturally leads to very complex two-scale computations (the well known ``finite element-square'' (FE$^2$)
approach) requiring reiterated solving the microproblems at all integration nodes of the discretized macroproblem. 
In \cite{Rohan-Lukes-CaS2024}, \cf \cite{Rohan-Lukes-PZ-porel,Rohan-Lukes-nlBiot2015}, the modelling framework for the homogenization-based computationally efficient modelling of weakly nonlinear problems has been proposed. Using the linearization approach based on the sensitivity analysis of the homogenized coefficients \wrt perturbed microconfigurations in response to the local macroscopic fields, the ``deformation-dependent'' homogenized coefficients are obtained using the classical periodic
homogenization \cite{Sanchez1980Book,Cioranescu2008a,Cioranescu-etal-UF-book2018}.


This paper builds on our recent work \cite{Rohan-Lukes-CaS2024} devoted to the peristalsis-driven flow in the homogenized porous medium with piezoelectric actuators generating the deformation (peristaltic) wave.  Therein the necessity of capturing the nonlinearity in the model of the electro-active
porous medium was demonstrated to be crucial to describe the ``propulsion'' effect generating the fluid transport. In the present study, we follow the same approach, such that the equilibrium is verified in the deformed configuration. Potentially, this is important to describe the interaction of the inflatable structure with surrounding objects (robotic applications), besides the perspective extension of the present model by the piezoelectric actuators enabling for the deformation driven fluid redistribution, as reported in \cite{Rohan-Lukes-CaS2024}.

\paragraph{Novelty and paper contribution}
\chE{We intend to provide a modelling tool for designing optimized structures which deform and develop controlable (programmable) interaction with surrounding objects in contact.}
The current work extends the above mentioned model \cite{Rohan-Lukes-CaS2024} by three important aspects: \textbf{a)} the heterogeneous structure involves two different porosity types: the connected channel network, allowing for the fluid redistribution in the whole structure, and the disconnected porosity constituted by periodic arrays of fluid inclusions which can be inflated, thus, increasing locally the structure volume; \textbf{b)} the two porosities are connected by ``fluid fibres'' equipped with automatic pressure gauge valves. These valves enable the fluid flux depending on the positive pressure difference, which induces the model nonlinearity (admission and ejection valves are considered, the latter allowing for deflation of the inclusion, if the overpressure in the inclusion is too high); \textbf{c)} moreover, semipermeable interfaces separating the fluid channels into sub-compartments are considered, which modify the microproblems to be solved in the microscopic representative cell, and influence the homogenized channels permeability. The nonlinearity requires adequate treatment to get numerical solutions for the discretized two-scale model. Up to our best knowledge, any similar model dealing with a periodic porous medium has not been reported in the literature. Nevertheless, as we illustrate in the numerical examples, the model and the considered type of the microstructure have a remarkable application potential, as we explain in the Conclusion.

 \paragraph{Paper organization}
 The fluid-structure interaction problem imposed in the heterogeneous two-phase continuum is introduced in Section~\ref{sec-problem}. For this, the mass conservation of the fluid in the two porosities interconnected by the valves is derived in Section~\ref{sec-equilibr}. This yields the weak formulation of the problem imposed in a fixed configuration (Section~\ref{sec-weak-form}),  which is subject of the homogenization reported in in Section~\ref{sec-homog}. Therein, using the periodic unfolding and formal analysis based on the asymptotic expansions (Section~\ref{sec-asymp}), the local problems yielding the characteristic responses (Section~\ref{sec-loc}), and the macroscopic model are derived (Section~\ref{sec-macro-lin}). Due to the valves embedded in the microstructure, the macroscopic problem is nonlinear, so that an incremental formulation is established in Section~\ref{sec-increment}. 
 Extension of the homogenization result towards problems with the equilibrium imposed in the deformed configuration is explained in Section~\ref{sec-nonlin-II}. For this, the homogenization result involving the homogenized coefficients expressed in terms of the characteristic responses of microstructures is interpreted in terms of the incremental formulation. Following the approach explained in \cite{Rohan-Lukes-CaS2024}, approximation of the
homogenized coefficients depending on the macroscopic state variables is reported in Section~\ref{sec-tan}. 
Numerical implementation of the inflatable metamaterial  using the finite element method is illustrated in Section~\ref{sec:numex}, in terms of several examples showing some important features of the proposed two-scale modelling approach, such as the functionality of the valves, of the semipermeable interfaces, and the effect of the nonlinearity arising from the Eulerian problem formulation leading to the deformation-dependent metamaterial parameters involved in the incremental formulation.
Technical issues related to the model derivations are presented in the Appendix, namely
\chE{the direct numerical simulation for the heterogeneous (nonhomogenized) medium, Appendix~\ref{sec-DNS},},
the sensitivity of the homogenized coefficients, Appendix~\ref{sec-SA}, and, Appendix~\ref{app-discr}, FE-the discretized formulation of the incremental problem reported in Sections~\ref{sec-nonlin-I} and \ref{sec-nonlin-II}. 
\chE{Finally, the survey of used notation in the paper is appended in Appendix~\ref{sec-Notation}}.



\paragraph{Notation}
We employ the following notation. Since we deal with a two-scale problem, we
distinguish the ``macroscopic'' and ``microscopic'' coordinates, $x$ and $y$,
respectively. We use $\nabla_x = (\pd_i^x)$ and $\nabla_y = (\pd_i^y)$ when
differentiation \wrt coordinate $x$ and $y$ is used, respectively, whereby
$\nabla \equiv \nabla_x$. By $\eeb{\ub} = 1/2[(\nabla\ub)^T + \nabla\ub]$ we
denote the strain of a vector function $\ub$, where the transpose operator
is indicated by the superscript ${}^T$. By $\ol{D} = D \cup \pd D$ we denote the closure of an open domain $D\subset \RR^3$ bounded by $\pd D$. The Lebesgue spaces of 2nd-power
integrable functions on an open bounded domain $D\subset \RR^3$ is denoted by
$L^2(D)$, the Sobolev space $\Wb^{1,2}(D)$ of the square integrable
vector-valued functions on $D$ including the 1st order generalized derivative,
is abbreviated by $\Hdb(D)$. Further $\Hpdb(Y)$ is the Sobolev space of
vector-valued Y-periodic functions (the subscript $\#$). \revE{The surface integrals are closed by $\dS,\dSx$, or $\dSy$, while volume integrals are written without an ``integration element'', like $\dd V$}. Positive part $a \in \RR$ is $\posPart{a} = \max\{0,a\}$.


\section{Micromodel of inflatable porous periodic structures}\label{sec-problem}

In the geometrical framework describing the fluid saturated porous medium at
the heterogeneity level, we introduce the fluid-structure interaction (FSI) problem
for which the two-scale problem is obtained by the homogenization method in
Section~\ref{sec-homog}.


We consider a quasi-static loading of a solid elastic skeleton interacting with
a viscous fluid saturating pores which are distributed as connected fluid channels and fluid-saturated ``inclusions'' (capsules)  in the solid phase. \chE{By $\veps = \ell/L$ we denote the scale parameter involved in the asymptotic analysis $\veps\rightarrow 0$, and  reporting the mico- and macr-scopic characteristic lengths, $\ell$ and $L$, respectively.}

 \chE{The skeleton deformation is desribed by the displacement field $\ub^\veps$. We employ the fluid velocity decomposition introduced by virtue of the relative (seepage) velocity $\wb^\veps $
\begin{equation}\label{eq-FS10}
  \begin{split}
    \wb^\veps = \vb^{f,\veps} - \dot{\wtilde\ub}^\veps\;,\\
  \end{split}
\end{equation}
where $\dot{\wtilde\ub}^\veps$ is the solid phase velocity field $\dot\ub^\veps$ extended from $\Om_s^\veps$ to pores $\Om_f^\veps$, and $\wb^\veps = \zerobf$ on the pore walls $\Gamma_\fsi^\veps$.}

\subsection{Geometry}\label{sec-geom}
\begin{figure}[ht]
  \centering
\includegraphics[width=0.7\linewidth]{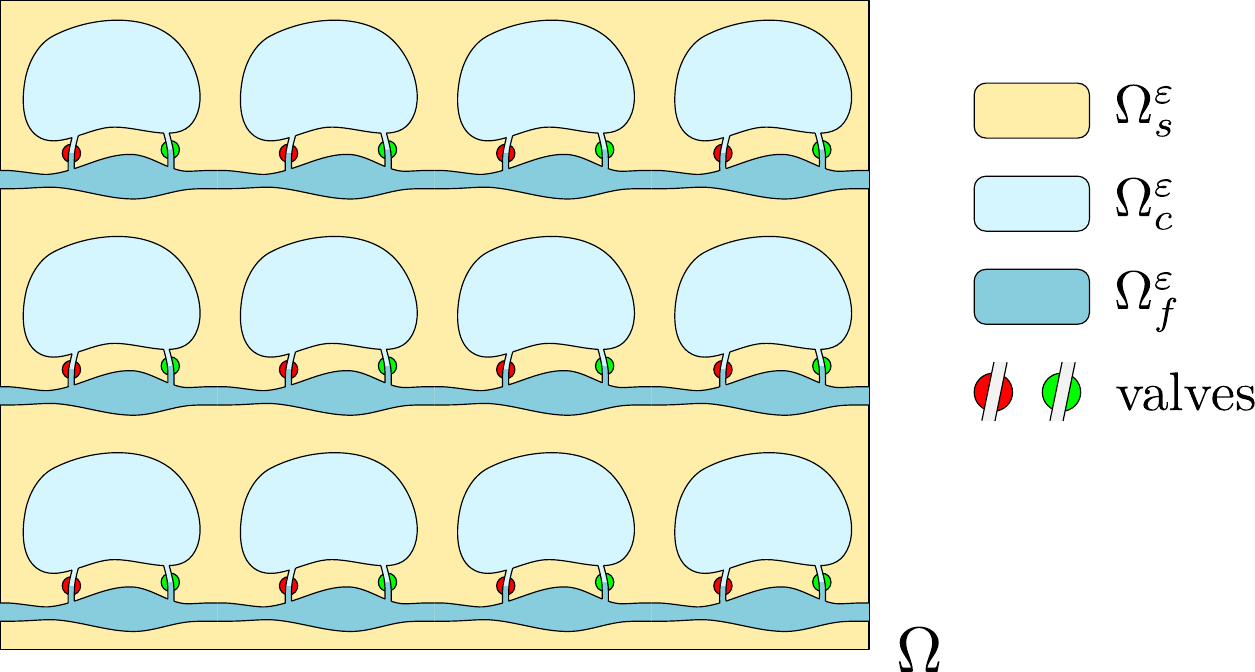}

  \caption{Decomposition of the periodic periodic structure with interconnected
  channels and individual inclusions connected to the channels by the valves.}
  \label{fig:domain-decomposition}
\end{figure}

The poroelastic medium occupies an open bounded domain $\Om \subset \RR^3$
with Lipschitz boundary $\pd \Om$.
The decomposition of the periodic microstructure is displayed in Fig.~\ref{fig:domain-decomposition}.
By $\Gamma_f^\veps$ we denote the solid-fluid interface, $\Gamma_f^\veps  = \ol{\Om_s^\veps} \cap \ol{\Om_f^\veps}$.
Domain $\Om$ is decomposed into the solid matrix,
$\Om_s^\veps$, 
and 
the fluid occupying the channels, $\Om_f^\veps$, and fluid inclusions, $\Om_c^\veps$.
\begin{equation}\label{eq-pzm1}
\begin{split}
\Om = \Om_f^\veps \cup \Om_c^\veps \cup \Om_s^\veps \cup \chE{\Gamma_f^\veps}\;,
\quad \Om_f^\veps \cap \Om_c^\veps \cap \Om_s^\veps = \emptyset\;.\\
\end{split}
\end{equation}
We assume that both the matrix $\Om_s^\veps$ and the fluid-filled channels $\Om_f^\veps$ are connected domains, whereby $\Om_c^\veps$ is the union of separated fluid inclusions. Below we employ the index sets $\Ical_c^\veps$ and $\Ical_f^\veps$, such that $\Om_c^\veps = \bigcup_{\kkk\in\Ical_c^\veps}{\Om_c^{\kkk,\veps}}$. An analogous decomposition holds for channels $\Om_f^\veps$, as defined below in \eq{eq-ifc1}.
The valve connect the two type porosities, $\Om_f^\veps$ and $\Om_c^\veps$.

\paragraph{External boundary decomposition}
We denote by $\pd_\ext\Om_s^\veps = \ol{\Om_s^\veps} \cap \pd\Om$ the
exterior boundaries of $\Om_s^\veps$. In analogy, we define
$\pd_\ext\Om_f^\veps$ as the exterior parts on the
boundaries of the fluid $\Om_f^\veps$. 

The following two splits are defined,
$\pd_\ext\Om_{s}^\veps = \Gamma_{\ub}^\veps \cup \Gamma_\sigma^\veps$ and
$\pd_\ext\Om_f^\veps = \Gamma_p^\veps \cup \Gamma_w^\veps$ 
such that the ``Dirichlet'' type  boundary conditions are prescribed on $\Gamma_{\ub}^\veps$ and $\Gamma_p^\veps$, whereas the surface traction forces and fluid fluxes are prescribed on $\Gamma_\sigma^\veps$ and $\Gamma_w^\veps$, respectively, where
\begin{equation}\label{eq-pzm1a}
\Gamma_\sigma^\veps = \pd_\ext\Om_{s}^\veps\setminus\Gamma_{\ub}^\veps\;, \quad 
\Gamma_w = \pd_\ext\Om_f^\veps\setminus\Gamma_p^\veps\;.
\end{equation}

\paragraph{Valves -- fibre junctions between the channels and the inclusions}
Each inclusion is connected to the fluid channel by two fluid fibres (thin channels) $\bar\Gamma_V^\veps$, $V = A,E$ equipped \chE{by ``admission'' (A) and ``ejection'' (E) } valves.
Let $\hat x_V^\veps \in \Gamma_\fsi^\veps\cap\bar\Gamma_V^\veps$ be a position where the fibre $\bar\Gamma_V^\veps$ flows out in the channel.
In analogy the fibre penetrates in the inclusion at a point $\hat x_V^\veps$. 
The effective (microscopic) fluid velocity relative to the solid phase describes
the flow through the fibre contributing to the mass balance in the channel by
\begin{equation}\label{eq-2}
  \begin{split}
    \bar w_A^\veps =  \kappa_A^\veps
    \posPart{\hat p_f^\veps(\hat x^A,\cdot) - p_c^\veps}\;,\\
    \bar w_E^\veps = \kappa_E^\veps
    \posPart{p_c^\veps - \hat p_f^\veps(\hat x^E,\cdot)  - \Dlt P_E}\;, 
\end{split}
\end{equation}
where $\kappa_A^\veps$ and $\kappa_E^\veps$ are permeabilities of the (open) admission and ejection channels, respectively; recall $\posPart{a}$ returns the positive part of $a$. The ejection pressure threshold $\Dlt P_E$ determines the  pressure drop $p_c^\veps - p_f^\veps >  \Dlt P_E$ enabling for the inclusion fluid evacuation. In \eq{eq-2}, $\hat p_f^\veps(\hat x^A,\cdot)$ and $\hat p_f^\veps(\hat x^E,\cdot)$ are evaluated at points $\hat x^A$, and  $\hat x^E$ on the channel surface $\Gamma_\fsi^\veps$ from where the ``fibres'' lead to the inclusion. Note that $p_c^\veps$ are spatially constant within any inclusion.

\todosol{definovat $\posPart{\cdot}$}

\begin{remark}\label{rem-distribution}
To avoid point values of the pressure $p_f^\veps$ in the channel, we use the Dirac distribution $\hat\delta_V^\veps(x)$, $x \in  \Gamma_\fsi^\veps$ defined at $\hat x_V^\veps$, such that the point value is established through $\intY_{\Gamma_\fsi}\hat\delta_V^\veps p_f^\veps \dS = \hat p_f^\veps(\hat x_V^\veps,\cdot)$.
\end{remark}

\paragraph{Semipermeable interfaces in the channels}
\chE{These iterfaces denoted by $\Gamma_f^\veps$ represent semi-permeable membranes, 
or a kind of ``sieves'' causing the pressure discontinuity, while the fluid velocity remains to be continuous.
Due to $\Gamma_f^\veps$ being embedded in the fluid channels, $\Om_f^\veps$ are decomposed into subdomains (sub-compartments) $\{\Om_f^{k,\veps}\}_k$, such that}
\begin{equation}\label{eq-ifc1}
  \begin{split}
    \Om_f^\veps & = \bigcup_{\kkk\in\Ical_f^\veps}{\Om_f^{\kkk,\veps}} \cup \Gamma_f^\veps\;,\quad
    {\Om_f^{\kkk,\veps}}\cap {\Om_f^{\kkk+1,\veps}} = \emptyset\;,\\
  \Gamma_f^{k,\veps} & = \pd{\Om_f^{k,\veps}} \cap\pd{\Om_f^{k+1,\veps}} \subset \Gamma_f^\veps\;,\\
  \Gamma_f^\veps & = \bigcup_\kkk{\Gamma_f^{\kkk,\veps}}\;.
  \end{split}
\end{equation}
\chE{Interfaces $\Gamma_f^{k,\veps}$ are characterized by permeability $\hat\kappa_f^\veps$ such that the pressure discontinuity is proportional to the transverse flow, }
\begin{equation}\label{eq-3}
  \begin{split}
     w_n^{k,\veps} & =  -\hat\kappa_f^\veps\jump{p_f^\veps}^{k,\veps}\;,\quad \mbox{ on }\Gamma_f^{k,\veps}\;,\\
    \mbox{ where }\quad \jump{q}^{k,\veps} & = q|_{k+1} - q|_{k}\;,\quad \mbox{ on }\Gamma_f^{k,\veps}\;.
\end{split}
\end{equation}

\begin{remark}\label{rem-channel}
  The decomposition \eq{eq-ifc1}, in the context of the subdomain ordering being related to the ordering of interfaces $\Gamma_f^{k,\veps}$, is pertinent for a 1D topology, but  can be generalized accordingly to the specific channel topologies in 2D and 3D. 
\end{remark}

\subsection{Micromodel -- linear setting}\label{sec-micromodel}
The constitutive equations for the stress $\sigmabf^\veps = (\sigma_{ij}^\veps)$ describe an elastic solid in
 $\Om_s^\veps$, and viscous compressible fluid in $\Om_f^\veps$
\begin{equation}\label{eq-1}
  \begin{split}
  \sigmabf^\veps(\ub^\veps) & = \Aop^\veps \eeb{\ub^\veps}\;,\\
\sigmabf_f^\veps & = -p_f^\veps \Ib + \Dop^{f,\veps} \eeb{\vb^{f,\veps}}\;, \quad
\end{split}
\end{equation}
where $\Aop^\veps=(A_{ijkl}^\veps)$ is the elasticity fourth-order symmetric
positive definite tensor of the solid, $A_{ijkl} = A_{klij} = A_{jilk}$, and  Newtonian fluid is described by the viscous stress rheology given by \chE{the viscosity tensor $\Dop_f^\veps = (D_{klij}^{f,\veps})$ defined through}
$D_{klij}^{f,\veps} = \mu^\veps (\dlt_{ik}\dlt_{jl} + \dlt_{il}\dlt_{jk} - (2/3)\dlt_{ij}\dlt_{kl})$, where  $\mu^\veps$ is the dynamic viscosity.

\subsubsection{Equilibrium and mass conservation}\label{sec-equilibr}
We disregard inertia effects in both the solid and the fluid phases. 
The static equilibrium in the solid  and the quasistatic approximation of the momentum balance in the fluid yield the following differential equations
\begin{equation}\label{eq-pzfl0}
\begin{split}
-\nabla\cdot\Aop^\veps \eeb{\ub^\veps}  & = \fb^{s,\veps}\;, \quad \mbox{ in }  \Om_s^\veps\;, \\
 -\nabla \cdot (\Dop^{f,\veps}\eeb{\vb^{f,\veps}}  + \nabla p_f^\veps & = \fb^{f,\veps} - (\wb^\veps \cdot\nabla)\vb^{f,\veps} \;,\quad \mbox{ in } \Om_f^\veps\setminus\Gamma_f^\veps\;.\\
\end{split}
\end{equation}
where $\gamma$ is the fluid compressibility. 
\chE{We consider slow flows, we neglect also the advection term in the fluid equilibrium \eq{eq-pzfl0}$_2$.}

\paragraph{Solid-fluid interface related to the flow in channels}
Due to presence of interfaces $\Gamma_f^{k,\veps}$ inducing the pressure discontinuity, it is necessary te respect the additional loading transmitted to the solid skeleton. Hence, extended solid-fluid interface is identic the compartment boundary $\pd \Om_f^{k,\veps}$, except of those for which $\pd \Om_f^{k,\veps}\cap \pd \Om \not = \emptyset$. Therefore, we define
\begin{equation}\label{eq-FS18}
  \begin{split}
    \Scal_\sigma^\veps(\sigmabf_f^\veps,\vb^\veps) & =  \sum_{k=1}^N\int_{\pd \Om_f^{k,\veps}} \nb^\sx\cdot\sigmabf_f^\veps\cdot\vb^\veps\dS - \int_{\pd_\ext \Om_f^\veps} \nb^\sx\cdot\sigmabf_f^\veps\cdot\vb^\veps\dS\;,\\
    \quad\mbox{ where }\quad
    \sigmabf_f^\veps  & = -{p}^\veps \Ib + \sigmabf_f^{\visc,\veps}\;,\quad 
    \chE{ \sigmabf_f^{\visc,\veps}  =  \Dop_f^\veps\eeb{\wb^\veps + \dot{\wtilde\ub}^\veps}\;.}
 \end{split}
\end{equation}
Note that the interface of the skeleton with the fluid inclusions is not involved in $\Scal_\sigma^\veps$.

\begin{remark}\label{rem-jump}
  \chE{Due to the velocity continuity on $\Gamma_f^{k,\veps}$, since its motion is described by $\dot{\wtilde\ub}^\veps$, hence $\jump{\dot{\wtilde\ub}^\veps} = 0$, by virtue of  \eq{eq-3} we get
\begin{equation}\label{eq-3a}
  \begin{split}
    \jump{\wb^\veps} = \jump{ \vb^{f,\veps}} = \zerobf\;, \quad \mbox{ and } \wb^\veps\cdot\nb^\kx = w_n^{k,\veps}  = -\hat\kappa_f^\veps\nb^\kx \jump{p_f^\veps}^{k,\veps}\;,
\end{split}
\end{equation}
where $\nb^\kx$ is the unit normal vector on $\Gamma_f^{k,\veps}$ outward to $\Om_f^{k,\veps}$, so pointing towards $\Om_f^{k+1,\veps}$, and $q|_k$ is the trace of $q$ on $\pd \Om_f^{k,\veps}$.
Due to the pressure drop on $\Gamma_f^{k,\veps}$, only the viscous part of the stress $\sigmabf_f^{\visc,\veps}$ is continuous on $\Gamma_f^{k,\veps}$, while the total stress is not,
\begin{equation}\label{eq-3b}
  \begin{split}
   \jump{\sigmabf_f^{\visc,\veps}}\cdot\nb^\kx = \zerobf\;.
\end{split}
\end{equation}
Consequently the discontinuity in the traction
$\bb^{k,\veps} = \jump{\sigmabf^\veps}\cdot\nb^\kx = -\jump{p_f^\veps}\nb^\kx$ is transmitted by the ``sieve'' anchored in the skeleton $\Om_s^\veps$. This effect is included in the definition of the virtual power  $\Scal_\sigma^\veps(\sigmabf_f^\veps,\vb^\veps)$ in \eq{eq-FS18}.
  }
\end{remark}

\chE{
\paragraph{Weak form of the fluid flow in channels}
Due to the pressure discontinuity on $\Gamma_f^{k,\veps}$, when $\hat\kappa$ is finite, the derivation of the weak form of the \lhs terms in \eq{eq-pzfl0}$_2$ yields (upon the multiplication by the test velocity $\thetabf^\veps$ vanishing on $\Gamma_\fsi^\veps$ and continuous on $\Gamma_f^{k,\veps}$),
\begin{equation}\label{eq-FS19a}
   \begin{split}
 \sum_{k=1}^N\int_{\Om_{f}^{k,\veps}} \left(\sigmabf_f^{\visc,\veps} 
 - p_f^\veps\Ib\right) : \eeb{\thetabf^\veps} \dV
 + \int_{\Gamma_{f}^\veps}\left(\jump{p_f^\veps}\thetabf^\veps\cdot\nb^\Gamma - \jump{\sigmabf_f^{\visc,\veps}}:\nb^\Gamma\otimes\thetabf^\veps \right) \dS \\
 =\int_{\Om_{f}^{\veps}} \left(\sigmabf_f^{\visc,\veps} 
 - p_f^\veps\Ib\right) : \eeb{\thetabf^\veps} \dV + \int_{\Gamma_{f}^\veps}\frac{1}{\hat\kappa^\veps} \wb^\veps\cdot\nb^\Gamma\otimes\nb^\Gamma \cdot \thetabf^\veps\dS\;,
 \end{split}
\end{equation}
where $\nb^\Gamma$ is the unit normal on $\Gamma_{f}^\veps$ and  \eq{eq-3a} was employed; note the normal  $\nb^\kx$ orientation established in Remark~\ref{rem-jump}. The jump of the viscous stress traction $\jump{\sigmabf_f^{\visc,\veps}}\cdot \nb^\Gamma$ vanishes due to \eq{eq-3b}. It is worth noting that for large values of $\hat\kappa_f$ (approaching $+\infty$) the standard model is recovered without any interface effect.
}


\paragraph{Mass conservation in inclusions}
  To derive the mass conservation in the fluid phase, we first consider the mass balance in any fluid inclusions $\Om_c^{k,\veps}$. Pressure $p_c^\veps \in \Qcal_c^\veps(\Om_c^\veps)$ is constant in each  $\Om_c^{k,\veps} \subset \Om_c^\veps$; the space of piecewise constant pressures is defined below, when introducing the weak formulation, see \eq{eq-5a}. 
Each $\Om_c^{k,\veps}$ is connected by the above mentioned valves with the fluid channels $\Om_f^\veps$, which is expressed by velocity $\bar w_V^\veps$ defined in \eq{eq-2} which, however, is defined pointwise. Therefore, in analogy with distributions $\hat\delta_V^\veps$ on the channel surface, we employ the (periodic) distributions $\bar\dlt_V^\veps$, $V = A,E$, located at $\hat x_V^\veps$ on surface $\pd \Om_c^{k,\veps}$ of any inclusion (labelled by $k$).
\begin{equation}\label{eq-Mc}
  \begin{split}
    \int_{\pd\Om_c^{k,\veps}} \nb^\cx\cdot\dot\ub^\veps + \gamma |\Om_c^{k,\veps}| \dot p_c^\veps & = \int_{\pd\Om_c^{k,\veps}} (\bar\dlt_A^\veps\bar w_A^\veps - \bar\dlt_E^\veps\bar w_E^\veps)\dS \\
 \int_{\Om_c^{k,\veps}} \nabla\cdot\dot{\wtilde{\ub^\veps}} + \gamma |\Om_c^{k,\veps}| \dot p_c^\veps    & = \int_{\pd\Om_c^{k,\veps}} \bar\dlt_A^\veps\kappa_A^\veps \posPart{\hat p_f^\veps - p_c^\veps} -  \int_{\pd\Om_c^{k,\veps}} \bar\dlt_E^\veps\kappa_E^\veps \posPart{p_c^\veps - \hat p_f^\veps - \Dlt P_E} \;,
  \end{split}
\end{equation}
where $\gamma$ is the fluid compressibility.
Note that $\dot\ub$ is the solid velocity (``dot'' means the material derivative) and $\wtilde\ub$ denotes a smooth extension of $\ub$ from the solid to the pore fluid domain.

\paragraph{Mass conservation in channels}
Recall that the channels are constituted by compartments $\Om_f^{k,\veps}$, $k = 1,\dots,N_\veps$, separated by $\Gamma_f^{k,\veps}$, as defined in \eq{eq-ifc1}. 
The mass conservation for a slightly compressible fluid is described by continuity equation,
\begin{equation}\label{eq-Mf1}
  \begin{split}
     \nabla\cdot\vb^{f,\veps} + \gamma \dot p_f^\veps & = 0\;,\\
    \nabla\cdot\wb^\veps + \nabla\cdot\dot{\wtilde{\ub}}+ \gamma \dot p_f^\veps & = 0\;,\mbox{ in } \Om_f^{k,\veps}\;.
  \end{split}
\end{equation}
Let us first disregard the presence of the valves connecting the channel compartments with the inclusions, see \eq{eq-2} (valves closed).
Considered \eq{eq-Mf1} in the weak sense, due to the discontinuity of  $p_f^\veps$ and the associated test function on the interface, the mass conservation within each $\Om_f^{k,\veps}$ yields
\begin{equation}\label{eq-Mf2}
  \begin{split}
 \int_{\Om_f^{k,\veps}}\gamma \dot p_f^\veps q^\veps   + \int_{\Om_f^{k,\veps}}(q^\veps \nabla\cdot\dot{\wtilde{\ub}} - \wb^\veps \cdot\nabla q^\veps ) + \int_{\Gamma_f^{k-1,\veps}}  q^\veps\wb^\veps\cdot\nb^{\kx} + \int_{\Gamma_f^{k,\veps}}  q^\veps\wb^\veps\cdot\nb^{\kx} = 0\;,
  \end{split}
\end{equation}
\chE{for all test pressures $q^\veps\in \Qcal_f^\veps(\Om_f^\veps)$, see below \eq{eq-5a}, which are discontinuous on any $\Gamma_f^{k,\veps}$},
whereby  $\wb^\veps$ is defined in \eq{eq-3}, and $\nb^{\kx}$ 
is the outward normal to $\Om_f^{k,\veps}$.
By virtue of Remark~\ref{rem-channel}, we can sum up \eq{eq-Mf2} over all subdomains labelled by $k$
\begin{equation}\label{eq-Mf3}
  \begin{split}
    \sum_{k=1}^N\left(\int_{\Om_f^{k,\veps}}\gamma \dot p_f^\veps q^\veps   + \int_{\Om_f^{k,\veps}}(q^\veps \nabla\cdot\dot{\wtilde{\ub}} - \wb^\veps \cdot\nabla q^\veps ) \right)+
    \sum_{k=1}^{N-1}  \int_{\Gamma_f^{k,\veps}} \kappa_f^\veps \jump{p_f^\veps}^{k,\veps} \jump{q^\veps}^{k,\veps} \\
    = \int_{\Gamma_\inout} \wb^\veps\cdot\nb^\fx q^\veps\;,
  \end{split}
\end{equation}
\chE{where $\Gamma_\inout$ represents the inlet and outlet surfaces of the connected porosity.}  
Now this equality must be modified due the valves defined in \eq{eq-2}. Let $\hat x_V^\veps \in \Gamma_\fsi^\veps$ be a position where the fibre $\bar\Gamma_V^\veps$ flows out in the channel. Using the Dirac distribution $\hat\delta_V^\veps(x)$, $x \in  \Gamma_\fsi^\veps$ defined at $\hat x_V^\veps$, $V=A,E$, flow through the fibre contributes to the mass balance in the channel by
\begin{equation}\label{eq-Mf4}
  \begin{split}
    \int_{\Gamma_\fsi^\veps} q^\veps \bar w_A^\veps \hat\delta_A^\veps \dS = \int_{\Gamma_\fsi^\veps} q^\veps \kappa_A^\veps \posPart{p_f^\veps - p_c^\veps} \hat\delta_A^\veps \dS\;,\\
     \int_{\Gamma_\fsi^\veps} q^\veps \bar w_E^\veps \hat\delta_E^\veps \dS = \int_{\Gamma_\fsi^\veps} q^\veps \kappa_E^\veps \posPart{p_c^\veps - p_f^\veps  - \Dlt P_E} \hat\delta_E^\veps \dS\;.
  \end{split}
\end{equation}
The flow through the admission and evacuation (ejection) fibres is oriented due to the valves, such that  \eq{eq-Mf3} is modified, as follows 
\begin{equation}\label{eq-Mf5}
  \begin{split}
  & \int_{\Om_f^{\veps}}\gamma \dot p_f^\veps q^\veps   + \int_{\Om_f^{\veps}}(q^\veps \nabla\cdot\dot{\wtilde{\ub}} - \wb^\veps \cdot\nabla q^\veps ) +
   \int_{\Gamma_f^{\veps}} \kappa_f^\veps \jump{p_f^\veps}^{\veps} \jump{q^\veps}^{\veps} \dS \\
 &  = \int_{\Gamma_\fsi^\veps} q^\veps \kappa_A^\veps \posPart{p_f^\veps - p_c^\veps} \hat\delta_A^\veps \dS
  - \int_{\Gamma_\fsi^\veps} q^\veps \kappa_E^\veps \posPart{p_c^\veps - p_f^\veps  - \Dlt P_E} \hat\delta_E^\veps \dS\\
 &   + \int_{\Gamma_\inout} \wb^\veps\cdot\nb^\fx q^\veps \dS\;.
  \end{split}
\end{equation}

\subsubsection{Boundary conditions (BCs)}\label{sec-bcs}
BCs prescribed for the solid involve the Dirichlet conditions concerning the fields $\ub^\veps,\vb^{f,\veps}$ and $\vphi^\veps$, whereas the Neumann type conditions comprise the boundary tractions $\sigmabf_s^\veps\cdot \nb = \bb^{s,\veps}$ acting on $\Gamma_\sigma$ of the solid surface, the fluid pressure $p^\veps$ representing the stress in the fluid acting on $\Gamma_p^\veps$ (assuming the dissipative part of $\sigmabf_f^\veps\cdot \nb$ vanishes), so that 
\begin{equation}
  \begin{split}\label{eq-BC1}
    \ub^\veps  = \ub_\pd \quad \mbox{ on } \Gamma_u^\veps\;,& \quad \sigmabf_s^\veps\cdot \nb = \bb^{s,\veps}\mbox{ on } \Gamma_\sigma\;, \\
       \vb^{f,\veps}  = 0  \quad \mbox{ on } \Gamma_w^\veps\;,& \quad  p^\veps  = p_\pd \quad \mbox{ on } \Gamma_p^\veps\;. \\
    \end{split}
\end{equation}

\subsubsection{Weak formulation for the quasistatic flow}\label{sec-weak-form}
The homogenization using the periodic unfolding method is applied to the weak formulation of the quasistatic fluid-solid interaction. The admissibility sets of displacements, fluid velocities, and fluid pressures in the two porosity types are defined in  terms of $H^1$ and $L^2$ spaces, whereby 
the boundary conditions \eq{eq-BC1} are respected,
\begin{equation}\label{eq-5a}
\begin{split}
\Ucalbf^\veps(\Om_{s}^\veps) & = \{\vb \in \Hdb(\Om_{s}^\veps)|\; \vb = \ub_\pd \mbox{ on }  \Gamma_{\ub}^\veps\}\;,\\
\Wcalbf_0^\veps(\Om_f^\veps,\Gamma_\fsi^\veps) & = \{\vb \in \Hdb(\Om_f^\veps)|\;\vb = \zerobf \mbox{ on } \Gamma_\fsi^\veps\cup \Gamma_{w}^\veps\}\;,\\
\Qcal_c^\veps(\Om_c^\veps) & = \{ q \in L^2(\Om_c^\veps)|\;q \mbox{ is constant in each } \Om_c^{k,\veps}\;,\; k \in I_c^\veps\}\;,\\
\Qcal_f^\veps(\Om_f^\veps)  & = \{q \in L^2(\Om_f^\veps)|\;q \in H^1(\Om_f^{k,\veps})\;,\; k \in \Ical_f^\veps\;,
q = p_\pd \mbox{ on }  \Gamma_p^\veps\}\;,
\end{split}
\end{equation}
recalling the porosity partitioning into fluid inclusions and fluid channel compartments, as introduced in \eq{eq-ifc1}.
\chE{Space $\Qcal_c^\veps(\Om_c^\veps)$ is introduced to adhere with general weak form involving ``test functions''. This enables to rewrite  \eq{eq-Mc} using one equality with test pressures serving as switches, as employed below in \eq{eq-wf3}.}

\chE{We use the particular results derived in Section~\ref{sec-equilibr}, namely \eq{eq-Mc} and \eq{eq-Mf5}. Obvious treatment of \eq{eq-pzfl0}$_1$ using test functions and integrating by parts, yields the weak from of the equilibrium in the solid part, whereas \eq{eq-FS19a} is employed to obtain the equilibrium in the fluid channel.}

\chE{Solutions of the quasistatic  problem for}
$(\ub^\veps(t,\cdot),\wb^\veps(t,\cdot),p_f^\veps(t,\cdot),p_c^\veps(t,\cdot)) \in \Ucalbf^\veps(\Om_{s}^\veps) \times \Wcalbf_0^\veps(\Om_f^\veps,\Gamma_\fsi^\veps)\times\Qcal_f^\veps(\Om_f^\veps)\times\Qcal_c^\veps(\Om_c^\veps)$ satisfy at $t>0$ the following weak forms of the conservation laws: \\
$\square$\emph{Quasistatic equilibrium in the solid and fluid phases},
\begin{equation}\label{eq-wf1}
\begin{split}
  &  \int_{\Om_{s}^\veps} \left(\Aop^\veps\eeb{\ub^\veps}  \right): \eeb{\vb^\veps} \dV
  -\Scal_\sigma^\veps(\sigmabf_f^\veps,\vb^\veps)
 + \int_{\pd \Om_c^\veps} p_c^\veps\nb^\sx\cdot\vb^\veps \dS
  =\int_{\Om_{s}^\veps} \fb^{s,\veps}\cdot\vb^\veps \dV + 
  \int_{\pd_\sigma\Om_{s}^\veps} \bb^{s,\veps}\cdot\vb^\veps \dV\;,\\
 &  \int_{\Om_{f}^\veps} \left(\Dop_f^\veps\eeb{\wb^\veps + \dot{\tilde\ub}^\veps}
   - p^\veps\Ib\right) : \eeb{\thetabf^\veps} \dV + \chE{\int_{\Gamma_{f}^\veps}\frac{1}{\hat\kappa^\veps} \wb^\veps\cdot\nb^\Gamma\otimes\nb^\Gamma \cdot \thetabf^\veps\dS}
   = \int_{\Om_{f}^\veps} \fb^{f,\veps}\cdot\thetabf^\veps \dV \;,
 \end{split}
\end{equation}
for all  $(\vb^\veps,\thetabf^\veps,q^\veps) \in  \Ucalbf_0^\veps(\Om_{s}^\veps) \times \Wcalbf_0^\veps(\Om_f^\veps,\Gamma_\fsi^\veps)$  (we recall the interface loading \eq{eq-FS18} due to the flow in $\Om_f^\veps$);\\
$\square$\emph{Mass conservation of the fluid in the channels},
\begin{equation}\label{eq-wf2}
\begin{split}
& \int_{\Om_f^{\veps}}\gamma \dot p^\veps q^\veps   + \int_{\Om_f^{\veps}}(q^\veps \nabla\cdot\dot{\wtilde{\ub}} - \wb^\veps \cdot\nabla q^\veps ) +
   \int_{\Gamma_f^{\veps}} \kappa_f^\veps \jump{p_f^\veps}^{\veps} \jump{q^\veps}^{\veps} \dS \\
  & = \int_{\Gamma_\fsi^\veps} q^\veps \kappa_A^\veps \posPart{p_f^\veps - p_c^\veps} \hat\delta_A^\veps \dS
  - \int_{\Gamma_\fsi^\veps} q^\veps \kappa_E^\veps \posPart{p_c^\veps - p_f^\veps  - \Dlt P_E} \hat\delta_E^\veps \dS\\
  &  + \int_{\Gamma_\inout} \wb^\veps\cdot\nb^\fx q^\veps \dS\;,
\end{split}
\end{equation}
for all $q^\veps\in\Qcal_{0f}^\veps(\Om_f^\veps)$, such that  $\Qcal_{0f}^\veps$ is defined as $\Qcal_{f}^\veps$ in \eq{eq-5a}, but with the zero BCs, \ie putting $p_\pd = 0$;\\
$\square$\emph{Mass conservation of the fluid in the inclusions},
\begin{equation}\label{eq-wf3}
\begin{split}  
   &\int_{\Om_c^\veps} q_c^\veps\left(\nabla\cdot\dot{\wtilde{\ub^\veps}}  + \gamma \dot p_c^\veps \right)
  =  
  \int_{\pd\Om_c^\veps} (w_A^\veps \bar\dlt_A^\veps - w_E^\veps \bar\dlt_E^\veps) \dS\\
  & = \int_{\pd\Om_c^\veps} \left(\kappa_A^\veps\bar\dlt_A^\veps \posPart{\hat p_f^\veps - p_c^\veps}q_c^\veps -  \kappa_E^\veps\bar\dlt_E^\veps \posPart{p_c^\veps - \hat p_f^\veps - \Dlt P_E}q_c^\veps\right) \dS \;,
\end{split}
\end{equation}
for all  $q_c^\veps \in \Qcal_c^\veps(\Om_c^\veps)$ (piecewise constant functions). Note that $\bar\dlt_{A,E}^\veps$ are periodic, providing support on $\pd\Om_c^{k,\veps}$ of each inclusion $k \in \Ical_c^\veps$.


\section{Homogenization}\label{sec-homog}
Following the upscaling approach based on the asymptotic homogenization by the unfolding method \cite{Cioranescu-etal-UF-book2018}, 
we consider the linear problem defined in Section~\ref{sec-weak-form}. In  \cite{Rohan-Lukes-CaS2024}, it has been shown how the nonlinear problem associated with the equilibrium in the Eulerian frame can be solve using  a sequence of
linear subproblems associated with an incremental formulation. \chE{We apply the asymptotic homogenization based on the periodic unfolding  $\Tuf{}$ operator introduced in \cite{Cioranescu2008a,Cioranescu-etal-UF-book2018}, which enables to present any function defined over a periodic structure in terms of the macroscopic (slow) coordinates $x$ and the microscopic (fast) coordinates $y$.}

\begin{remark}
\chE{Although, in this paper, we consider small deformations of the solid and the linear constitutive law \eq{eq-1}, we provide also an extension of the model to consider  the nonlinearity of the FSI on the interface due to the deforming pores, thereby modifying the hydraulic permeability and other effective properties of the homogenized medium. Consistently with our previous work \cite{Rohan-Lukes-CaS2024}, where the linearized incremental formulation was derived in detail by differentiation of the equilibrium and mass conservation equations in the Eulerian frame, in Section~\ref{sec-nonlin-II}, we provide an extension of the homogenized model obtained for the new micromodel, as described above.}
  \end{remark}

\subsection{Periodic structure}\label{sec-periodic}
\begin{figure}[ht]
	\centering
  \includegraphics[width=0.95\linewidth]{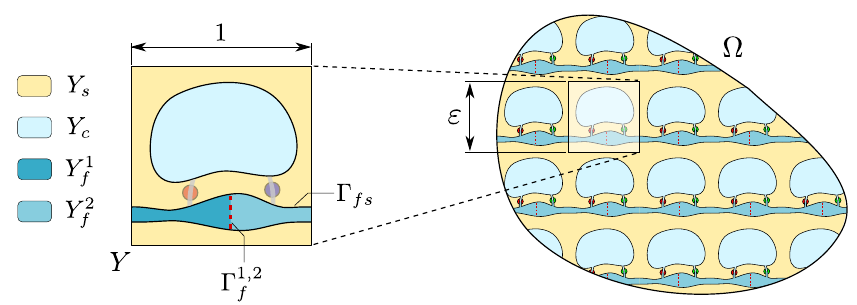}
  \caption{\chE{Left:} The reference unit periodic cell $Y$ and its decomposition, fibers equipped with the vales, connecting the two porosities; $\Gamma_f^{1,2}$ is the semipermeable interface \chE{in the channel $Y_f$}. Right: The periodic structure with interconnected channels and individual inclusions connected to the channels by the valves. \chE{(Note that these valves are not geometrically distinguished in $Y$, see Remark~\ref{rem-distribution}.)} }
  \label{fig-Y-per}
\end{figure}

\paragraph{Reference periodic unit cell}
The heterogeneous medium is generated as a periodic lattice by a reference periodic cell $Y$
decomposed into three non-overlapping subdomains: the solid skeleton $Y_s$, the fluid channel $Y_f$ and the inclusion $Y_c$, see
Fig.~\ref{fig-Y-per}, whereby $|Y|\approx 1$. 
\begin{equation}\label{eq-6}
\begin{split}
  Y &= Y_s \cup Y_f \cup Y_c \cup \Gamma_Y\;, \quad
   Y_i \cap Y_j = \emptyset\, \mbox{ for } j\not = i \mbox{ with } i,j \in \{f,c,s\}\;,\\
   \Gamma_\fsi &= \pd Y_s \cap \pd Y_f\;.
\end{split}
\end{equation}
Channel $Y_f$ is subdivided into compartments $Y_f^k$, here labelled by superscript $k$, in the context of Remark~\ref{rem-channel}, separated by interfaces $\Gamma_f^{k,k+1}$, such that
\begin{equation}\label{eq-6b}
  \ol{Y}_f =\bigcup_k  \ol{Y}_f^k\;,\quad Y_f^k \cap   Y_f^{k+1} = \emptyset\;, \quad
  \Gamma_f^{k,k+1} =  \ol{Y}_f^k\cap  \ol{Y}_f^{k+1}\;.
\end{equation}
In what follows, by $\Gamma_f$ we refer to the union of interfaces $\Gamma_f^{k,k+1}$ within one cell, so in $\ol{Y_f}$. 

\begin{figure}[ht]
	\centering
  \includegraphics[width=0.95\linewidth]{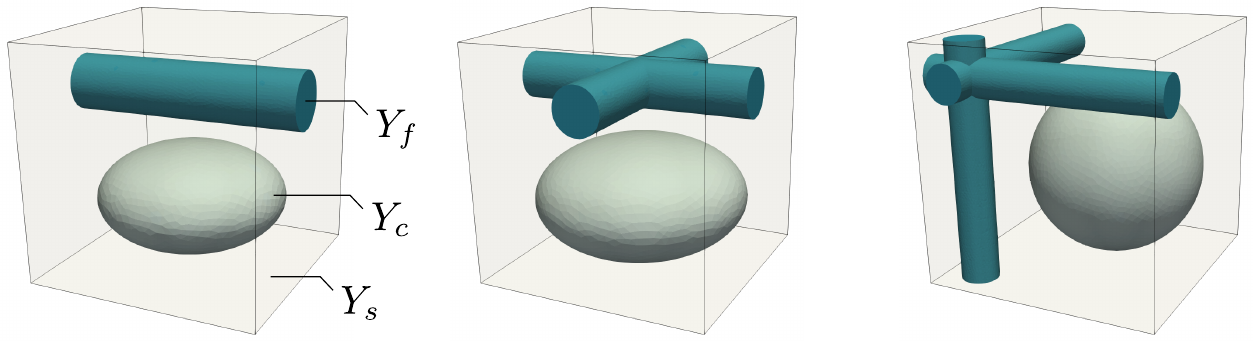}
  \caption{The 3D unit periodic cells permitting flow in one, two, and three directions.}
  \label{fig-Y-3D}
\end{figure}

%
Figure~\ref{fig-Y-per} (left) illustrates a 2D microstructure with reduced connectivity, permitting the fluid flow in only one direction. However, this representation is purely
illustrative as the connecting channels can be arbitrary in one or more directions as
shown in Fig.\ref{fig-Y-3D} for 3D microstructures.

\paragraph{Scaling in the fluid flow}
The parameters describing the interface permeabilities, $\kappa_A^\veps,\kappa_E^\veps$, and $\kappa_f^\veps$ require a specific scaling in the context of the asymptotic analysis. By virtue of the valve property and the ``no flow'' in the inclusions, it is adequate to consider a progressively small  $\kappa_{A,E}^\veps$ with $\veps\rightarrow 0$. On contrary, $\kappa_f^\veps$ is intended to reduce the flow, but not to separate the channel compartments $\Om_f^{k,\veps}$ on $\Gamma_f^\veps$. Therefore, we consider 
\begin{equation}\label{eq-kappa-eps}
\begin{split}
  \Tuf{\kappa_V^\veps(x)} & = \veps \bar\kappa_V(x)\;, \quad V = A,E\;,\\
  \Tuf{\kappa_f^\veps(x)} & = \chE{\frac{1}{\veps}}\hat\kappa_f(x)\;.
\end{split}
\end{equation}
\chE{Recalling the periodic unfolding operator $\Tuf{}$, the valve permeabilities $\bar\kappa$ can only vary with the macrostructure.}
The fluid viscosity scaling, $\mu^\veps = \veps^2 \bar\mu$, is conventional to respect the no-slip BC on pore walls $\Gamma_\fsi^\veps$. 

\subsection{Asymptotic analysis -- limit model}\label{sec-asymp}
Truncated asymptotic expansions can be justified by the a~priori estimates of solutions to the problem \eq{eq-wf1}-\eq{eq-wf3}; due to the topology with two disconnected compartments, two macroscopic pressure fields, $p_f^0$ and $p_c^0$ are involved,
\begin{equation}\label{eq-FS14a}
\begin{split}
\Tuf{\ub^{\veps}} \approx \ub^{R\veps}(x,y) & := \ub^{0}(x) + \veps \ub^{1}(x,y)\;,  \\
\Tuf{p_f^\veps}\approx p_f^{R\veps}(x,y) & := p_f^{0}(x) + \veps p_f^{1}(x,y) \;,\\
\Tuf{p_c^\veps}\approx p_c^{R\veps}(x,y) & := p_c^{0}(x)\;,\\
\Tuf{\wb^\veps} \approx \wb^{R\veps}(x,y) & := \hat\wb(x,y)\;.
\end{split}
\end{equation}
The test functions 
are given accordingly,
\begin{equation}\label{eq-FS14c}
\begin{split}
  \Tuf{\vb^\veps(x)} & = \vb^0(x) + \veps\vb^1(x,y)\;,\\
  \Tuf{q_f^\veps}\approx p_f^{R\veps}(x,y) & := q_f^{0}(x) + \veps q_f^{1}(x,y) \;,\\
  \Tuf{q_c^\veps}\approx p_c^{R\veps}(x,y) & := q_c^{0}(x)\;,\\
    \Tuf{\thetabf^\veps(x)} & = \hat\thetabf(x,y)\;.
\end{split}
\end{equation}

All the two-scale functions are $Y$-periodic in the second variable $y$ and for almost all $t>0$
\begin{equation}\label{eq-FS14d}
\begin{split}
  \ub^1(\cdot,t), \vb^1 & \in L^2(\Om; \Hpdb(Y_s))\;,\quad \ub^0(\cdot,t), \vb^0 \in \Hdb(\Om)\;,\\
\hat\wb(\cdot,t),\hat\thetabf & \in L^2(\Om;\HpdbO(Y_f))\;,\\
p_f^1(\cdot,t),q_f^1 & \in L^2(\Om; H_\#^1(Y_f\setminus\Gamma_f))\;,\quad p_f^0(\cdot,t), q_f^0 \in H^1(\Om)\;, \\
  \quad p_c^0(\cdot,t), q_c^0 & \in H^1(\Om)\;, 
\end{split}
\end{equation}
where 
\begin{equation}\label{eq-Hsp}
  \begin{split}
    \HpdbO(Y_f) & = \{\vb \in \Hpdbav(Y_f)|\; \vb = \bmi{0}  \mbox{ on } \Gamma_{fs}\}\;,\\
    H_\#^1(Y_f\setminus\Gamma_f)) & = \{q \in L^2(Y_f)|\; q \chi_k \in H_\#^1(Y_f^k)\;\mbox{ for each } Y_f^k \subset Y_f\}\;,
\end{split}
\end{equation}
recalling the channel decomposition.

In the next sections, convergence of all terms dealing with interfaces, as involved in the micromodel, is derived using the asymptotic expansions substituted into the unfolded integrals. Then, in Section~\ref{sec-loc}, the local problems providing the characteristic responses are introduced. The macroscopic equations are presented in Section~\ref{sec-macro-lin}, where also expressions for the involved effective model parameters -- the homogenized coefficients, are given.

\subsubsection{Limit of the interface terms}
The micromodel equations involve expressions related to the pore walls $\Gamma_\fsi^\veps$, the valves (connecting to surfaces $\pd\Om_c^{k,\veps}$ and $\Gamma_\fsi^\veps$), and the channel interfaces $\Gamma_f^\veps$.
By virtue of \eq{eq-kappa-eps} and the expansions \eq{eq-FS14a}, \eq{eq-FS14c}, we get for $p_f^\veps,q_f^\veps \in \Qcal_f^\veps(\Om_f^\veps)$
\begin{equation}\label{eq-va_ff}
\begin{split}
  \int_{\Gamma_f^\veps}\kappa_f^\veps\jump{p_f^\veps} \jump{q_f^\veps} = \int_\Om \frac{1}{\veps}\int_{\Gamma_f}  \hat\kappa_f\jump{\Tuf{p_f^\veps}} \jump{\Tuf{q_f^\veps}}\\ \rightarrow \int_\Om \int_{\Gamma_f}\hat\kappa_f\jumpG{p_f^1}\jumpG{q_f^1}\;,
\end{split}
\end{equation}
with $p_f^1,q_f^1 \in H_\#^1(Y_f\setminus\Gamma_f))$,
This will be employed in the limit of \eq{eq-wf2}.
There the terms related to the valves converge, as follows
\begin{equation}\label{eq-va_fcfA}
\begin{split}
  & \int_{\Gamma_f^\veps}\kappa_A^\veps\posPart{p_f^\veps- p_c^\veps}q_f^\veps \dS  = \int_\Om \frac{1}{\veps}\int_{\Gamma_\fsi} \veps\bar\kappa_A\posPart{\Tuf{\hat p_f^\veps}- \Tuf{p_c^\veps}}\Tuf{q_f^\veps} \hat\dlt_f \dSy\\
  & \rightarrow \int_\Om  \bar\kappa_A\posPart{p_f^0 - p_c^0}q_f^0\;,
\end{split}
\end{equation}
and, in analogy, for the ejection valves involving $\kappa_E^\veps$. To treat \eq{eq-wf3}, we need
\begin{equation}\label{eq-va_fcfB}
\begin{split}
  & \int_{\pd\Om_c^\veps}  \kappa_E^\veps\bar\dlt_E^\veps \posPart{p_c^\veps - \hat p_f^\veps - \Dlt P_E}q_c^\veps\dS\\
  & = \int_\Om \frac{1}{\veps}\intY_{\pd Y_c}\veps\bar\kappa_E\bar\dlt_E\posPart{\Tuf{p_c^\veps}- \Tuf{\hat p_f^\veps}-\Dlt P}\Tuf{q_c^\veps} \bar\dlt_E \dSy\\
  & \rightarrow  \int_\Om  \bar \kappa_E\posPart{p_c^0 - p_f^0- \Dlt P}q_c^0\;,
\end{split}
\end{equation}
and, in analogy, for the admission valves involving $\kappa_A^\veps$.

Finally we consider limit of the interface virtual work $\Scal_\sigma^\veps$ associated with the channels flow. Although the loading of the skeleton is modified due to the pressure discontinuity, formally the limit expressions is identic with the one obtained in \cite{Rohan-Lukes-CaS2024}, \cf \cite{Rohan-zamp2020}, for a continuous pressure.
In general, each compartment $\Om_f^{k,\veps}$ is the image of $Y_f\setminus\Gamma_f$, hence the Stokes theorem applied to surface $\pd\Om_f^{k,\veps}$ yields the divergence in $Y_f\setminus\Gamma_f$

 The virtual work of $\sigmabf_f^\veps\cdot\nb^\sx$ traction on $\Gamma_\fsi^\veps$ 
  \begin{equation}\label{eq-FS7c}
  \begin{split}
    &\Scal_\sigma^\veps(\sigmabf_f^\veps,\vb^\veps) = \sum_{k=1}^N\int_{\pd \Om_f^{k,\veps}} \nb^\sx\cdot\sigmabf_f^\veps\cdot\vb^\veps - \int_{\pd_\ext \Om_f^\veps} \nb^\sx\cdot\sigmabf_f^\veps\cdot\vb^\veps\dS \\
    &=  - \int_{\Om_f^\veps} \nabla\cdot (p^\veps \vb^\veps) +\int_{\Om_f^\veps}\nabla\cdot( {\veps^2}\bar\Dop_f \eeb{\wb^\veps}\cdot\vb^\veps) - \int_{\pd_\ext \Om_f^\veps} \nb^\sx\cdot\sigmabf_f^\veps\cdot\vb^\veps\dS \\
    & \rightarrow \int_\Om  \phi_f\nabla_x\cdot(p^0\vb^0)
    +  \int_\Om\vb^0\cdot\intY_{Y_f\setminus\Gamma_f}  \nabla_y p_f^1 - \int_\Om p_f^0 \intY_{\Gamma_\fsi} \vb^1\cdot \nb^\sx \\
    & - {\int_\Om  \vb^0\cdot\intY_{Y_f\setminus\Gamma_f}\nabla_y\cdot\bar\Dop_f\eeby{ \hat\wb}}  + \int_{\pd_\sigma \Om} \phi_f  \wtilde\bb^f \cdot\vb^0\;,
\end{split}
  \end{equation}
  noting that the surface integral of $\vb^1\cdot \nb^\sx$ on $\pd Y_f\setminus\Gamma_f$ is equivalent to the one on $\Gamma_\fsi$ due to the continuity of $\vb^1$ on $\Gamma_f$.
  Using the local momentum conservation for the quasistatic flow, upon integration in $Y_f\setminus\Gamma_f$
\begin{equation}\label{eq-FS19}
  \begin{split}
     {\intY_{Y_f\setminus\Gamma_f}\nabla_y\cdot\bar\Dop_f\eeby{ \hat\wb} -\intY_{Y_f\setminus\Gamma_f}\nabla_y p_f^1} 
    & = - (\fbav^f-\nabla_x p_f^0)\phi_f \;,
\end{split}
\end{equation}
\chE{in \eq{eq-FS7c},} we can eliminate the viscous part of the stress which is not involved in the limit expression,
  \begin{equation}\label{eq-FS20}
  \begin{split}
  &  \tilde\Scal_\sigma^\veps(\tilde\sigmabf_f^\veps,\vb^\veps)  \rightarrow \\
    & \int_\Om {\phi_f p_f^0\nabla\cdot\vb^0-  \int_\Om p_f^0 \intY_{\Gamma_\fsi} \vb^1\cdot \nb^\sx}
    + \int_\Om \phi_f\fbav^f\cdot\vb^0 +  \int_{\pd_\sigma \Om} \bar\phi_f  \wtilde\bb^f \cdot\vb^0\;,
\end{split}
  \end{equation}
  where $\bar\phi_f$ is the boundary porosity. Note that on $\Gamma_p^\veps$, the
prescribed pressure $p_\pd$ represents the boundary traction stress acting on
the fluid part. 
The condition prescribing the
averaged traction stress $\bar\phi_f\bar\bb^f = -\nb p_\pd$ on $\pd_p\Om$
emerges in \eq{eq-FS20} due to the convergence of $\bb^{f,\veps}$ acting on
$\Gamma_p^\veps$.

  The traction due to the inclusion pressure acting on the surface $\pd \Om_c^\veps$ yields
  \begin{equation}\label{eq-FS20c}
  \begin{split}
  &   \int_{\pd \Om_c^\veps} p_c^\veps\nb^\sx\cdot\vb^\veps \dS \rightarrow -\int_\Om \phi_c p_c^0\nabla\cdot\vb^0 +  \int_\Om p_c^0 \intY_{\pd Y_c} \vb^1\cdot \nb^\sx\dSy\;.
\end{split}
  \end{equation}


\subsubsection{Limit of the mass conservation}
Using \eq{eq-va_ff} and \eq{eq-va_fcfA}, we pass to the limit in \eq{eq-wf2} to obtain the two-scale mass conservation in the channels,
\begin{equation}\label{eq-lim-Mf}
  \begin{split}
    & \int_\Om \left( \nabla_x\cdot\left(\phi_f\dot\ub^0 +  \intY_{Y_f}\hat\wb\dVy\right) + \intY_{Y_f}\nabla_y\cdot \dot{\wtilde{\ub}}^1\dVy
 + \gamma \dot p_f^0\right) q_f^0 \dVx \\
& 
 - \int_\Om  \bar\kappa_A\posPart{p_f^0 - p_c^0}q_f^0 + \int_\Om  \bar \kappa_E\posPart{p_c^0 - p_f^0- \Dlt P}q_f^0 = 0\;,\\
 & \int_\Om \left(\int_{\Gamma_f}\hat\kappa_f\jumpG{p_f^1}\jumpG{q_f^1}-\intY_{Y_f}\hat\wb\cdot\nabla_yq_f^1 \dVy\right)\dVx = 0\;,
 \end{split}
\end{equation}
for all  $q^0 \in L^1(\Om)$ and $q_f^1 \in L^2(\Om;H_\#^1(Y_f\setminus\Gamma_f))$.
Then we consider the mass conservation in the inclusions.
Recalling piecewise constant $q^\veps \in \Qcal_c^\veps(\Om_c^\veps)$ and using \eq{eq-va_fcfB}, in the limit, \eq{eq-wf3} yields,
\begin{equation}\label{eq-lim-Mc}
  \begin{split}
    0 = 
    \gamma\int_\Om \phi_c \dot p^0 q_c^0 +   \int_{\Om} q_c^0 \left( \phi_c \nabla_x\cdot \dot\ub^0- \intY_{Y_s} \nabla_y\cdot \dot\ub^1\right)\\
   + \int_\Om  \bar\kappa_A\posPart{p_f^0 - p_c^0}q_c^0 -\int_\Om  \bar \kappa_E\posPart{p_c^0 - p_f^0- \Dlt P}q_c^0\;,
\end{split}
\end{equation}
for all $q_c^0 \in L^2(\Om)$.

\subsection{Local problems in the reference cell}\label{sec-loc}
Using the particular convergence results obtained above, the homogenization is applied in analogy with the treatment of the peristaltic induced flow \cite{Rohan-Lukes-PZ-porel}. For the sake of brevity, we shall employ the bilinear form associated with the elastic energy and the viscous dissipation power,
\begin{equation}\label{eq-bilin}
\begin{split}
\aYs{\ub}{\vb}  = \intY_{Y_s}\eeby{\vb} : \Aop\eeby{\ub} \dVy\;, \quad \ub,\vb \in \Hpdb(Y_{s})\;,\\
\aYf{\hat\wb}{\hat\thetabf} = \intY_{Y_f}\eeby{\hat\thetabf} : \Dop_f\eeby{\hat\wb} \dVy
\chE{+\intY_{\Gamma_f} (\hat\kappa_f)^{-1} \nb\otimes\nb : \hat\thetabf \otimes\hat\wb}
\;\;, \quad \hat\wb,\hat\thetabf \in \HpdbO(Y_f)\;,
\end{split}
\end{equation}
\chE{where $\nb$ is the unit normal vector on the semiperable membrane $\Gamma_f$.}

\subsubsection{Local equilibrium in the solid and fluid phases} 
The equilibrium equations \eq{eq-wf1} related to both solid and fluid attain quite analogous form to the one reported in \cite{Rohan-Lukes-CaS2024} from where one can simply deduce the resulting two-scale limit equations. We pursue the homogenization procedure explained in \cite{Rohan-Lukes-PZ-porel} and  report separately the local and the global subproblems.

The local one in the solid is obtained 
using the expansions \eq{eq-FS14a} and \eq{eq-FS14c} for vanishing test displacement $\vb^0 = 0$ in \eq{eq-FS14c}, such that, due to the modified interaction term \eq{eq-FS20}-\eq{eq-FS20c} in the limit, \eq{eq-wf1}$_1$ yields
\begin{equation}\label{eq-L*S1}
\begin{split}
& \intY_{\Om \times Y_s}
\eeby{\vb^1} : \Aop \GrxyS{\ub^0}{\ub^1}\dVxy\\
& =  -\int_\Om p_f^0 \intY_{\Gamma_\fsi} \vb^1 \cdot \nb^\sx \dSy \dVx
- \int_\Om p_c^0 \intY_{\pd Y_c} \vb^1\cdot \nb^\sx\dSy \dVx\;,
\end{split}
\end{equation}
for all $\vb^1 \in \Lb^2(\Om;\Hpdb(Y_s))$.

In the fluid, since $\hat\wb$ is a two-scale function, see \eq{eq-FS20}, the limit Stokes flow problem obtained from \eq{eq-wf1}$_2$ is local, being constraint by the incompressibility condition \eq{eq-lim-Mf}$_2$,
\begin{equation}\label{eq-L*F1}
\begin{split}
 \intY_{\Om} \left( \aYf{\hat\wb}{\hat\thetabf} 
 - \intY_{Y_f}p_f^1\nabla_y\cdot\hat\thetabf \dVy\right)\dVx 
 &  = \int_\Om(\fbav^f -  \nabla_x p_f^0 )\dVx \cdot\intY_{Y_f}\hat\thetabf\dVy\;,\\
\int_{\Gamma_f}\hat\kappa_f\jumpG{p_f^1}\jumpG{q_f^1}\dSy-\intY_{Y_f}\hat\wb\cdot\nabla_yq_f^1 \dVy & = 0\;,
 \end{split}
\end{equation}
for all $\hat\thetabf \in L_0^\infty(\Om;\HpdbO(Y_f))$, $q_f^1 \in  L^2(\Om;H_\#^1(Y_f\setminus\Gamma_f))$, where $L_0^\infty(\Om;\HpdbO(Y_f))$ is the space of two-scale functions with compact support in $\Om$ (vanishing on $\pd\Om$).

\subsubsection{Characteristic response in the fluid channel $Y_f$}
In the fluid compartment, the  characteristic responses $(\wb^k, \pi^k)$, $k=1,2,3$ can be defined,
such that (the summation in $k$ applies),
\begin{equation}\label{eq-A10a}
\begin{split}
  \hat\wb & = \psibf^k(\what f_k - \pd_k^x p_f^0)\;,\quad
 p_f^1  = \pi^k(\what f_k - \pd_k^x p_f^0)\;,
\end{split}
\end{equation}
where $(\psibf^{k},\pi^k)$ are solutions of the following problem: Find $\psibf^{k}\in\HpdbO(Y_f)$ and \chE{$\pi^{k}\in H_\#^1(Y_f\setminus\Gamma_f)/\RR$}, such that
\begin{equation}\label{eq-S3}
  \begin{split}
    \aYf{{\psibf}^k}{\vb} +\ipYf{{\nabla_y \pi}^k}{\cdot\vb}  & = \ipYf{\onebm_k}{\vb}\;, \\ 
    -\ipYf{\nabla_y q}{\cdot{\psibf}^k}  + \ipGG{\hat\kappa_f\jumpG{\pi^k}}{\jumpG{q}} & = 0\;, 
\end{split}
\end{equation}
for all $\vb \in \HpdbO(Y_f)$ and \chE{$q \in H_\#^1(Y_f\setminus\Gamma_f)$}.

The mean flow velocity is expressed through the permeability $\Kb = (K_{ij})$ which is defined by
\begin{equation}\label{eq-S7a}
\begin{split}
  K_{ij} = \intY_{Y_f}\psi_j^i =\aYf{{\psibf}^j}{{\psibf}^i} + \ipGG{\hat\kappa_f\jumpG{\pi^j}}{\jumpG{\pi^i}}   \;.
\end{split}
\end{equation}
\chE{Note that $\hat\kappa_f$ is involved also in the first term due to \eq{eq-bilin}.}

\subsubsection{Characteristic response in the solid skeleton $Y_s$}
By virtue of \eq{eq-L*S1}, the fluctuating (two-scale) displacements $\ub^1$ can be decomposed using the characteristic autonomous responses -- the correctors $\omegabf^{ij},\omegabf^f,\omegabf^c$
\begin{equation}\label{eq-L2b}
\begin{split}
  \ub^1(x,y) & = \omegabf^{ij}e_{ij}^x(\ub^0)  - p_f^0 \omegabf^f - p_c^0 \omegabf^c\;.\\
\end{split}
\end{equation}
Hence, upon substituting in \eq{eq-L*S1}, the following  characteristic problems are identified:
\begin{itemize}
\item \emph{(strain correctors)} Find
$\omegabf^{ij}\in \Hpdb(Y_{s})$ for any $i,j = 1,2,3$
satisfying
\begin{equation}\label{eq-mip1}
\begin{split}
\aYs{\omegabf^{ij} + \Pibf^{ij}}{\vb}  = 0\;, \quad \forall \vb \in  \Hpdb(Y_{s})\;,\\
\end{split}
\end{equation}
where  $\Pibf^{ij} = (\Pi_k^{ij})$, $i,j,k   = 1,2,3$ with components $\Pi_k^{ij} = y_j\delta_{ik}$ represents the homogeneous displacement field due to the principal macroscopic strain modes.
\item  \emph{(pressure correctors)} Find
$\omegabf^P\in \Hpdb(Y_{s})$, for the two compartments $P = f,c$,
satisfying
\begin{equation}\label{eq-mip2}
\begin{split}
\aYs{\omegabf^c}{\vb}  = -\intY_{\pd Y_c} \vb\cdot \nb^\sx \dSy\;, \quad
\forall \vb \in  \Hpdb(Y_s) \;,\\
\aYs{\omegabf^f}{\vb}  = -\intY_{\Gamma_\fsi} \vb\cdot \nb^\sx \dSy\;, \quad
\forall \vb \in  \Hpdb(Y_s) \;.
\end{split}
\end{equation}
\end{itemize}


\subsection{Macroscopic problem in fixed reference configuration}\label{sec-macro-lin}
Although the inflation effect is naturally associated with the geometrical nonlinearity due to the inflation pressure acting on deformed surfaces, we first introduce the model which results from the treatment of the problem in a fixed reference configuration. Later, in Section~\ref{sec-nonlin-II}, we provide an extension of the model to respect these additional nonlinear effects.

The limit macroscopic equations are obtained from the two-scale limit equations
with the only non-vanishing macroscopic test functions $\vb^0,q_f^0$ and $q_c^0$.
The macroscopic state of the poroelastic medium is described by $(\ub^0,p_f^0,p_c^0)$ satisfying the quasistatic equilibrium of the homogenized fluid-saturated solid skeleton, and the mass conservation in the two fluid compartments: in the channels and in the inclusions. 

We employ the admissibility sets for  $\ub^0 \in \Ucalbf(\Om)$ and $p_f^0\in \Qcal(\Om)$ defined to respect the BCs on $\pd\Om$, see \eq{eq-BC1},
\begin{equation}\label{eq-UPspaces}
\begin{split}
\Ucalbf(\Om) = \{\ub \in \Hdb(\Om)|\; \ub = \ub_\pd \mbox{ on }\pd_u\Om\}\;,\\
\Qcal(\Om) = \{ p \in H^1(\Om)|\; p = p_\pd \mbox{ on }\pd_p\Om\}\;.
\end{split}
\end{equation}
Accordingly, the test function spaces $\Ucalbf_0(\Om)$ and $\Qcal_0(\Om)$, associated with $\Ucalbf(\Om)$ and $\Qcal(\Om)$, respectively, are defined by virtue of prescribing the zero Dirichlet boundary conditions satisfied by the test functions.

The limit macroscopic equilibrium is obtained from \eq{eq-wf1}$_1$; upon substituting there the truncated expansions \eq{eq-FS14a} and \eq{eq-FS14c} with vanishing test displacement $\vb^1 = 0$, and passing to the limit $\veps\rightarrow 0$. Using the interface integral convergence \eq{eq-FS20} and \eq{eq-FS20c}, we get  
\begin{equation}\label{eq-H1*S1}
\begin{split}
\int_{\Om} &
\aYs{\ub^1-\Pibf^{kl}e_{kl}^x(\ub^0)}{\Pibf^{ij}} e_{ij}^x(\vb^0) \dVx \\
& -\int_{\Om} (p_f^0 \phi_f + p_c^0 \phi_c)\nabla_x\cdot \vb^0 \, \dVx  =
\int_{{\Om}}\fbav \cdot \vb^0 \, \dVx + \int_{\pd {\Om}}\bbav \cdot \vb^0\dSx\;,
\end{split}
\end{equation}
for all $\vb^0 \in \Ucalbf_0(\Om)$. Then the fluid mass conservation in the channels \eq{eq-lim-Mf}$_1$ and the one in the inclusions \eq{eq-lim-Mc} is involved, with the only non-vanishing test functions $q_c^0$ and $q_f^0$. This completes the system of macroscopic equations in which the homogenized coefficients can be identified due to the characteristic responses treated in Section~\ref{sec-loc}.

\subsubsection{Homogenized coefficients}
Besides the permeability introduced in \eq{eq-S7a} to express the average  fluid velocity relative to the solid skeleton (the Darcy law)
\begin{equation}\label{eq-darcy}
\begin{split}
\wb^0 = \intY_{Y_f} \hat\wb \dVy = \intY_{Y_f}\psibf^i  \dVy (\what{f}_i -\pd_i^x p_f^0) = \Kb  (\fbav^f-\nabla_x p_f^0)\;,
\end{split}
\end{equation}
the macroscopic equations involve further homogenized coefficients of the generalized poroelasticity.
These are identified in \eq{eq-lim-Mf}$_1$,\eq{eq-lim-Mc} and \eq{eq-H1*S1} upon substituting there the decomposition \eq{eq-L2b}. Using the micro-problems \eq{eq-mip1}-\eq{eq-mip2}, we get
\begin{equation}\label{eq-HC1}
\begin{split}
A_{klij} & = \aYs{\omegabf^{ij} + \Pibf^{ij}}{\omegabf^{kl}+\Pibf^{kl}}\;, \\
B_{ij}^P & = \aYs{\omegabf^P}{\Pibf^{ij}} + \phi_P\delta_{ij} = \intY_{\Gamma_P} \omegabf^{ij} \cdot \nb^\px \dSy + \phi_P\delta_{ij} \\
& = \intY_{Y_P}\nabla_y\cdot \wtilde\omegabf^{ij} + \phi_P\delta_{ij}\;,\\
M^{PQ} & = -\intY_{\Gamma_Q} \omegabf^P\cdot \nb^\sx \dSy + \phi_P\gamma \delta_{PQ}
= \aYs{\omegabf^P}{\omegabf^Q} 
+ \phi_f\gamma \delta_{PQ}\;,
\end{split}
\end{equation}
where $P,Q  = f,c$ refers to the fluid compartments, such that $\Gamma_c \equiv \pd Y_c$ and  $\Gamma_f \equiv \Gamma_\fsi$. Note that $\delta_{ff} = \delta_{cc} = 1$, whereas $\delta_{fc} = \delta_{cf} = 0$.
Clearly, $M^{cf} = M^{fc}$ and $M^{cc} > 0$ and $M^{ff} > 0$.

\subsubsection{Macroscopic problem formulation}
The weak formulation is obtained from the three macroscopic equations, \eq{eq-lim-Mf}$_1$,\eq{eq-lim-Mc} and \eq{eq-H1*S1}, upon substituting there the homogenized coefficients. For simplicity, we may consider ``zero'' initial conditions. 
For any time $t > 0$, find  $(\ub^0(t,\cdot),p_f^0(t,\cdot),p_c^0(t,\cdot)) \in \Ucalbf(\Om)\times \Qcal(\Om) \times L^2(\Om)$ which satisfies
\begin{equation}\label{eq-S25-fc}
\begin{split}
&\int_{\Om} \left (\Cop  \eeb{\ub^0}
-  p_f^0  \Bb^f -  p_c^0  \Bb^c \right ):\eeb{\vb^0}  = \int_{\Om} \fbav \cdot \vb^0
+ \int_{\pd \Om}\bbav \cdot \vb^0\dSx\;, \quad \forall \, \vb^0 \in
\Ucalbf_0(\Om) \;,\\
&\int_{\Om} q_f^0\left (\Bb^f :\eeb{\dot\ub^0} + M^{ff} \dot p_f^0 + M^{fc} \dot p_c^0 \right) + \int_{\Om} \Kb\left(\nabla_x p_f^0 -\what{\fb}^f\right) \cdot \nabla_x q_f^0 \\
&+  \int_{\Om}q_f^0\left ( \bar\kappa_A\posPart{p_f^0 - p_c^0} - \bar\kappa_E\posPart{p_c^0 - p_f^0- \Dlt P}\right)
= 0\;,\quad \forall q_f^0 \in \Qcal_0(\Om)\;,\\
&\int_{\Om} q_c^0\left (\Bb^c :\eeb{\dot\ub^0} + M^{cf} \dot p_f^0 + M^{cc} \dot p_c^0 \right) \\
&-  \int_{\Om}q_c^0\left ( \bar\kappa_A\posPart{p_f^0 - p_c^0} - \bar\kappa_E\posPart{p_c^0 - p_f^0- \Dlt P}\right)
= 0\;,\quad \forall q_c^0 \in L^2(\Om)\;,\\
\end{split}
\end{equation}
Remarks: The boundary conditions for the flow can be considered, as follows:
$p_f^0 = \bar p_\inlet$ and $p_f^0 = \bar p_\outlet$ on the inlet $\pd_\inlet \Om$ and $\pd_\outlet \Om$, respectively. The solid can be clamped at $\pd_\inlet \Om$, for instance.

The macroscopic problem \eq{eq-S25-fc} is nonlinear due to the cut-off functions associated with valves opening and closing. Also the ``valve permeabilities'' $\hat\kappa_f$ and $\bar\kappa_A$, $\bar\kappa_E$ can be controlled by any of the macroscopic field variables, such as by pressure $p_c^0$.

Moreover, we aim to capture also the nonlinearity arising due to the equilibrium satisfied in the deformed configuration. It is therefore necessary to consider an incremental formulation in the sense of the Newton-Raphson iterations, as described in \cite{Rohan-Lukes-CaS2024}.


\subsection{Nonlinear setting -- incremental formulation}\label{sec-increment}
\chE{In the rest of the paper we drop the superscript ${}^0$ labeling the macroscopic variables, so $\ub \equiv \ub^0$ and $p_a \equiv p_a^0$, $a = f,c$.}
Due to the two nonlinearities mentioned above, the macroscopic model \eq{eq-S25-fc} governing the state variable $\ssb = (\ub,\ul{p})$, where $\ul{p} = (p_f,p_c)$ can be interpreted in the context of
the perturbation $\dlt\ssb = (\dlt\ub,\dlt\ul{p})$ and the current solution approximation $\dlt\ssb$,
such that $\ssb = \bar\ssb + \dlt\ssb$ and $\dlt\ul{p} = (\dlt p_f,\dlt p_c)$. Since the problem is time dependent, by virtue of the incremental problem, it is convenient to consider a time-discretized reformulation related to time levels $t_k = k \Delta t$, $k=0,1,2,\dots$.

We consider time interval 
$[t_{k-1}, t_k]$ and use the abbreviated notation $t:=t_k$, $t\prevstep := t_{k-1}$.
By $\ssb$
we shall refer to the unknown fields at time $t = t_k \in\{t_i\}_i$,
whereas $\ssb^\tp \approx \ssb(t^\tp,\cdot)$ 
thus being associated with time $t - \Delta t$.
Using the implicit approximation $\dot\ssb \approx (\ssb - \ub^\tp)/\Delta t$.
problem  \eq{eq-S25-fc} can be discretized in time. 
The nonlinear problem is solved using a Newton-type iterative method, therefore,
we introduce the residual-based incremental formulation.

\subsubsection{Incremental problem I: valve non-linearity only}\label{sec-nonlin-I}
For the sake of clarity, we shall first consider the nonlinearity only due to the valves. Using the time discretization in \eq{eq-S25-fc}, the residual functional is given by ($\qb = (\vb,\ul{q})$ denotes the test functions)
\begin{equation}\label{eq-S25-psi}
\begin{split}
& \Psi^t(\ssb;\qb) =   
\int_{\Om} \left (\Cop  \eeb{\ub}
-  p_f \Bb^f -  p_c \Bb^c \right ):\eeb{\vb}  - \int_{\Om} \fbav \cdot \vb
- \int_{\pd \Om}\bbav \cdot \vb\dSx   \\
& + \int_{\Om} q_f\left (\Bb^f :(\eeb{\ub} - \eeb{\ub^\tp}) + \ul{M}^f (\ul{p} - \ul{p^\tp})\right) + \Dlt t \int_{\Om} \Kb\left(\nabla_x p_f -\what{\fb}^f\right) \cdot \nabla_x q_f \\
&+  \int_{\Om}q_f\left ( \bar \kappa_A\posPart{p_f - p_c} - \bar \kappa_E\posPart{p_c - p_f- \Dlt P}\right) \\
& +\int_{\Om} q_c\left (\Bb^c :(\eeb{\ub} - \eeb{\ub^\tp}) + \ul{M}^c (\ul{p} - \ul{p^\tp})\right) \\
&-  \int_{\Om}q_c\left ( \bar \kappa_A\posPart{p_f - p_c} - \bar \kappa_E\posPart{p_c - p_f- \Dlt P}\right)\;,
\end{split}
\end{equation}
where $\ul{M}^P \ul{p} = M^{Pf}p_f + M^{Pc} p_c$ for $P = f,c$.

\paragraph{Abstract incremental problem}
We employ $\Scalbf(\Om) =  \Ucalbf(\Om)\times\Qcalbf(\Om)$ and $\Scalbf_0(\Om) =  \Ucalbf_0(\Om)\times\Qcalbf_0(\Om)$, where $\Qcalbf(\Om) = \Qcalbf_0(\Om) = \Qcal_0(\Om) \times L^2(\Om)$.
For a given previous time step solution $\ssb^\tp$, find $\ssb=(\ub,p) \in \Scalbf(\Om)$ 
which satisfies
 \begin{equation}\label{eq-S25TD}
   \Psi^t(\ssb;\qb) = 0\quad\mbox{  for any }\quad \qb \in \Scalbf_0(\Om)\;,
 \end{equation}
evaluated at time $t$ (hence the superscript in $\Psi^t$) and depending on the state $\ssb^\tp$ by virtue of the time discretization. 
It can be solved by the Newton-Raphson method which uses through successive iterations $\ssb^i$, $i = 0,1,2,\dots$ based on the
straightforward decomposition introduced in terms of
the recent state approximation $\ssb^i \equiv \bar\ssb \approx \ssb$ and the correction $\delta \ssb$, such that
\begin{equation}\label{eq-S50}
\begin{split}
\ssb^i \rightarrow \ssb\quad \mbox{ for } \ssb^{i+1} = \ssb^i + \dlt\ssb\;,
\end{split}
\end{equation}
whereby $\ssb^\tp:=\ssb$ for the next time level computing.
Within the $i$-th iteration, given $\bar\ssb \equiv \ssb^{i-1}$, the first order expansion in $\ssb$ leads to the obvious approximation of problem \eq{eq-S25TD} at $\bar\ssb$,
\begin{equation}\label{eq-rif2}
  \begin{split}
    0  & = \Psi^t(\ssb;\qb) \approx \Psi^t(\bar\ssb;\qb)  + \delta_\ssb  \Psi^t(\bar\ssb;\qb)\circ\dlt\ssb\;,
\end{split}
\end{equation}
to hold for any $\qb\in\Scalbf_0(\Om)$. In order to compute $\delta_\ssb  \Psi^t$ by differentiating \eq{eq-S25-psi} at $\ssb := \bar\ssb$ is straightforward, since the only nonlinear terms requiring a special care are those associated with cut-off function $[~]_+$. 
We define $\gamma_A(\ul{p})$ and  $\gamma_E(\ul{p})$, such that
\begin{equation}\label{eq-sts6i-V1}
\begin{split}
  \gamma_A(\ul{p}) =  \posPart{\sgn(p_f - p_c)} =  \left \{
  \begin{array}{ll}
    1 & \mbox{ when } p_f - p_c > 0\;,\\
    0 & \mbox{ otherwise }   
  \end{array} \right. \;,\\
  \gamma_E(\ul{p}) =  \posPart{\sgn(p_c - p_f - \Dlt P_E)} =  \left \{
  \begin{array}{ll}
    1 & \mbox{ when } p_c - p_f - \Dlt P_E > 0\;,\\
    0 & \mbox{ otherwise }   
  \end{array} \right. \;, 
\end{split}
\end{equation}
so that
\begin{equation}\label{eq-S25-psi-diff}
\begin{split}
& \dlt_\ssb\Psi^t(\ssb;\qb)\circ\dlt\ssb =   
\int_{\Om} \left (\Cop  \eeb{\dlt\ub}
- \dlt p_f \Bb^f - \dlt p_c \Bb^c \right ):\eeb{\vb}  + \Dlt t \int_{\Om} (\Kb\nabla_x \dlt p_f ) \cdot \nabla_x q_f \\
& \quad
+ \int_\Om q_f\big(\Bb^f:\eeb{\delta\ub}+  \ul{M}^f \dlt \ul{p}\big)
+\Delta t\int_\Om \nabla q_f\cdot\left ( \gamma_A(\ul{\bar p})\bar \kappa_A + \gamma_E(\ul{\bar p})\bar \kappa_E\right) (\dlt p_f - \dlt p_c)
\\
& \quad
+ \int_\Om q_c\big(\Bb^c:\eeb{\delta\ub}+  \ul{M}^c \dlt \ul{p}\big)
- \Delta t\int_\Om \nabla q_c\cdot\left ( \gamma_A(\ul{\bar p})\bar \kappa_A + \gamma_E(\ul{\bar p})\bar \kappa_E\right)(\dlt p_f - \dlt p_c)
\;.
\end{split}
\end{equation}
Note that $\gamma_A(\ul{p}) + \gamma_E(\ul{p}) \in \{0,1\}$ by virtue of the admission and ejection valve. To treat properly the non-smoothness at zero, \ie $\dlt [0]_+$ would yield the inclusion  $\gamma_{A,E}(\ul{p}) \in [0,1]$. A similar treatment to the one employed when dealing with ``complementarity problems'' leading to the non-smooth equation, see \eg \cite{Rohan-Heczko-CaS} where a nen-smooth version of the Newton method has been employed.

\subsubsection{Incremental problem II: equilibrium at deforming configuration}\label{sec-nonlin-II}
Problem \eq{eq-S25TD} is now considered in the context of additional geometric nonlinearity due to the the deforming microstructure which affects also the homogenized coefficients  $\Hop$  in the generic sense, involved in the macroscopic problem; in what follows, the ``tilde'' notation $\tilde \Hop$ means the the dependence of $\Hop$ on the state variables, namely on the local macroscopic strain $\eeb{\ub}$ and pressures $\ul{p}$. We pursue the approach proposed and validated in \cite{Rohan-Lukes-CaS2024}. Stress  $\sigmabf$ is given by the poro-hypoelastic constitutive law,
\begin{equation}\label{eq-E2b}
  \begin{split}
    \sigmabf(\ssb(t,x)) & = \int_0^t \dlt_{\ssb} \sigmabf({\ssb}(\tau,x))\circ\hat\ssb(x)\dd \tau\;,\quad \ssb(t,x) = \hat\ssb(x) t\;,\\
    \dlt_{\ssb} \sigmabf\circ\dlt\ssb(\tau,x) & = \tilde\Cop  \eeb{\dlt\ub}
- \dlt \ul{p}  \tilde{\ul{\Bb}}\;,
\end{split}
\end{equation}
where the homogeneity of the relationship $\dlt\ssb(\tau,x) = \hat\ssb(x) \dd \tau$ reflects the solid elasticity. In analogy with \eq{eq-S25-psi}
\begin{equation}\label{eq-rif1}
  \begin{split}
  & \Psi^t(\ssb;\qb)  =
\int_{\Om}\sigmabf(\ssb(t,x)):\eeb{\vb}\dVx  - \left( \int_{\Om} \wtilde\fb \cdot \vb \dVx
+ \int_{\pd_\sigma \Om}\wtilde\bb \cdot \vb\dSx \right)\\
& +\int_{\Om} \ul{q}\left (\tilde{\ul{\Bb}} :\eeb{\ub-\ub\prevstep} +  (\ul{p} - \ul{p}\prevstep) \tilde{\ul{M}} \right)\dVx \\
& + \int_{\Om} \frac{\Delta t}{\bar\mu}\tilde\Kb\left(\nabla_x p_f -\fbav^f\right) \cdot \nabla_x q^f \dVx\\
&+  \int_{\Om}q_f\left ( \tilde \kappa_A\posPart{p_f - p_c} - \tilde \kappa_E\posPart{p_c - p_f- \Dlt P}\right) \\
&-  \int_{\Om}q_c\left ( \tilde \kappa_A\posPart{p_f - p_c} - \tilde \kappa_E\posPart{p_c - p_f- \Dlt P}\right)\;,
\end{split}
\end{equation}
where the notation $(\ul{\tilde{M}} \ul{p})\ul{q} = \sum_{P = f,c}q^P (\tilde{M}^{Pf}p_f + \tilde{M}^{Pc} p_c)$, and in analogy for the expression involving $\ul{\tilde\Bb}$. The external forces $\wtilde\fb$ and $\wtilde\bb$ are obtained, as in ``the linear case'', corresponding to $\fbav$ and $\bbav$, respectively, however now can depend on the local micro-deformation (note they reflect the microstructure decomposition into $Y_s,Y_f$ and $Y_c$).


As proposed in \cite{Rohan-Lukes-nlBiot2015} and further elaborated in \cite{Rohan-Lukes-CaS2024}, all the perturbed coefficients $\tilde \Hop(\eeb{\ub},\ul{p})$ can be approximated using the first
order expansion formulae which results from the deforming cell due to microdeformation by the strain $\eebx{\ub^0} + \eeby{\ub^1}$. Since $\ub^1$ is given by \eq{eq-L2b}, the generic form of the Taylor expansion attains the form involving the ``total shape variation'' $\delta_\shp$,
\begin{equation}\label{eq-S24}
\begin{split}
\tilde \Hop(\eeb{\ub},\ul{p}) & \approx  \Hop^0 + \delta_{\eb} \Hop^0 : \eeb{\ub} + \sum_{P=f,c}\delta_{p_P} \Hop^0 p_P \;,\\
(\delta_{\eb}\Hop^0)_{ij} & := \big(\pd_{\eb} (\delta_\shp \Hop^0\circ{\ub^\mic})\big)_{ij} = \delta_\shp \Hop^0 \circ {(\omegabf^{ij} + \Pibf^{ij})}\;,\\
\delta_{p_P} \Hop^0 & := \pd_{p_P} (\delta_\shp \Hop^0\circ{\ub^\mic}) = \delta_\shp \Hop^0 \circ (-\omegabf^P)\;.
\end{split}
\end{equation}
Coefficients $\Hop^0$ are computed using \eq{eq-HC1} and \eq{eq-S7a} for the
reference ``initial'' configuration  and $\delta_\shp \Hop^0$ are
sensitivities which are derived for a general ``design velocity field'' in \cite{Rohan-Lukes-CaS2024}, Appendix C. The particular formulae employed for the homogenized coefficients involved in \eq{eq-rif1} are presented in \Appx{sec-SA}.

Iterative scheme \eq{eq-rif2} requires the total variation $\delta\Psi^t$ which can be expressed in terms of tangential incremental coefficients $\ol{\Hop}$ defined below, 
see Section~\ref{sec-tan}.
Now differentiation of the residual functional in \eq{eq-rif1} using \eq{eq-S24} yields,
\begin{equation}\label{eq-sts6i-V2}
  \begin{split}
    & \dlt_\ssb\Psi^t(\ssb;\qb)\circ\dlt\ssb 
    = \int_\Om \left(\ol{\Cop}\eeb{\delta\ub} - \sum_{P=f,c}\delta p_P \ol{\Bb}^P\right):\eeb{\vb}\\
    & \quad
      - \int_\Om \big(\pd_\eb{\fbav}\circ\eb(\delta\ub)+\sum_{P=f,c}\pd_{p_P}{\fbav}\circ\delta p_P\big)\cdot \vb
    -\int_{\pd\Om}\big(\pd_\eb{\bbav}\circ\eb(\delta\ub)+\sum_{P=f,c}\pd_{p_P}{\bbav}\circ\delta p_P\big)\cdot \vb \\
    & \quad
    +  \frac{\Delta t}{\bar\mu}\int_\Om \nabla q_f\cdot
    \left(\ol{\Kb}\nabla \delta p_f + \ol{\Gb}:\eeb{\delta\ub} + \sum_{P=f,c}\ol{\Qb}_{fP}\delta p_P
    \right)\\
    & \quad
    + \int_\Om q_f\big(\ol{\Db_f}:\eeb{\delta\ub}+  \sum_{P=f,c}\ol{M}_{fP}\delta p_P\big)
    +\Delta t\int_\Om \nabla q_f\cdot\left ( \gamma_A(\ul{\ppstep})\bar\kappa_A + \gamma_E(\ul{\ppstep})\bar\kappa_E\right) (\dlt p_f - \dlt p_c)
    \\
    & \quad
    + \int_\Om q_c\big(\ol{\Db_c}:\eeb{\delta\ub}+  \sum_{P=f,c}\ol{M}_{cP}\delta p_P\big)
    - \Delta t\int_\Om \nabla q_c\cdot\left ( \gamma_A(\ul{\ppstep})\bar\kappa_A + \gamma_E(\ul{\ppstep})\bar\kappa_E\right)(\dlt p_f - \dlt p_c)
    \;,
  \end{split}
\end{equation}
%
where all the expression for $\ol{\Hop}$ are given below. 
We recall that $\gamma_A(\ul{p})$ and  $\gamma_E(\ul{p})$ were defined in \eq{eq-sts6i-V1}.

\subsubsection{{Consistent homogenized coefficients for incremental formulation}}\label{sec-tan}
The expansion formula applied to each HC involved in the residual $\Psi^t$, see \eq{eq-rif1}, yields expressions of the following type, where $s^1,s^2$ represent any two ``state variables'', and $\tilde\Hop^j$, $j=1,2$ are 2 different coefficients \chE{(no summation in $j$, or other repeated indices unless the $\sum_l$ is used)},
\begin{equation}\label{eq-barH1}
  \begin{split}
    \tilde\Hop^j(s^1,s^2)s^j & \approx \left(\Hop^j + \sum_l\pd_l\Hop^j\circ s^l\right) s^j
  = \left(\Hop^j + \sum_l\pd_l\Hop^j\circ (\bar s^l + \dlt s^l)\right) (\bar s^j + \dlt s^j) \\
  & \approx (\Hop^j + \sum_l\pd_l\Hop^j\circ \bar s^l)\bar s^j
  + (\Hop^j + \sum_l\pd_l\Hop^j\circ \bar s^l)  \dlt s^j + \sum_l(\pd_l\Hop^j\circ\dlt s^l)\bar s^j \;,
  \end{split}
\end{equation}
where $\dlt s^i$ means the perturbation of $s^i$ and $\bar s^i$ is the reference value,
and $\pd_l\Hop^j\circ s^l$ means the differential of \chE{$\Hop^j$} due to the perturbation by $s^l$.
The higher order term $\sum_l\pd_l(\Hop^j\circ\dlt s^l)\dlt s^j$ was neglected. The differentiation $\pd_i$ \wrt the perturbation $\dlt s^i$  now yields \chE{($\dlt_{ij}$ is the Kronecker symbol, indeed)}
\begin{equation}\label{eq-barH2}
  \begin{split}
  \pd_i(\tilde\Hop^j s^j) = \pd_i \tilde\Hop^j(s^1,s^2)s^j & \approx \chE{\delta_{ij}(\Hop^i + \sum_l\pd_l\Hop^i\circ \bar s^l)} + (\dlt_i\Hop^j\circ(~))\bar s^j \;,\\
 \mbox{ hence }\quad   \dlt_i \tilde\Hop^js^j \circ \dlt s^i & = \chE{\delta_{ij}}(\Hop^i + \sum_l\pd_l\Hop^i\circ \bar s^l)\dlt s^i
    + (\dlt_i\Hop^j\circ\dlt s^i)\bar s^j\\
    & =  \chE{\delta_{ij}}(\Hop^i + \ol{\pd}\Hop^i)\dlt s^i + \dlt_i(\Hop^j\bar s^j)\circ\dlt s^i\;,\\
    \mbox{ where }\quad\ol{\pd}\Hop^i & =  \sum_l\pd_l\Hop^i\circ \bar s^l\;.
  \end{split}
\end{equation}
In the context of the actual HCs involved in  \eq{eq-rif1}, we can interpret $\ol{\pd}\Hop^i$ as the obvious Taylor expansion (in the 2nd expression, the generic form of $\Hop^0$ is applied to elasticity modulus $\Cop^0$), so that
\begin{equation}\label{eq-S40i}
  \begin{split}
    \ol{\pd}\Hop^0 &= \pd_\eb \Hop^0\circ\eeb{\bar\ub } + \sum_{P=f,c}\pd_{p_P} \Hop^0\circ \bar p_P \;,\\
        \ol{\pd}\Cop^0 &= \pd_\eb \Cop^0\circ\eeb{\bar\ub } + \sum_{P=f,c}\pd_{p_P} \Cop^0\circ \bar p_P \;,
\end{split}
\end{equation}
where the partial sensitivities are given in \eq{eq-S24}.
Implementation of \eq{eq-barH2} then leads to \eq{eq-sts6i-V2}, where the generic $\ol{\Hop}$ tensors involved,
\begin{equation}\label{eq-S41i}
  \begin{split}
    \ol{\Cop}& = \Cop^0 + \ol{\pd}\Cop^0 + \pd_\eb (\Cop^0\eeb{\bar\ub })\circ \circpare{~} \\
    & \quad  - \sum_{P=f,c} \pd_\eb (\Bb^P \bar p_P)\circ \circpare{~} 
    \;,\\
  \\ 
   \ol{\Bb}^Q& = \Bb^Q + \ol{\pd}\Bb^Q
   +\sum_{P=f,c} \pd_{p_Q} (\Bb^Q\bar p_P)\circ \circparpP{~}{Q} \\
    & \quad -\pd_{p_Q} (\Cop^0\eeb{\bar\ub })\circ \circparpP{~}{Q} 
   \;,
  \\
   \ol{\Db}_Q & = \Bb^Q + \ol{\pd}\Bb^Q 
   + \pd_\eb \Bb^Q :(\eeb{\bar\ub }-\eeb{\upstep})\circ \circpare{~} \\
   & \quad + \sum_{P=f,c} \pd_\eb M_{fP} (\bar p_P - \ppstep_P)\circ \circpare{~} 
   ,\\
   \ol{M}_{fc} & =  M_{fc} + \ol{\pd}M_{fc} + \sum_{P=f,c}\pd_{p_c}  M_{fP}(\bar p_P - \ppstep_P)\circ \circparpP{~}{c} \\
   & \quad + \pd_{p_c} \Bb^F :(\eeb{\bar\ub }-\eeb{\upstep})\circ \circparpP{~}{c} 
   \;,\\
  \ol{\Kb} & = \Kb^0 + \ol{\pd}\Kb^0\;,\\
  \ol{\Gb} & = \pd_\eb \Kb^0 (\nabla \bar p_f - \fb)\circ \circpare{~}\;, \\
  \ol{\Qb}_Q& = \pd_{p_Q} \Kb^0 (\nabla \bar p_f - \fb))\circ \circparpP{~}{Q}\;.
  \end{split}
\end{equation}

\subsubsection{\newSecE{Iterative scheme for solving the nonlinear two-scale problem}}\label{sec-iter-nonlin}
With reference to the above introduced incremental formulation we shall consider the subproblems for time incremets, now discretized also in space using the FE method. Hence, the weak-form problem \eq{eq-S25TD} with the residual \eq{eq-rif1}, and the incremental subproblem \eq{eq-rif2} with the differential \eq{eq-sts6i-V2} are replaced by the matrix formulation arising  from the FE discretization: at each time level $t_k$, $k = 1,2,\dots$, the nonlinear problem to find the state $\sbm$ such that
\begin{equation}\label{eq-fed1}
  \begin{split}
    \fbm^k(\sbm,\sbm\prevstep) & = \zerobf\;,\quad \mbox{ with } \sbm \approx \ssb^k\;,\sbm\prevstep \approx \ssb^{k-1} \;.
\end{split}
\end{equation}
It is solved by the following Newton-Raphson iterations: starting with $\sbm := \sbm\prevstep$, compute  $\dlt\sbm$, such that
 \begin{equation}\label{eq-fed2}
  \begin{split}   \Gbm^k(\sbm,\sbm\prevstep)\dlt\sbm  & =  -\fbm^k(\sbm,\sbm\prevstep)\;,\quad
\Gbm^k(\sbm,\sbm\prevstep) = \pd_\sbm\fbm^k(\sbm,\sbm\prevstep)\;,\\
\sbm & := \sbm + \dlt\sbm\;,\quad \sbm^{k}\approx\sbm
\end{split}
 \end{equation}
 where $\Gbm^k$ is the Jacobian matrix of the mapping $\sbm\mapsto\fbm^k(\sbm,\sbm\prevstep)$. Subprobem \eq{eq-fed2} constitute the main step of Algorithm~\eq{alg:micro-macro},
 where the micro-problems  are solved ``off-line'' and independently of the specific macroscopic problem. In the micro-problem stage, the characteristic responses $\psibf^k$,  $\omegabf^{ij}$, $\omegabf^F$, and $\omegabf^C$ are solved using the initial, unperturbed reference cell $Y$, then the homogenized coefficients $\Hop^0$ and also the sensitivieties $\delta_\sbm \Hop^0$ are evaluated using \eq{eq-S24}, \eq{eq-S40i}, and \eq{eq-S41i}. These are used in the macrsoscopic stage of the algorithm to solve the macro-problem by the interations, whereby the Jacobian $\Gbm^k(\tilde\sbm,\sbm\prevstep)$ is updated in each iteration using $\Hop^0$ and $\delta_\sbm \Hop^0$, such that the tangential coefficients $\tilde \Hop(\sbb|_e) = \Hop^0 + \delta_{\sbm|_e} \Hop^0 \circ \sbb|_e$ are defined for each element of the FE discretization.

\begin{algorithm}[!ht]
 \caption{Solution of the two-scale nonlinear problem}\label{alg:micro-macro}
    \begin{algorithmic}
      \Procedure{Pre-computing -- Micro-problems}{}
        \State Compute characteristic responses $\psibf^k$,  $\omegabf^{ij}$, $\omegabf^F$, $\omegabf^C$ from \eq{eq-S3}, \eq{eq-mip1}, \eq{eq-mip2}
        \State Evaluate homogenized coefficients $\Hop^0$, \eq{eq-HC1}
        \State Evaluate coefficient sensitivities $\delta_\sbm \Hop^0$, \eq{eq-HC4a} -- \eq{eq-HC3}, \eq{eq-S11}
      \EndProcedure
      \Procedure{Macro-problem}{} 
        \State Set initial state: $\sbb^0 \leftarrow \sbb^{init}$
        
        \For{all loading steps $k=1,2,\dots$}
            \State Given: $\widehat{\fb}^k$, $\widehat{\bb}^k$, $\sbb^{k-1}$
            \State Set: $\sbb\prevstep \leftarrow \sbb^{k-1}$,
            \State Set initial guess: $\sbb \leftarrow \sbb^{k-1}$
            \Do
            \State Update perturbed coefficients:
            \State \quad $\tilde \Hop(\sbb|_e) = \Hop^0 + \delta_{\sbb|_e} \Hop^0 \circ \sbb|_e$ for each element $e$. 
            \State Evaluate residual vector $\fbm^k(\sbm,\sbm\prevstep)$, see \eq{eq-S25-psi}.
            \State Calculate state increment $\delta \sbb$ by solving \eq{eq-fed2},
            $$\Gbm^k(\sbm,\sbm\prevstep)\dlt\sbm   =  -\fbm^k(\sbm,\sbm\prevstep)$$
            \State Update state: $\sbb \leftarrow \sbb + \delta \sbb$
            \doWhile{$\vert \delta \sbb \vert > tolerance $}
            \State Set: $\sbb^{k} \leftarrow \sbb$
            \State (Optional step) For a given $\varepsilon_0$ and $\sbb$
            \State \quad reconstruct microscopic responses in the local cells $Y(x)$
            \State \quad for selected positions $x \in \Om$. 
        \EndFor
      \EndProcedure
    \end{algorithmic}
\end{algorithm}

\input{aux-infl-simulations.tex}
\input{infl_simulations_dns.tex}

\section{Conclusion}

In this work, we proposed a new model describing fluid flow and inflation of a
double porosity fluid saturated medium. The model is derived using the
two-scale homogenization of the fluid-structure interaction in the periodic
microstructure consisting of the elastic skeleton incorporating the fluid
channels and the fluid inclusions. These two porosities are connected by
automatic admission and ejection valves which enable to inflate, or deflate the
inclusions depending on the difference of the pressure between the two
porosities. The asymptotic homogenization of the micromodel yields a kind of
the Biot continuum described by the displacement and the two pressures related
to the double porosity. Optionally, the fluid channels are subdivided into
compartments mutually separated by semipermeable interfaces yielding the
pressure discontinuity at the pore level, affecting the characteristic response
of the microstructure by the induced pressure jumps on the interfaces, thereby,
also the hydraulic permeability of the macroscopic Darcy law, \chE{\cf \cite{GrisoRohan2014} where a similar pressure discontinuity was treated in the context of double-porosity media.} Letting the
interface permeability depend on the solution (the local strain, or the
pressure jump), can lead to the ``valve effect'' at the micro-level, in contrast
with the valves separating the two porosities. Besides the nonlinearity induced
by the valves, we have considered also the nonlinearity associated with the
Eulerian approach to respect the deforming microconfigurations as the reference
frame to describe the equilibrium and the mass conservation. Using the
treatment proposed in \cite{Rohan-Lukes-CaS2024}, this nonlinearity yields
deformation-dependent effective properties of the homogenized material.
Efficient numerical treatment is based on the sensitivity analysis based
approximation of the homogenized coefficients involved in the equations of the
incremental subproblem. In the context of the small deformation, this enables
to avoid recalculating the local problems in the deformed RVEs otherwise needed
when dealing with nonlinear two-scale problems. At the same time, this
treatment can be applied in a case of slowly varying periodicity of the
microstructures.


Efficiency of the proposed Algorithm \ref{alg:micro-macro} for the two-scale nonlinear problem has been tested.
Solving the corrector problems in the microscopic domain $Y$ defined in
Section~\ref{sec:numex_dns} and evaluating the homogenized coefficients takes
approximately 6 seconds on a PC for the finite element mesh consisting of 9800
tetrahedral elements and 2115 nodes. The computational algorithm at the
macroscopic level involves the time-stepping loop with embedded iterations. The
computational time for 100 time steps is about 128 seconds, but almost 80\% of
that time is consumed by the ``Optional Step'' -- the reconstruction of the displacement, fluid pressure
and velocity fields at the microscopic level using the local characteristic
responses combined with the macroscopic solutions. In contrast, the DNS
solution time was approximately 1750 seconds. Thus, the speedup achieved by the
computational homogenization is nearly 14 times.

The presented numerical examples illustrate the functionality of the considered
smart double porosity medium to be inflated by the fluid pressure prescribed on
the inlet/outlet boundary. We have also demonstrated how this ``porous
metamaterial'' can be used to generate bending of a beam constituted by this
material. \chE{The valves have been introduced in the microstructure using the phenomenological approach. Their refined modelling at a sub-scale is challenging and also an important issue to approach realistic applications of the inflatable metamaterial model.}

\chE{Recent advances in 3D printing technologies open new possibilities in design of porous metamaterials intended for fluid transport, or shape morphing with various potential applications.
This paper offers a homogenization-based  modelling approach which will enable a systematic design of such complex structures. Moreover, the proposed model can be extended towards locally controlable metamaterials due to electromechanical transducers  \cite{Rohan-Lukes-CaS2024,Rohan-Lukes-PZ-porel}, which may support new ideas in construction of soft robots, or medical care devices -- inflatable orthosis with programmable pressure (in time and space) to heal vascular deceases, or just soft adaptable and wearable human-machine interfaces, like smart gloves. We have in mind new designs and design concepts, for which we propose the modelling \& computational tools. Specific structures and applications will be considered in our forthcoming research.}

\paragraph{Acknowledgment}
The research has been supported by the grant project GA~22-00863K and also in a part by project GA~24-12291S of the Czech Science Foundation.

\input{aux-Appx-inflatable-R.tex}

\bibliographystyle{elsarticle-harv}


\end{document}

%% file: aux-infl-simulations.tex
\section{Numerical simulations}\label{sec:numex}

Using the derived two-scale model of the homogenized inflatable structures we
aim to illustrate some interesting properties and functionality of this
two-phase porous medium. To this purpose we consider a particular
microstructure of the porous material which constitutes one of the bi-layer
``macroscopic'' structure. It attains a form of a prism with square
cross-section consisting of the porous layer which is is attached to an elastic
substrate -- the second layer, see Fig.~\ref{fig:num-macro-geom}. Such a
construction inspired by the bi-metallic beam which bends in response to the
temperature changes; in the present case, as reported in
Section~\ref{sec:numex-pressure-driven}, we shall explore analogous behavior in
response to the inflation pressure. Besides the deformation induced by the
pressure, it is of interest to observe the time-space distribution of the
inflation pressure in the inclusions and the fluxes through the admission and
ejection valves. In Section~\ref{sec:numex-nonlin}, we consider the
modification of the model concerning the geometric nonlinearity due to the
equilibrium imposed in the deformed configuration, as explained in
Section~\ref{sec-nonlin-II}, which leads to the deformation-dependent
homogenized coefficients involved in the macroscopic incremental formulation.
Finally, in Section~\ref{sec:numex-channel-interface}, we study the variation
of the homogenized hydraulic permeability as a function of the permeability of
a semipermeable membrane embedded in the fluid channel.

\subsection{Preliminary settings and modelling options}\label{sec:problem-settings}

The reported numerical studies involve data on the microstructure, parameters
of the specimen and applied loads in the virtual testing. The numerical
approximation schemes are based on discretization schemes using the finite
element in space and finite differences in time. All the mathematical model
parameters are specified just here below.

\begin{figure}[ht]
    \centering
    \includegraphics[width=0.8\linewidth]{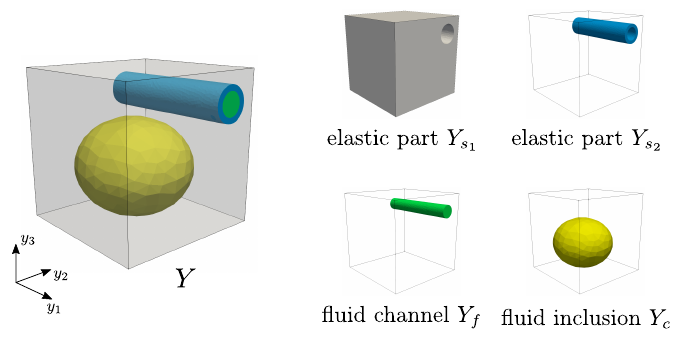}
    \caption{Decomposition of the representative periodic cell $Y$ into subparts $Y_{s_2}$, $Y_{s_2}$, and $Y_c$.}
    \label{fig:num-micro-geom}
\end{figure}

\paragraph{Geometry -- domains, boundaries}
The microscopic domain $Y = ]0,1[^3$ is composed of the solid part $Y_s$, fluid
channel $Y_f$, and fluid inclusion $Y_c$. To better control the inflation
effect due to the increase of the inclusion volume, it proves to be
advantageous to suppress the inflation of the skeleton due to the increased
pressure in the fluid channel $Y_f$. Therefore, to reduce this undesirable
effect, but to retain the compliance needed to inflate the inclusions, we
consider a stiffer shell $Y_{s_2}$ of the channel and a softer material in the
rest of the solid domain $Y_{s_1}$ in which the inclusions are embedded ($Y_s =
Y_{s1} \cup Y_{s2}$), as depicted in Figure~\ref{fig:num-micro-geom}.

\begin{figure}[ht]
    \centering
    \includegraphics[width=0.99\linewidth]{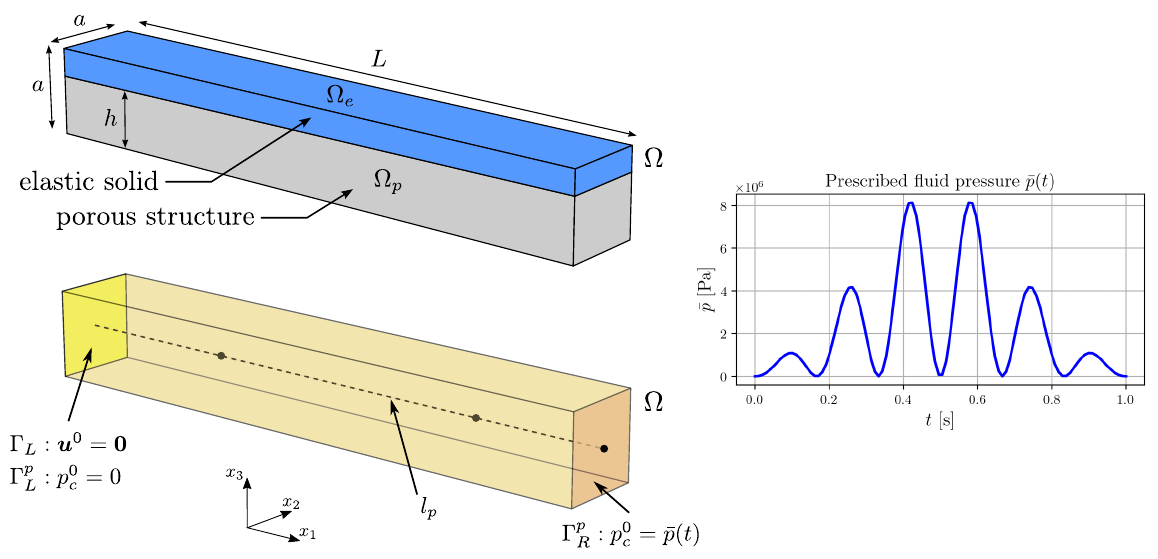}
    \caption{Macroscopic domain $\Omega = ]0,L[ \times ]0,a[^2$
             and the applied pressure $\bar{p}(t)$ in time interval $t\in\langle0, 1\rangle$\,s.
             Three point with relative positions $x_p = 0.27$, 0.75, and 1.0 are marked on
             on the probe line $l_p$.}
    \label{fig:num-macro-geom}
\end{figure}

At the macroscopic level, we consider a prismatic specimen -- block sample $\Om
\subset \RR^3$ with the square cross-section $a \times a$ of dimensions and
$a=0.005$ \,m, and the length $L = 0.1$ \,m, such that $\Om = \omega \times
]0,a[$, with $\omega = ]0,L[ \times ]0,a[ \subset \RR^2$ see
Figure~\ref{fig:num-macro-geom} left, consisting of the porous part $\Omega_p =
\omega \times ]0,h[$ and an elastic part $\Omega_e = \omega \times ]h,a[$ (the
substrate), where $h = 5/7 a$. To specify BCs, the
boundary is decomposed, as follows: $\pd \Om = \Gamma_L \cup \Gamma_R \cup
\Gamma_S$, the three disjoint parts representing the left and right sides, and
the block ``sleeve'' $\Gamma_S$, where $\Gamma_L = \{x \in \pd \Om|\,x_1 =
0\}$, $\Gamma_R = \{x \in \pd \Om|\,x_1 = L\}$. Moreover, we have to
distinguish the portions of $\Gamma_L, \Gamma_R$, and $\Gamma_S$ constituting
boundaries of $\Om_p$ and $\Om_e$. Henceforth, we introduce $\Gamma_K^s$, with
$K := L,R,S$ and $s := p,e$. Finally, the two subparts share the interface
$\Gamma_{ep} = \omega \times (0,0,h) = \pd\Omega_e \cap \pd\Omega_p$.

\paragraph{Material parameters of the microstructure}
The material properties
of the elastic skeleton and the fluid constituting the bi-phase microstructures are introduced in 
Tab.~\ref{tab:num-material-properties}.
Besides, for the elastic substrate we take the Young's modulus $E = 2$\,GPa and the Poisson's
ratio $\nu = 0.4$, Further parameters related to the valves 
comprise permeability parameters $\bar{\kappa}_A$ and $\bar{\kappa}_E$ which (both) are
set to $10^{-7}$\,Pa$^{-1}$s$^{-1}$, and the pressure gauge $\Delta P = 3$\,MPa. Note that the characteristic responses
of the microstructure are computed with the rescaled viscosity $\bar\mu$
obtained for a given microstructure size $\varepsilon_0 = 0.0025$.

\paragraph{Boundary conditions and loading}
Although we consider a modified problem due to the elastic part $\Om_e$, we
do not need to specify all details for the associated modified formulation. Obviously, it is straightforward to
consider a degenerate porous medium with the porosities $\phi_f = \phi_c = 0$
yielding an elastic (possibly homogeneous) medium For the macroscopic problem
\eq{eq-S25-fc}, the Dirichlet type BCs are specified accordingly the
admissibility sets $\Ucalbf(\Om)$ and $\Qcal(\Om)$, see \eq{eq-UPspaces}, where
$\pd_u\Om \equiv \Gamma_L$ and $\pd_p\Om \equiv \Gamma_L^p \cup \Gamma_R^p$
(recall the meaning of superscript ${}^p$ referring to the porous part of
geometric entities). On $\Gamma_L$, the structure is fully anchored and $p_c^0 = 0$ on $\Gamma_L^p$.
On the opposite side, $\Gamma_R^p$, the porous structure is loaded by a time-varying pressure $p_c^0 = \bar{p}(t)$.
By the consequence of the boundary decomposition,
$\pd_\sigma \Om = \Gamma_S \cup \Gamma_R$; zero traction and volume forces are
considered, so $\fbav = \bbav = \zerobf$. The porous specimen $\Om_p$ is
impermeable on $\Gamma_S^p$, including the interface $\Gamma_{ep}$. The
prescribed inflation pressure on $\Gamma_R^p$ is given by 
\begin{equation}\label{eq:num_boundary_pressure}
    \bar{p}(t) = \sin^2(k\pi t)
        \frac{1}{\sqrt{2\pi b^2}} \exp^{-\frac{(t - c)^2}{2b^2}},
\end{equation}
where we chose $k=2/3$, $b=1/5$, and $c=1/2$, see
Figure~\ref{fig:num-macro-geom} right.

\paragraph{Numerical models and implementation}
The model has been implemented and all the computations and simulations have been made
 in the Python based FE package {\it SfePy}: Simple Finite Elements in Python \cite{Sfepy_2019}. The Python code used to generate results shown in the following parts has been published in the
public repository \cite{Lukes_code_2025}.
The macroscopic problem \eq{eq-S25-fc} and the microscopic subproblems for the
unknown correctors \eq{eq-S3}, \eq{eq-mip1}, \eq{eq-mip2} are spatially
discretized using the \emph{finite element} (FE) method. Lagrangian elements are employed
for spacial approximations at both scales: the displacement and fluid pressure fields
are approximated by the piecewise linear Q1 elements, while the velocity field by
the piecewise quadratic Q2 elements. The macroscopic problem is solved by a time-stepping algorithm in the interval $t \in \langle0, 1\rangle$\,s. The time derivatives in the macroscopic
equations \eq{eq-S25-fc} are approximated by the backward finite differences with the time step $\Delta t = 0.01$\,s.
The nonlinear subproblem is solved iteratively within each time increment (of the time-stepping loop) to account for
the nonlinearity introduced through the admission and ejection valves, as described in Section~\ref{sec-geom}, and by the varying coefficients defined in
Section~\ref{sec-increment}. Despite the macroscopic nonlinearity, due to the sensitivity analysis based approximation of the  homogenized coefficients in the macroscopic incremental algorithm, the two-scale
computational algorithm remains fully decoupled: all the micro-level problems involving the characteristic responses, \eq{eq-S3}, \eq{eq-mip1}, and
\eq{eq-mip2},  evaluation of the homogenized coefficients and their sensitivity analysis \eq{eq-S40i}-\eq{eq-S41i} are made off-line independently of solving a specific macroscopic problem. 


\paragraph{Results presentation, post-processing}
Due to the fluid viscosity being scaled by $\veps^2$, the simulations must be interpreted for a give microstructure size; we consider $\varepsilon_0 = 0.0025$ in all the simulations.
 In order to visualize the spatial
distribution of the selected macroscopic quantities, we introduce a line probe
$l_p$ aligned in the $x_1$ direction through the sample (the axis of the porous prism) and define a relative coordinate $x_p$ that takes the value 0 at $\Gamma_L$ and 1 at $\Gamma_R$. So, in the examples, we refer to the points on the $l_p$ by coordinates $x_p \in [0,1]$.
Using the characteristic responses, by virtue of formulae \eq{eq-A10a} and \eq{eq-L2b}, the two-scale functions can be evaluated and the displacement, fluid pressure and velocity fields in the microstructure can be reconstructed. 


\begin{table}
    \begin{center}
    \begin{tabular}{|l|c|}
        \hline
        {elastic part $Y_{s_1}$} & E = 20\;MPa, $\nu$ = 0.49\\
        \hline
        {elastic part $Y_{s_2}$} & E = 200\;MPa, $\nu$ = 0.3\\
        \hline
        {fluid in $Y_f$ and $Y_c$ (water)} & $\gamma = 1.0 / (2.15 \cdot 10^9)$\;Pa$^{-1}$, $\mu^\veps = 8.9 \cdot 10^{-4}$\;Pa\;s\\
        \hline
    \end{tabular}
    \end{center}
    \caption{Material properties of the porous structure.
             Note that $\bar\eta = \bar\eta^\veps/\veps$ is computed for a given scale $\veps:=\veps_0$.}
    \label{tab:num-material-properties}
\end{table}

\subsection{Pressure driven inflation}\label{sec:numex-pressure-driven}
As stated above, the simulation of the inflation due to varying ``inlet'' pressure was performed for the model parameters and numerical approximations described above, see Table~\ref{tab:num-material-properties}. In this study, we consider the macroscopic model \eq{eq-S25-fc} involving  constant homogenized coefficients defined by expressions \eq{eq-S7a}, \eq{eq-HC1}. Below, in Section~\ref{sec:numex-nonlin},
we compare simulations performed with this ``Lagrangian configuration based'' model with those obtained due to the ``Eulerian configuration based'' model.

\begin{figure}[ht]
    \centering
    \includegraphics[width=0.49\linewidth]{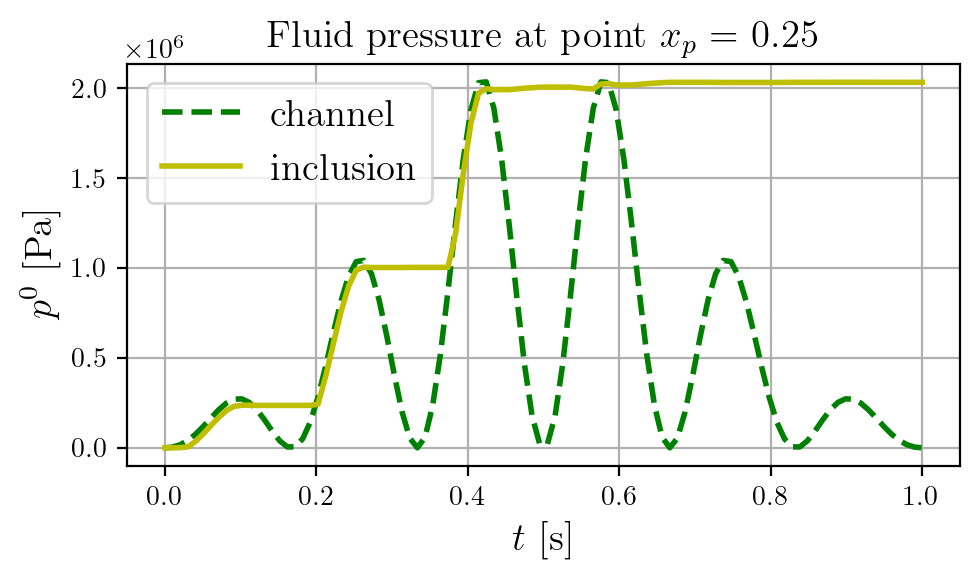}
    \hfill
    \includegraphics[width=0.49\linewidth]{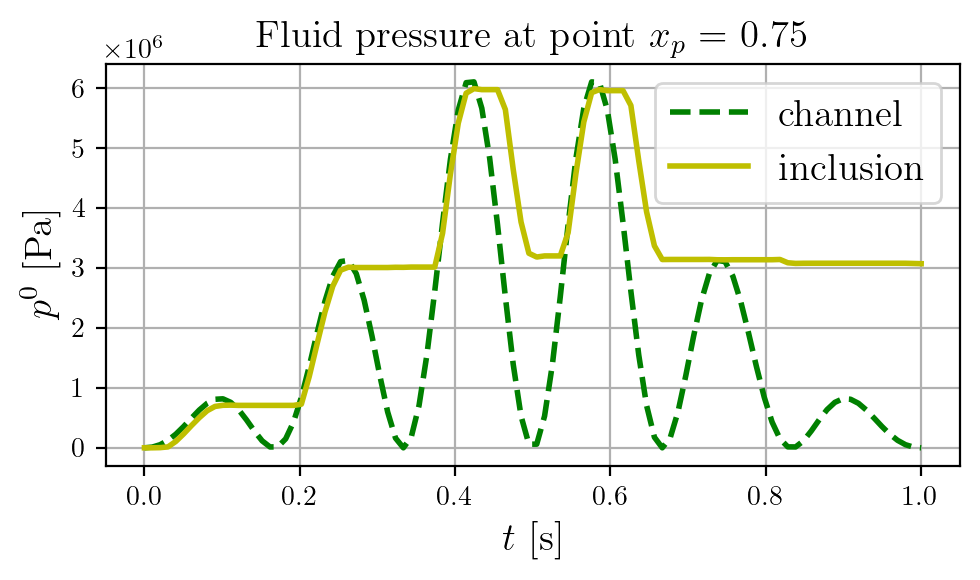}
    \caption{Time evolution of the channel ($p^0_c$) and inclusion ($p^0_f$)
             macroscopic pressures at two distinct macroscopic points $x_p=0.25, 0.75$.}
    \label{fig:num-2mat-pressure}
\end{figure}

\begin{figure}[ht]
    \centering
    \includegraphics[width=0.49\linewidth]{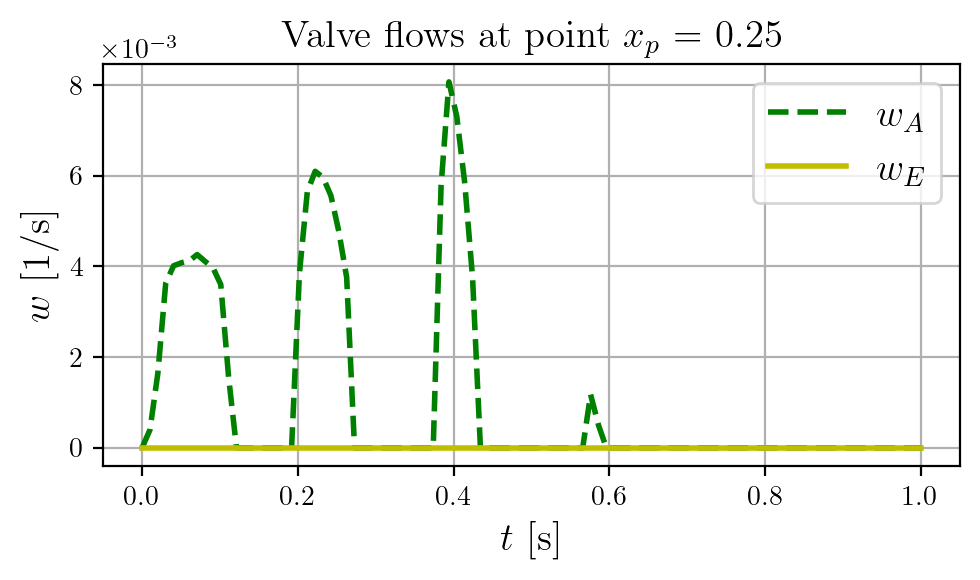}
    \hfill
    \includegraphics[width=0.49\linewidth]{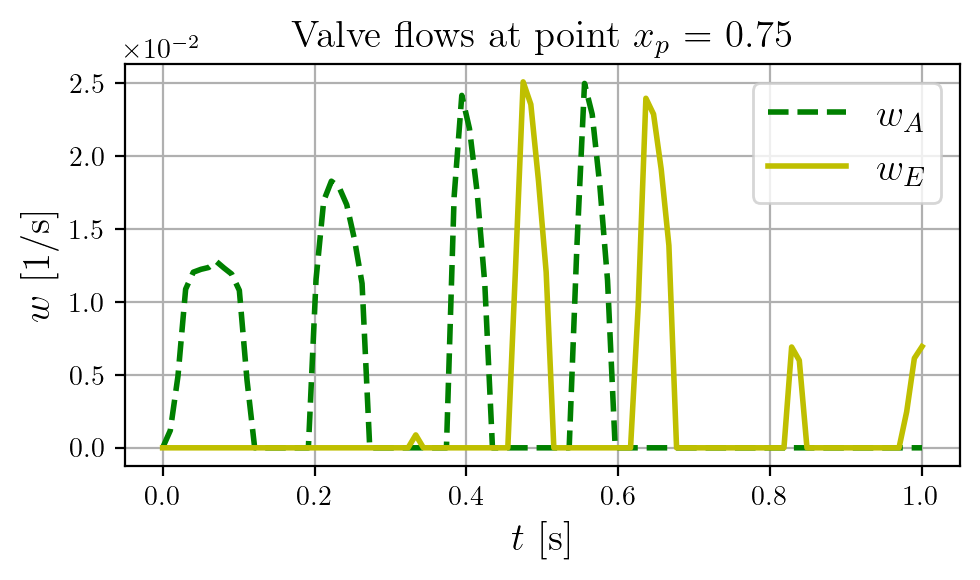}
    \caption{Time evolution of valve fluxes $w_A$ and $w_E$
             at two distinct macroscopic points $x_p=0.25, 0.75$.}
    \label{fig:num-2mat-valve-flux}
\end{figure}

\begin{figure}[ht]
    \centering
    \includegraphics[width=0.49\linewidth]{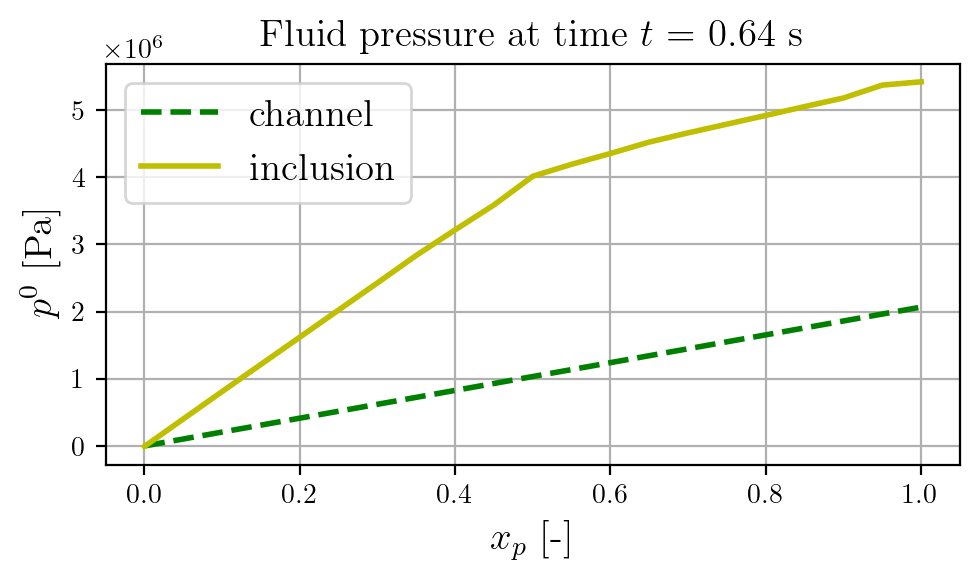}
    \hfill
    \includegraphics[width=0.49\linewidth]{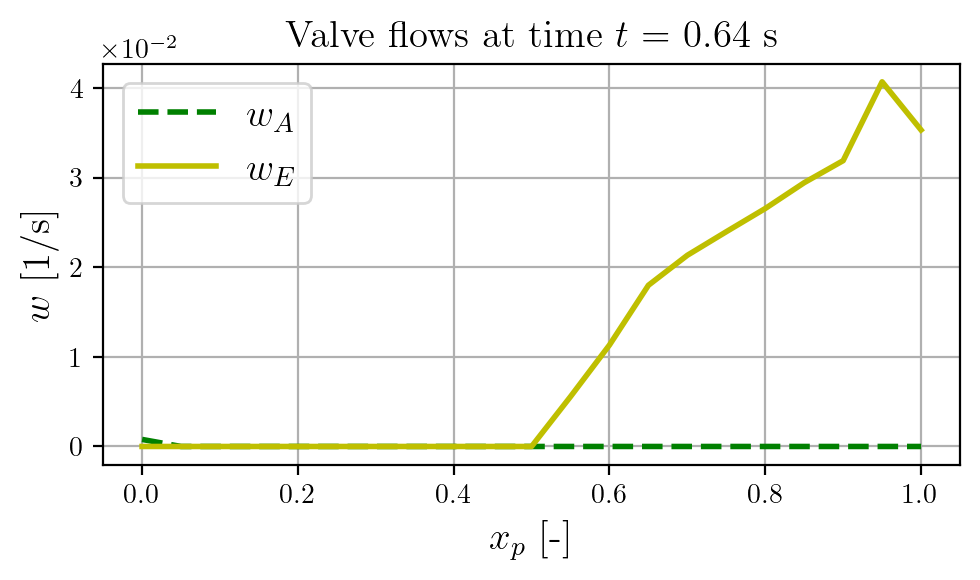}
    \caption{Spatial distribution of pressures $p^0_c$, $p^0_f$ (left)
             and valve flows $w_A$, $w_E$ (right) along line $l_p$.}
    \label{fig:num-2mat-pressure-flux}
\end{figure}

\begin{figure}[ht]
    \centering
    \includegraphics[width=0.49\linewidth]{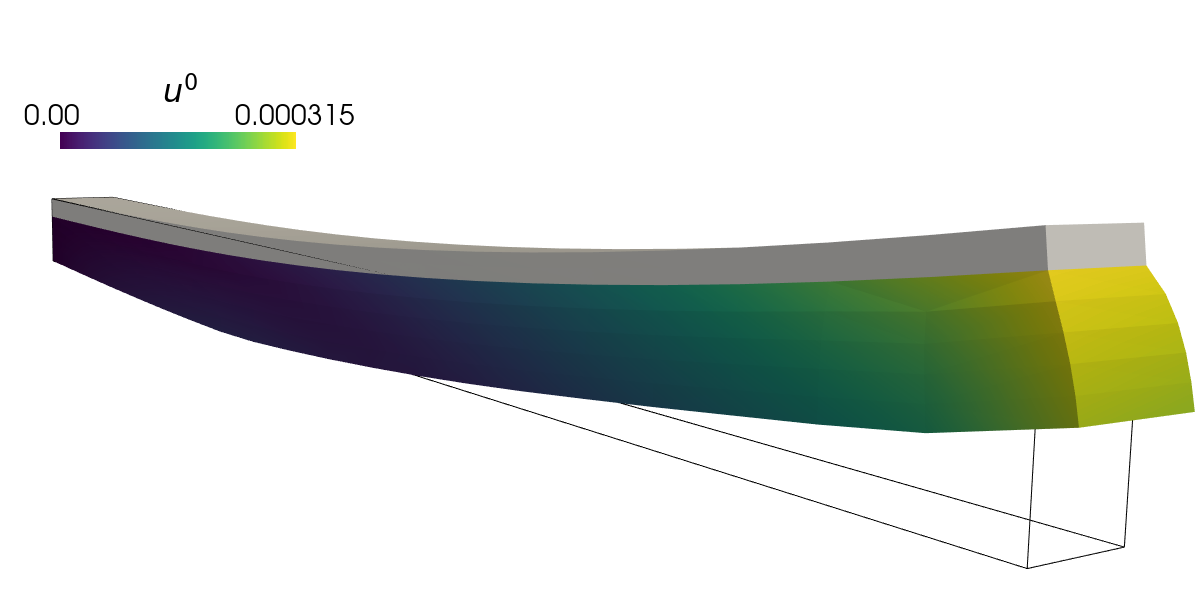}
    \hfill
    \includegraphics[width=0.49\linewidth]{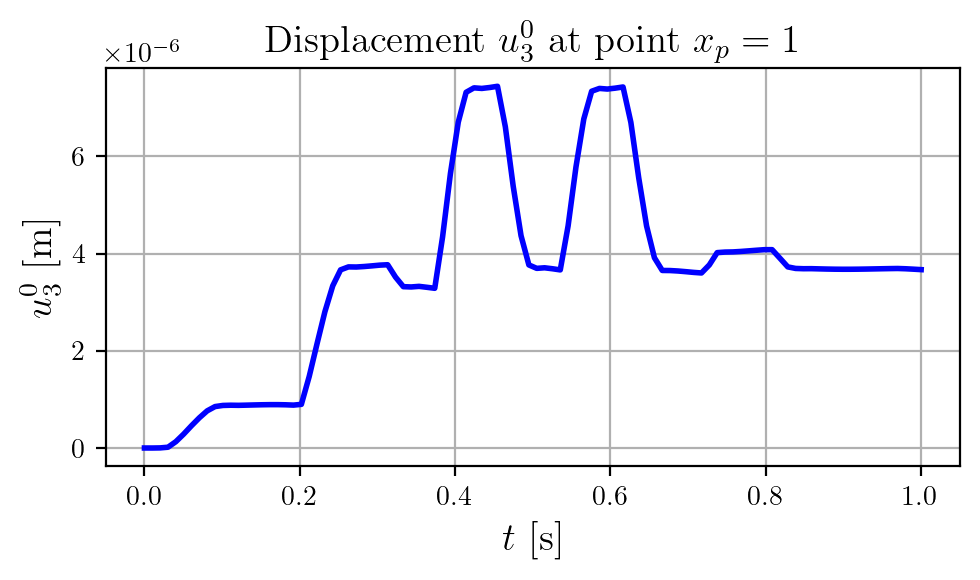}
    \caption{Spatial distribution of macroscopic the  displacements $\ub^0$ in time $t=1$\,s (left) and
              time evolution of the displacement of the sample end point $x_p = 1$ (right).}
    \label{fig:num-2mat-u_def}
\end{figure}

\begin{figure}[ht]
    \centering
    \includegraphics[width=0.49\linewidth]{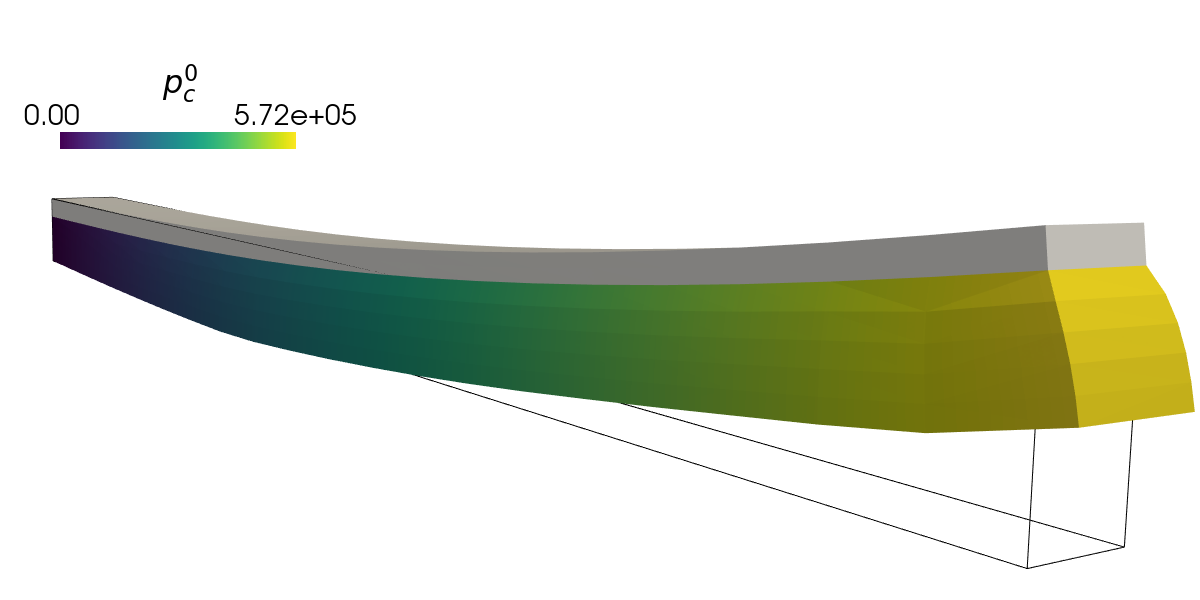}
    \hfill
    \includegraphics[width=0.49\linewidth]{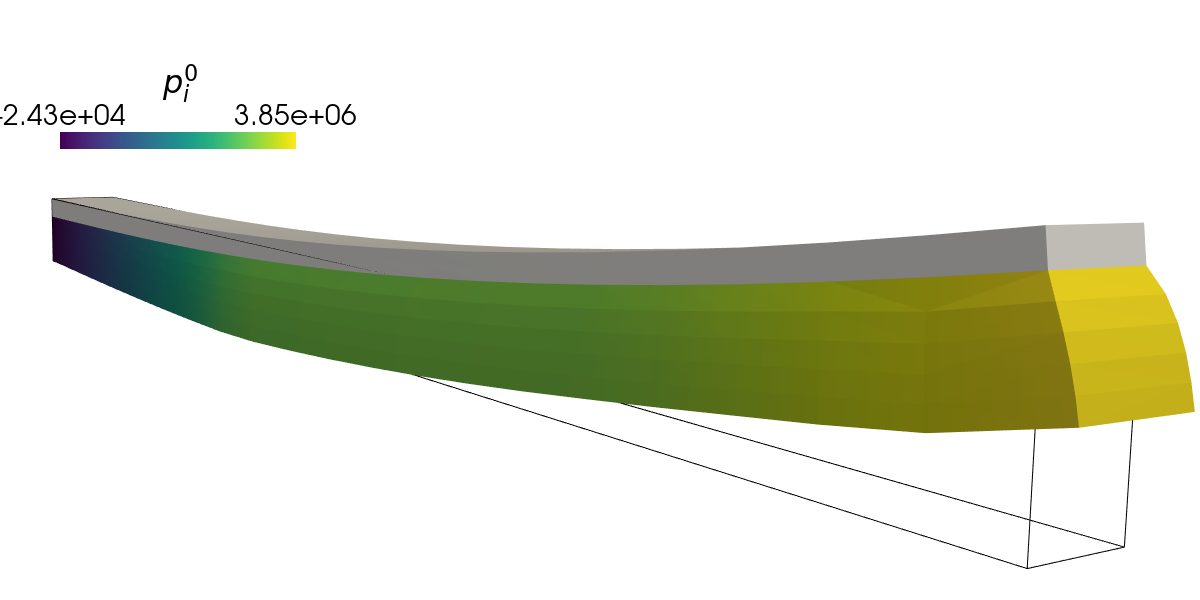}
    \caption{Spatial distribution of macroscopic pressures $p^0_c$, $p^0_f$ in time $t=0.8$\,s.}
    \label{fig:num-2mat-p_def}
\end{figure}

The macroscopic responses to pressure loading are presented in
Figures~\ref{fig:num-2mat-pressure}--\ref{fig:num-2mat-p_def}. The time
evolution of the channel pressure $p_c^0$ and inclusions pressure $p_i^0$ at
two points on the line probe is plotted in Fig.~\ref{fig:num-2mat-pressure}.
Whereas at places of the lower pressure $p^0_c$, the inclusions are gradually
inflated by pressure waves in the fluid channels, see Fig.~\ref{fig:num-2mat-pressure} left,
at places with the higher channel pressure the inflation and deflation of the inclusions
are changing due to the chosen value of $\Delta P = 3 \cdot 10^6$\,Pa,
see Fig.~\ref{fig:num-2mat-pressure} right.
Recalling the limit mass conservation laws \eq{eq-lim-Mf} and \eq{eq-lim-Mc}, the fluid flows between the channels and inclusions are characterized by the
quantities $w_A = \bar{\kappa}_A \left[p_f^0 - p_c^0\right]_+$ and $w_E = \bar{\kappa}_E
\left[p_c^0 - p_f^0 - \Delta P\right]_+$, which are shown in
Fig.~\ref{fig:num-2mat-valve-flux} for two distinct points on the $l_p$ probe.
The zero flux $w_E$ from the inclusions to channels in
Fig.~\ref{fig:num-2mat-valve-flux} left indicates pure inflation at point
$x_p=0.25$ while the varying values $w_A$ and $w_E$, in
Fig.~\ref{fig:num-2mat-valve-flux} (right), mean the inflation and deflation
phases at point $x_p=0.75$.

The spatial distribution of pressures $p_c^0$, $p_f^0$ and flows $w_A$, $w_E$
in time $t=0.8$\,s along the probe line is plotted in
Fig.~\ref{fig:num-2mat-pressure-flux}. Figures \ref{fig:num-2mat-u_def} left
and \ref{fig:num-2mat-p_def} show the pressures and displacements in time
$t=0.8$\,s in the deformed macroscopic sample. Note that the pressure fields
are defined only in the porous part $\Omega_p$, while the displacement field is
valid in the whole domain $\Omega$. In Fig. \ref{fig:num-2mat-u_def}, the depicted time
variation of  displacement $u^0_z := u^0_3$ at the sample extremity ($x_p = 1$ on $\Gamma_R$)  reveals its correlation with the variation of the inclusion pressure in the whole structure.

\begin{figure}[ht]
    \centerline{
        \includegraphics[width=0.98\linewidth]{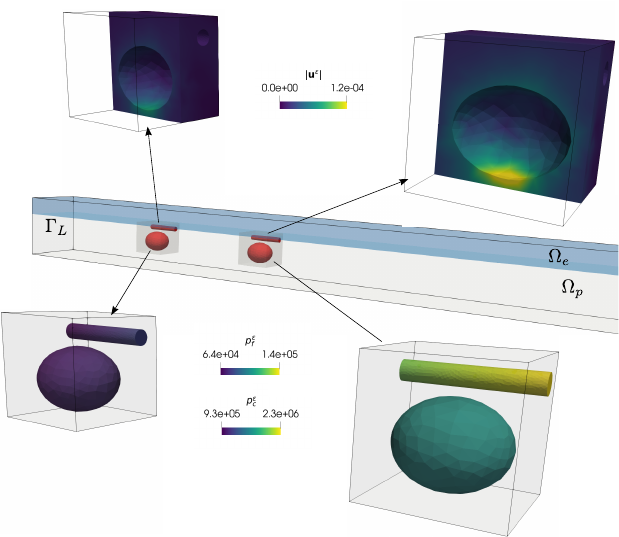}
    }
    %
    \caption{Reconstructed pressure fields $p_c^\varepsilon$, $p^\varepsilon_f$, and $\ub^\varepsilon$ in time $t=0.8$\,s.}
    \label{fig:num-2mat-reconstructions}
\end{figure}

Finally, for illustration, the reconstructed pressure fields $p_c^\varepsilon$,
$p_i^\varepsilon$, and $\ub^\varepsilon$ for the finite $\varepsilon_0$
at two points  
of the sample are shown in Fig.~\ref{fig:num-2mat-reconstructions}.
%

\subsection{Simulation with deformation-dependent homogenized coefficients (HC)}\label{sec:numex-nonlin}
\begin{figure}[ht]
    \centerline{
        \includegraphics[width=0.49\linewidth]{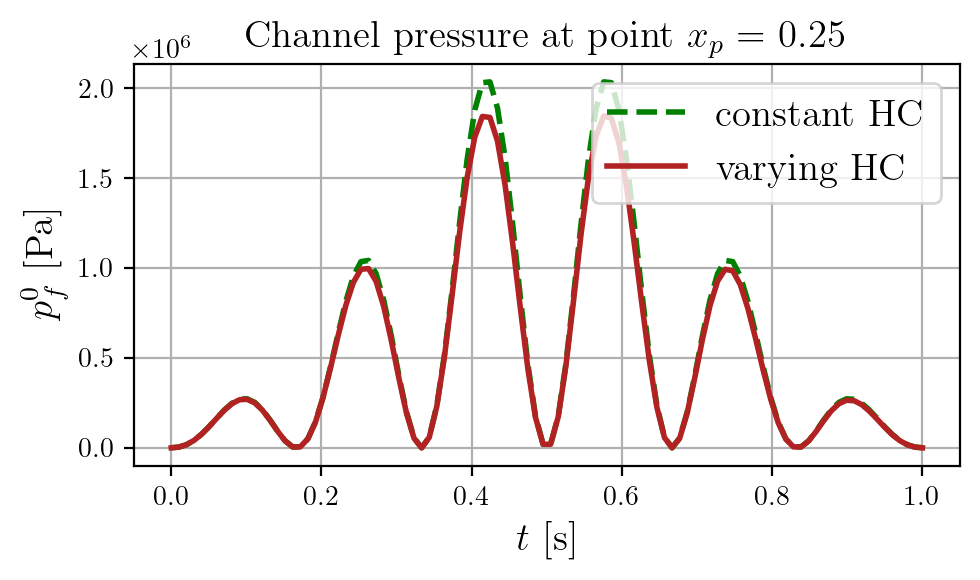}
        \hfil
        \includegraphics[width=0.49\linewidth]{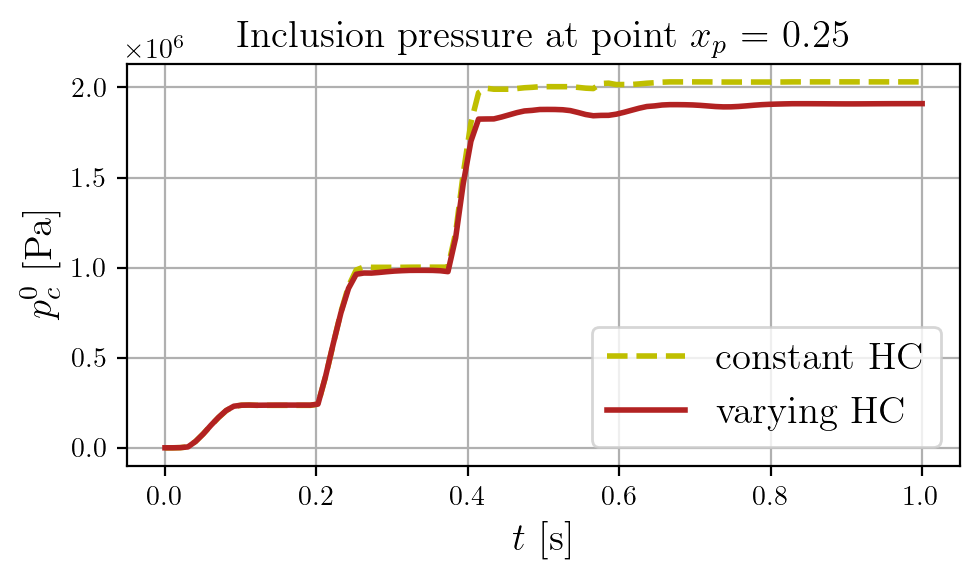}
        }
    \centerline{
        \includegraphics[width=0.49\linewidth]{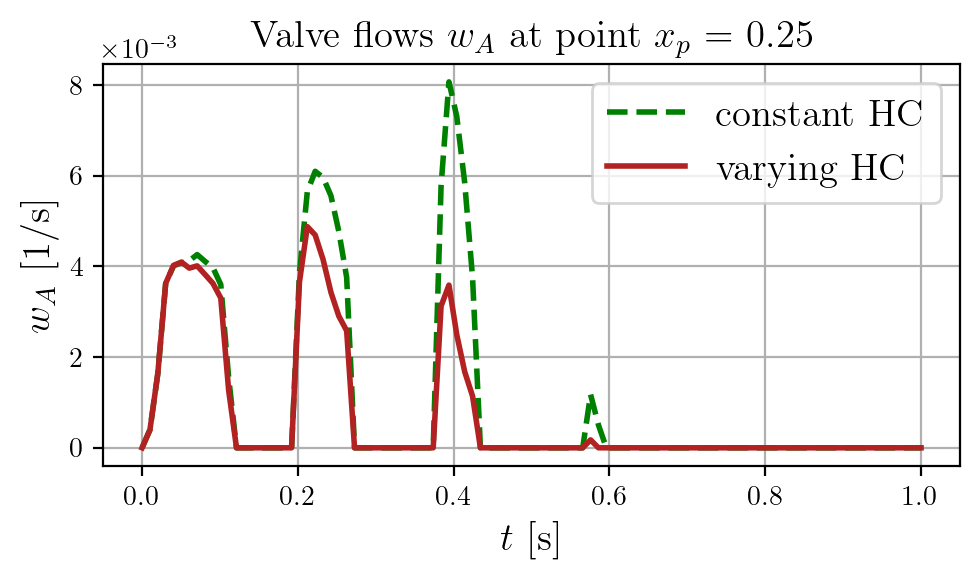}
        \hfill
        \includegraphics[width=0.49\linewidth]{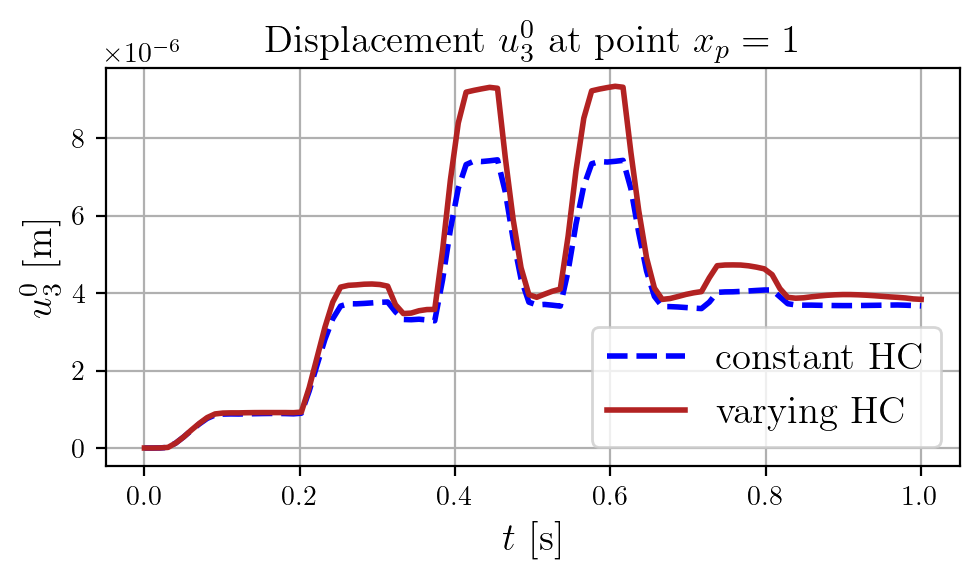}
    }
    \caption{Comparison of two models with constant (``constant HC'') and
             strain- and pressure-dependent (``varying HC'') homogenized coefficients:
             fluid channel pressures $p^0_f$ (top left),
             inclusion pressures $p^0_c$ (top right),
             valve flows $w_A$ (bottom left),
             and displacements $u^0_3$ at point $x_p=1$ (bottom right).}
    \label{fig:num-2mat-nonlin}
\end{figure}

We now consider the same inflation problem as the one reported above, however, we use model \eq{eq-rif1} which reflects the influence of the equilibrium and mass conservation associated with the deformed configuration, as explained in Section~\ref{sec-nonlin-II}.
Due to the sensitivity analysis based approximation of the HC involved in the nonlinear problem, the computational costs is only slightly increased compared to the simulations with fixed reference configurations. Recall, that the model is nonlinear any way due to the valves and requires an iterative algorithm  to compute solutions.
For sensitivities of the HC, the details are given in \Appx{sec-SA}.
Using this modified ``Eulerian configuration based'', we obtain solutions, referred to as ``E-model'' , characterized by varying HC. In 
Figure~\ref{fig:num-2mat-nonlin}, responses computed by the ``Eulerian configuration based'' and ``Lagrangian configuration based'', referred to by ``L-model'' are compared. 
%
As can be expected, the largest differences from the linear model occur at
locations of large macroscopic pressures and strains. In
Figure~\ref{fig:num-2mat-nonlin}, we compare the following macroscopic fields,
from the top-left to the bottom-right: fluid channel pressures $p^0_f$,
inclusion pressures $p^0_c$, valve flows $w_A$, and displacements $u^0_3$ at
point $x_p=1$. The results reveal that the pressure is overestimated by the L-model compared the corresponding solution of the E-model. Conversely, the displacements $u_3^0(x_p = 1,t)$ predicted by the L-model are smaller than those computed using the E-model. In other words, the E-model is more compliant that the L-model.

\subsection{Semipermeable interfaces in fluid channels}\label{sec:numex-channel-interface}
\begin{figure}[ht]
    \centerline{
        \includegraphics[width=0.32\linewidth]{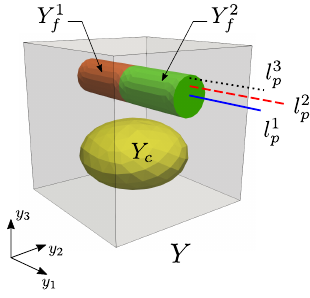}\hfill
        \includegraphics[width=0.6\linewidth]{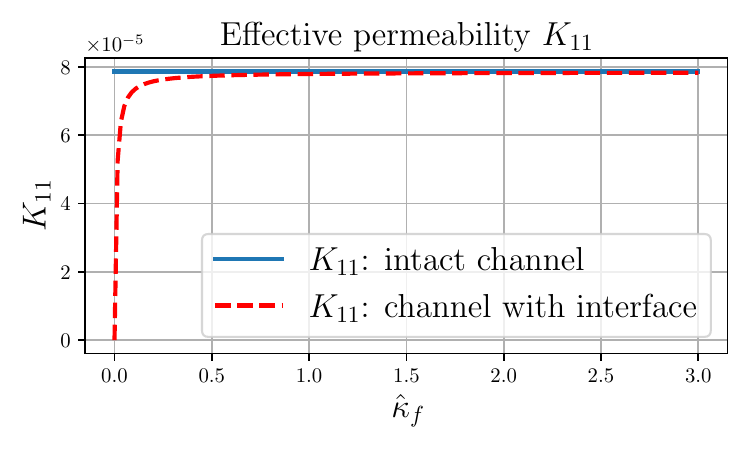}
    }
    \caption{Periodic cell $Y$ representing discontinuous fluid channel
             and the dependence of the homogenized permeability $K_{11}$ on
             the membrane permeability parameter $\hat\kappa_f$.}
    \label{fig:micro_fvalve}
\end{figure}

\begin{figure}[ht]
    \centerline{
        \includegraphics[width=0.49\linewidth]{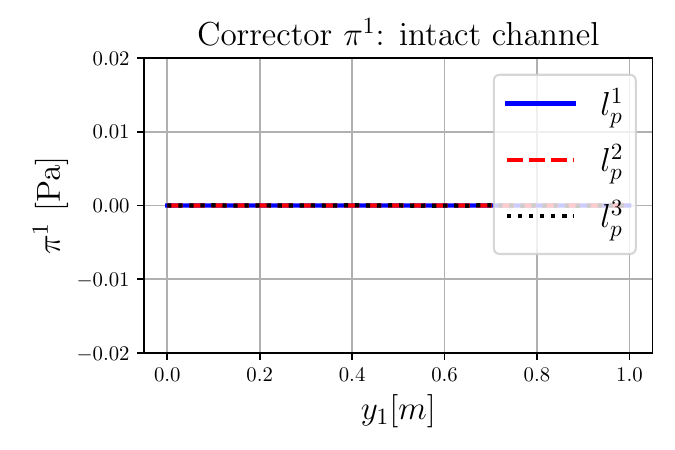}\hfill
        \includegraphics[width=0.49\linewidth]{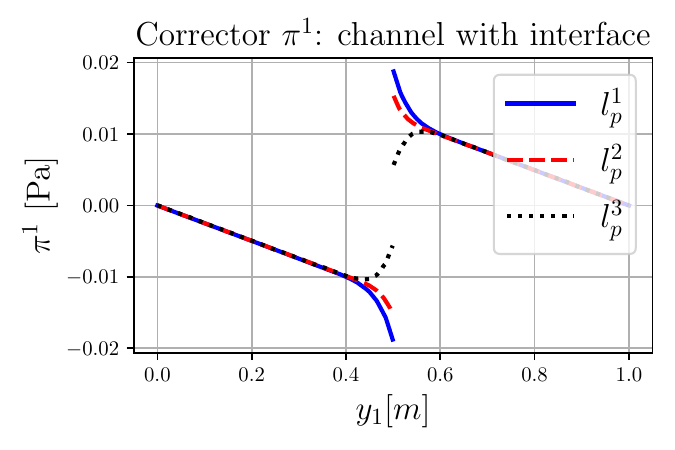}
    }
    \caption{Spatial distribution of corrector $\pi^1$ along line probes $l_p^1$, $l_p^2$, $l_p^3$ in $Y_f$
            for the intact fluid channel (left) and for the channel with semipermeable interface (right).}
    \label{fig:micro_fvalve_corr}
\end{figure}

In this section, we demonstrate the influence of the membrane permeability
parameter $\hat\kappa_f$ on the effective permeability tensor $\Kb$ defined by \eq{eq-S7a}.
Calculations are performed for $\hat\kappa_f \in [ 0, 3]$
on a periodic cell where the fluid channel $Y_f$ 
is divided
into two parts, as depicted in Figure~\ref{fig:micro_fvalve} left, where the associated domains are marked as red and green). Otherwise the microstructure is the same as in all the above reported simulations.
From the plot in Figure~\ref{fig:micro_fvalve} right, it can be seen that the
resulting homogenized permeability $K_{11}$ in the $y_1$ direction is zero for
$\hat\kappa_f = 0$, which means a closed channel due to impermeable interface $\Gamma_f$. For increasing $\hat\kappa_f$, permeability $K_{11}$ converges to the value
obtained intact channel, so without any interface $\Gamma_f$.
It is worthy to mention that interface imposing the pressure discontinuity disproportional to $\hat\kappa_f$. When this interface permeability is infinite (actually, no interface), pressure is continuous across $\Gamma_f$. Moreover, since the channel is straight, the Poiseuille flow is recovered, such that the pressure field is a linear, depending on $x_1$ only, hence the correction due to $p^1$ is zero, as seen in Figure~\ref{fig:micro_fvalve_corr}, left. On contrary, smaller $\hat\kappa_f$ is, larger discontinuity of $p^1$ is induced in the vicinity of $\Gamma_f$; correspondingly, the largest jump in the corrector $\pi^k$ appears on the channel axis, whereas it vanishes when approaching the wall $\Gamma_\fsi$, see Fig.~\ref{fig:micro_fvalve_corr}, right.

%% file: infl_simulations_dns.tex
%
%
%
%


\subsection{\newSecE{Validation of homogenized model}}\label{sec:numex_dns}

\begin{figure}[ht]
    \centering
    \includegraphics[width=0.49\linewidth]{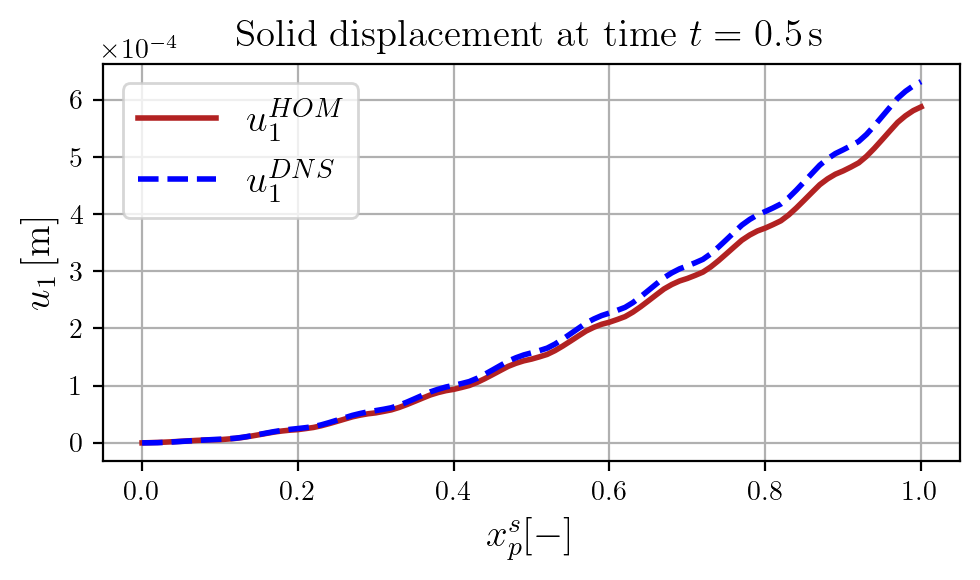}\hfil
    \includegraphics[width=0.49\linewidth]{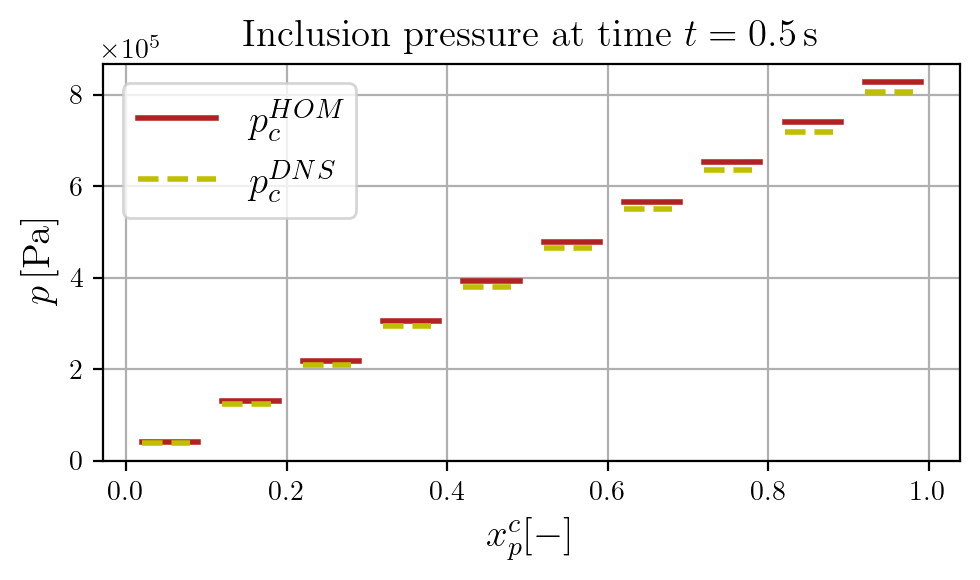}\\
    \includegraphics[width=0.49\linewidth]{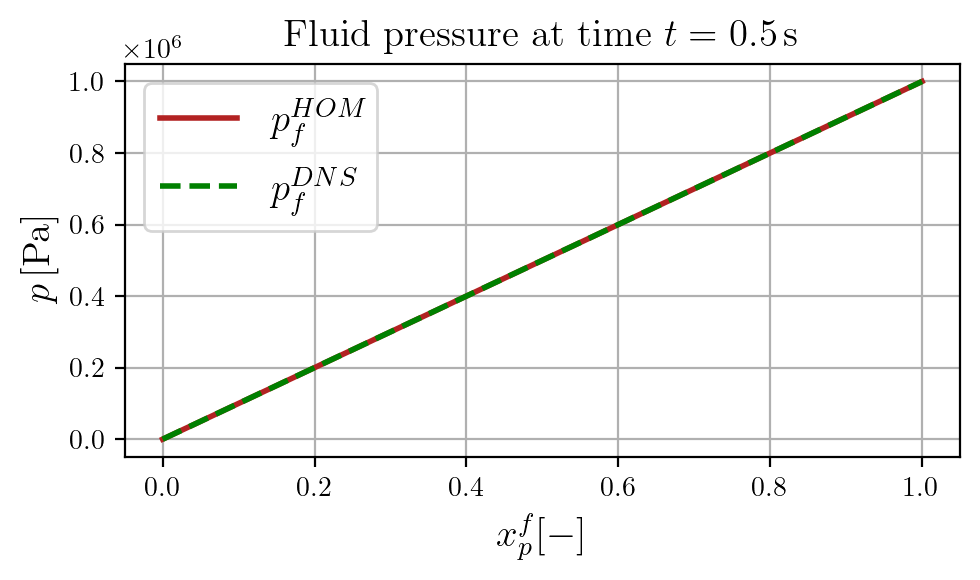}\hfil
    \includegraphics[width=0.49\linewidth]{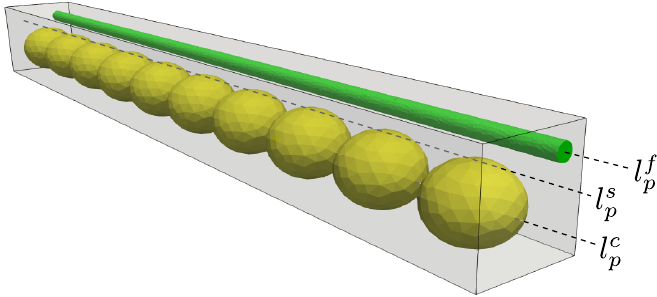}\\
    \caption{Comparison of the direct numerical simulation (DNS) and the homogenized model (HOM).
             Spatial distributions of displacements $u_1$, inclusion pressures $p_c$,
             and channel pressures $p_f$ along line probes
             $l_p^s$ (solid), $l_p^c$ (inclusion), $l_p^f$ (channel), and $l_p$ (homog. model)
             at time $t=0.5$\,s.}
    \label{fig:num-dns-line-probes}
\end{figure}

A direct numerical simulation (DNS) of an inflatable heterogeneous structure is
performed to validate the proposed homogenized model. The finite element model
used for this validation is build up by periodically repeating the periodic
unit cell, depicted in Fig.~\ref{fig:num-micro-geom}, along the $x_1$-axis.
Assuming the microstructure size $\veps_0 = 0.01$, the resulting sample
comprises of 10 copies of the RVE and has dimensions of $0.1 \times 0.01 \times
0.01$\,m, see bottom right subfigure in Fig.~\ref{fig:num-dns-line-probes}. For
simplification, the rigid compartment $Y_{s_2}$ is eliminated by assigning it
the same elasticity parameters as in the part $Y_{s_1}$. The boundary
conditions similar to those depicted in Fig.~\ref{fig:num-micro-geom} are
applied to the sample faces $\Gamma_L$: $\ub_0 = \zerobf$, $p_c = 0$,
$\Gamma_R$: $p_c = 10^6 \sin(\pi t)$, while periodic displacement conditions
are enforced on the remaining four faces.
This configuration represents a pseudo 1D problem,
which enables a straightforward comparison with the corresponding homogenized
model with the identical dimensions and boundary conditions at the macroscopic
level as the DNS.

\begin{figure}[ht]
    \centering
    \includegraphics[width=0.49\linewidth]{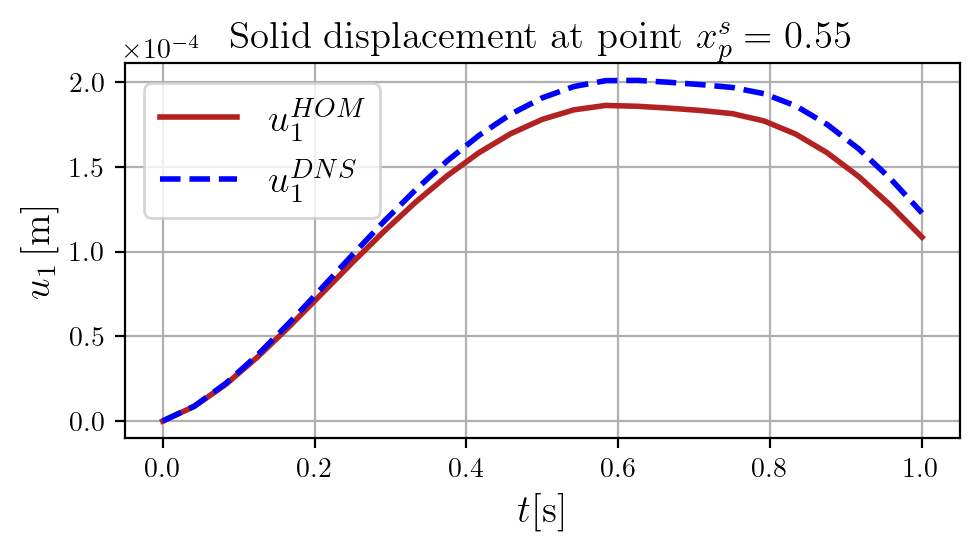}\hfil
    \includegraphics[width=0.49\linewidth]{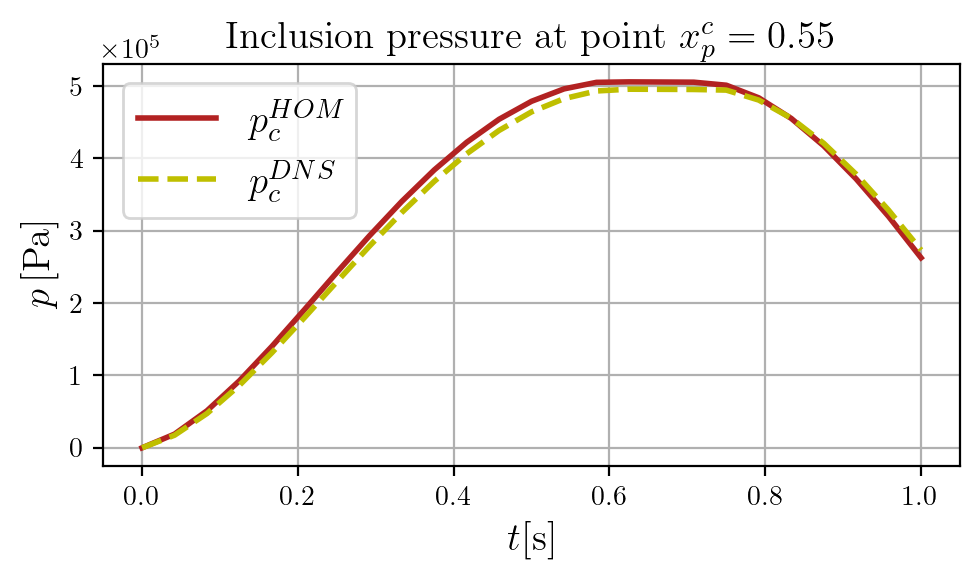}\\
    \includegraphics[width=0.49\linewidth]{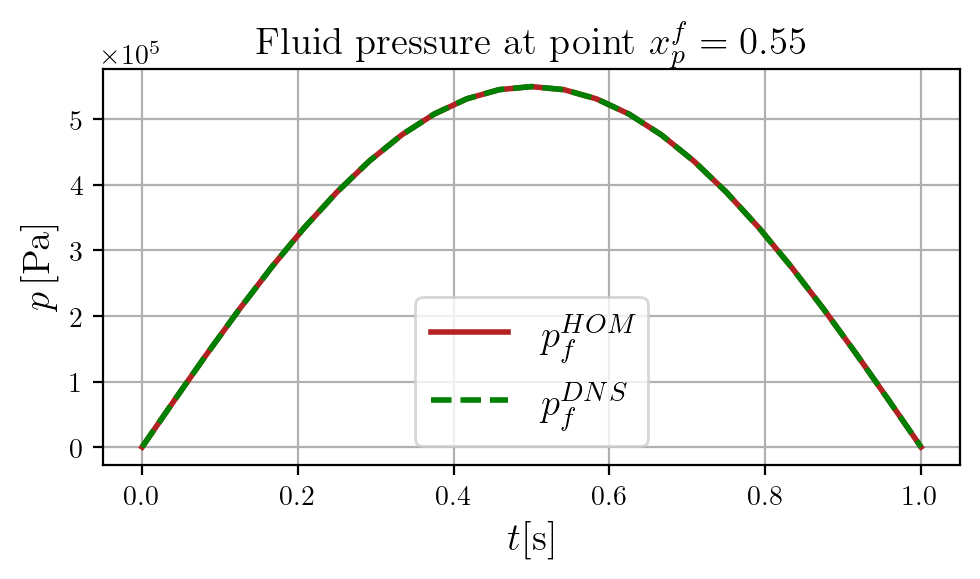}\hfil
    \includegraphics[width=0.49\linewidth]{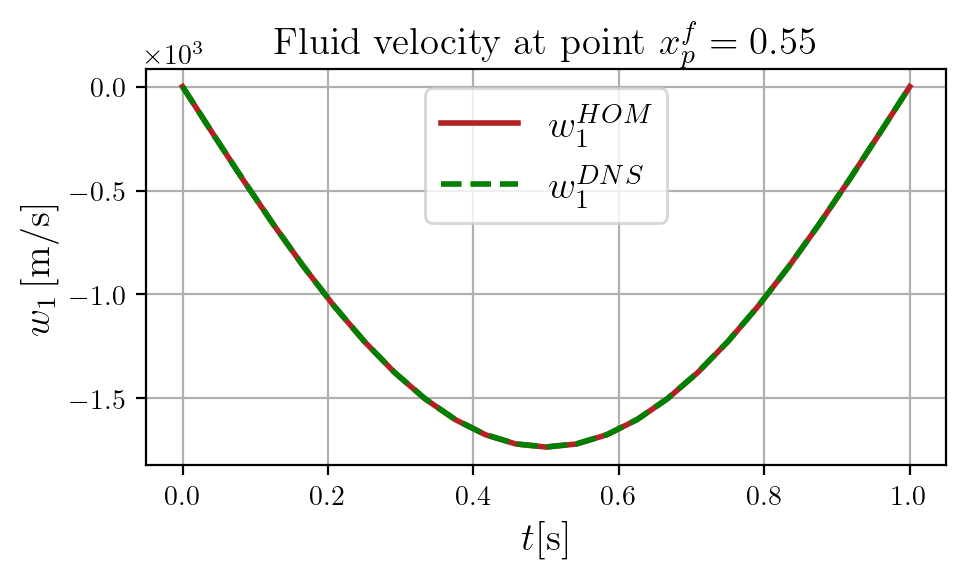}\\
    \caption{Comparison of the direct numerical simulation (DNS) and the homogenized model (HOM).
             Time evolution of displacements $u_1$, inclusion pressures $p_c$,
             channel pressures $p_f$, and channel velocity $w_1$ at point $x_p = 0.5$.}
    \label{fig:num-dns-time-evolution}
\end{figure}

Comparison of the results of DNS and the two-scale simulation is presented in
Figures \ref{fig:num-dns-line-probes} and \ref{fig:num-dns-time-evolution}. The
results of the homogenized model are denoted by the superscript $^{HOM}$ and
are were obtained by reconstructing the fields at the micro-level using
expansions \ref{eq-FS14a} for a finite value of $\veps_0$.
Fig.~\ref{fig:num-dns-line-probes} shows the spatial distribution of the
solution at a fixed time $t=0.5$\,s along the probe lines $l_p^s$, $l_p^f$, and
$l_p^c$ for the DNS and and along the probe line $l_p$ in the case of the
homogenized model. The time histories of the quantities at the selected point
$x_p$ of the line probes are shown in Fig.~\ref{fig:num-dns-time-evolution}.

Both simulations exhibit very good agreement, with only minor discrepancies
observed near the right boundary of the test sample, particularly in the
displacements $u_1$ and inclusion pressures $p_i$. These variations can be
explained by different model assumptions: while the homogenized model enforces
periodicity at the microscopic level, such periodic boundary conditions are not
imposed in the direct numerical simulation at the left and right sample boundaries.

\subsection{\newSecE{Steady state test -- stability}}\label{sec:num-steady}
In the present paper, besides the effects of deforming configuration being respected by the homogenized ``poro-piezo-hypoelastic'' material model defined in terms of the deformation dependent incremental homogenized coefficients, another nonlinearity due to the valves is included. This latter nonlinearity in combination with the dissipation due to the flow naturally rises questions about the uniqueness of the steady state solutions and their stability. These theoretical questions are beyond the scope of this publication, however, it is of interest to demonstrate a kind of stability of the developed computational algorithm to attain a steady state $\sbb(t) = \bar\sbb$. In particular, we are interested to show that, if the load and the BCs attain  constant spatial distributions for times $t > t_L > 0$ with a given $t_L$, there exists $t_s > t_L$, such that $\|\sbb(t) - \sbb(t_s)\| < \delta$ for any $t>t_s$ with a small $\delta> 0$. 
A rigorous proof may be difficult to obtain, but our numerical simulations confirm convergence to such steady states. Since the model involves the dissipation due to the viscous flow, it is not clear, if the steady state can be determined uniquely as a solution of the stationary problem (time independent), \ie without solving a nonstationary problem to achieve a steady state. Intuitively, this depends on the ejection pressure gauge $\Delta P$; when  $\Delta P > 0$, there is the uncertainty interval for the inclusion pressure $p_c(x) \in [p_f(x), [p_f(x) + \Delta P]$ for a.a. $x \in \Om$.

\begin{figure}[t]
    \centering
\begin{tabular}{cc}
    \includegraphics[width=0.49\linewidth]{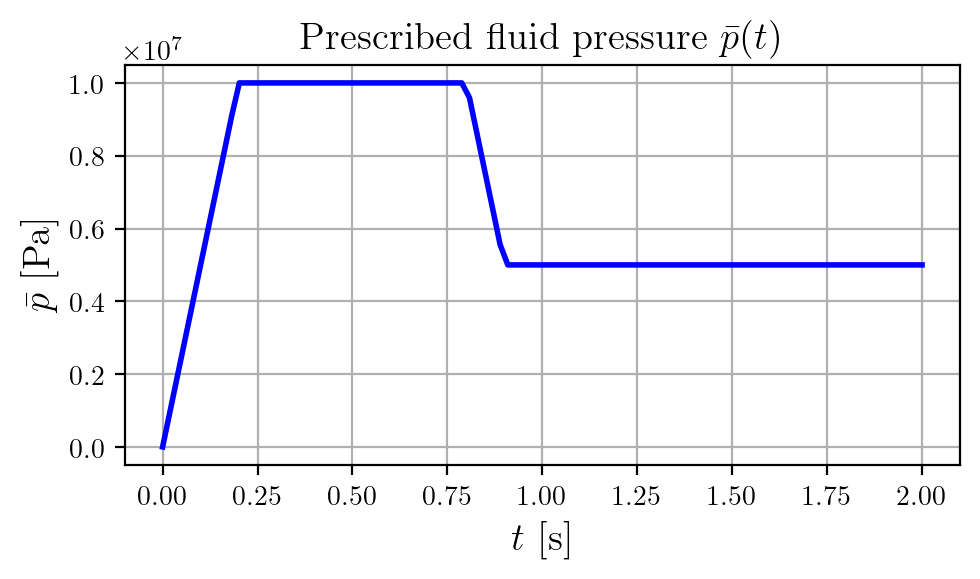}
    &
   \includegraphics[width=0.49\linewidth]{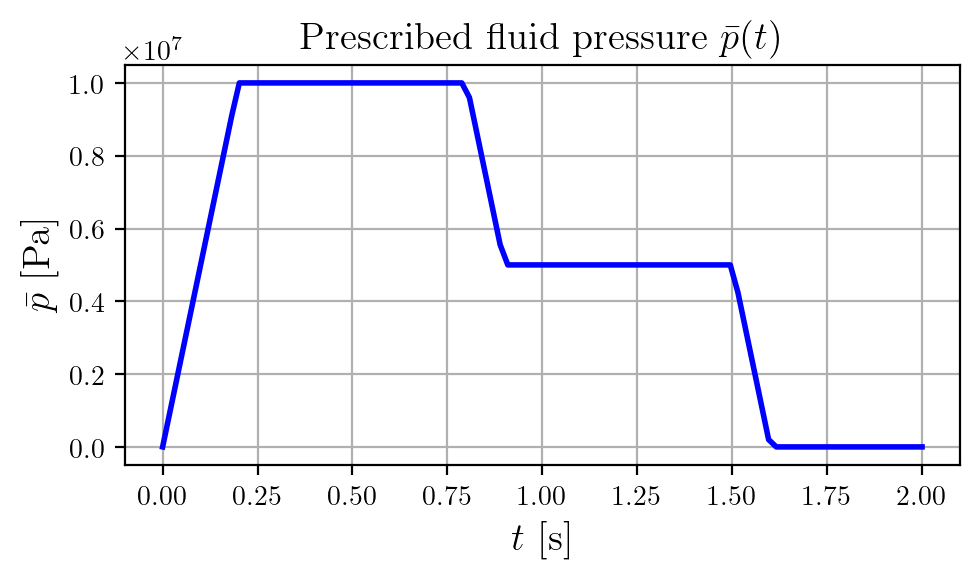}
\end{tabular}    
    \caption{Loading for the inflation tests in Fig.~\ref{fig-ss1} (left) and Fig.~\ref{fig-ss2} (right). Prescribed the channel pressure $p_f(\hat x=L,t) = \bar p(t)$ on the ``right'' end.}
\label{fig-ss0}
\end{figure}

\begin{figure}[t]
    \centering

\begin{tabular}{cc}
    \includegraphics[width=0.49\linewidth]{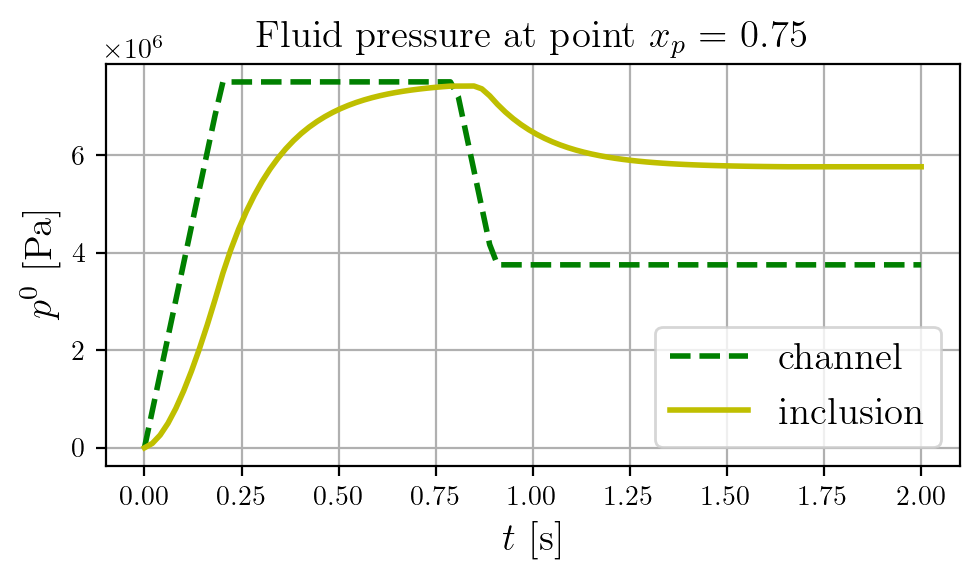}
    &
    \includegraphics[width=0.49\linewidth]{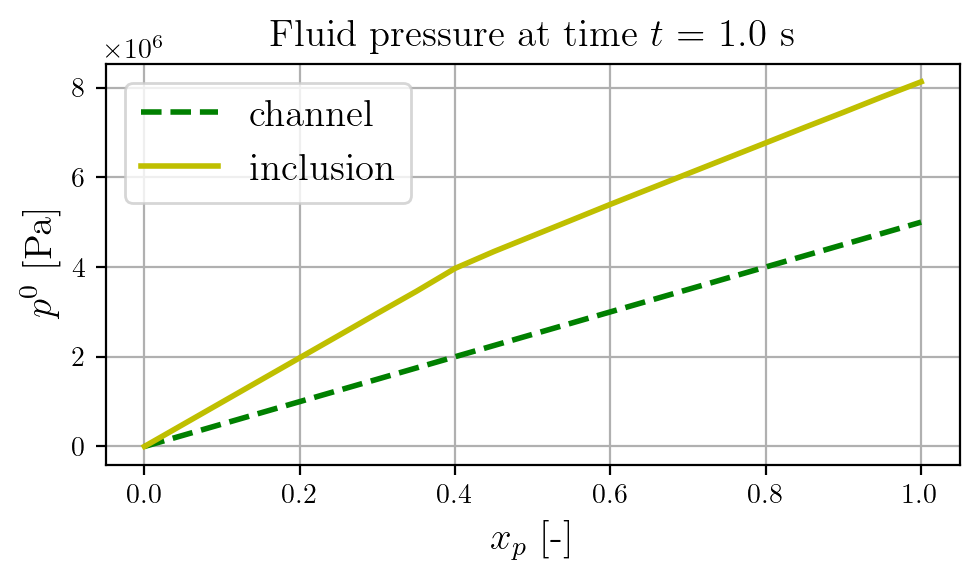}\\
    \includegraphics[width=0.49\linewidth]{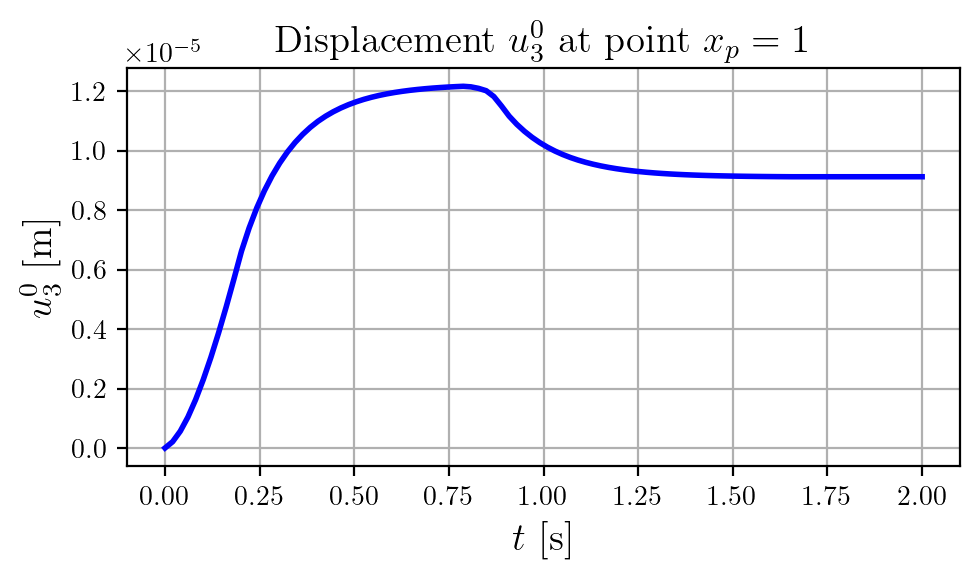}
    &
    \includegraphics[width=0.49\linewidth]{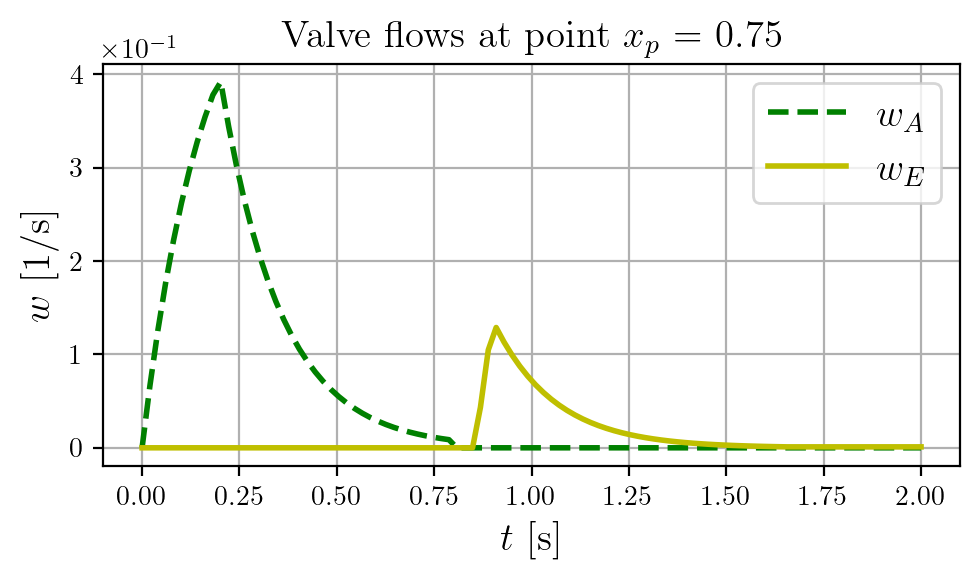}
\end{tabular}    
    \caption{Loading by the inflation pressure $\bar p$ according to in Fig.~\ref{fig-ss0} (left), testing the convergence to the steady state. Time evolution of $p_f(t,\hat x)$ and $p_c(t,\hat x)$ at $\hat x = 0.75$ (top-left); spatial distribution of $p_f(\hat t,x)$ for $\hat t = 1$ s (top-right);
    time evolution of displacement $u_3(t,L)$ at the sample end (bottom-left); flow through the valves $w_A(t,\hat x)$ and $w_E(t,\hat x)$ (bottom-right).}
    \label{fig-ss1}
\end{figure}

Behaviour of the inflatable structure in response to step-change in the channel pressure $p_f(\hat x,t)$ at $\hat x = L$, see Fig.~\ref{fig-ss0}, is demonstrated in two examples.
In Fig.~\ref{fig-ss1}, we illustrate the convergence to a steady state and its stability using the example presented in Section~\ref{sec:problem-settings} without the elastic part $\Om_e$, \ie $h=a$. The BCs are the same except of the pressure loading on the ``inlet'':  on $\Gamma_R^p$, the porous structure is loaded by a ``step-wise ramp-and-hold'' shaped pressure function $p_f(t,L) = \bar p f_R(t)$, where $\bar p = 10$ MPa and 
the ``double ramp'' function is $f_R(t) = f_\wedge(t) - 0.5 f_\vee(t)H_+(t - \hat t_2)$ with
$f_\wedge(t) = t/ \hat t_1 - (t - \hat t_1) H_+(t - \hat t_1)$,
and $f_\vee(t) = (t - \hat t_2)/ (\hat t_3-\hat t_2) - (t - \hat t_3) H_+(t - \hat t_3)$, see Fig.~\ref{fig-ss0} (left), where $\hat t_1=0.1$ s, $\hat t_2=0.4$ s,, and  $\hat t_3=0.45$ s.

\begin{figure}[t]
    \centering
    \includegraphics[width=0.49\linewidth]{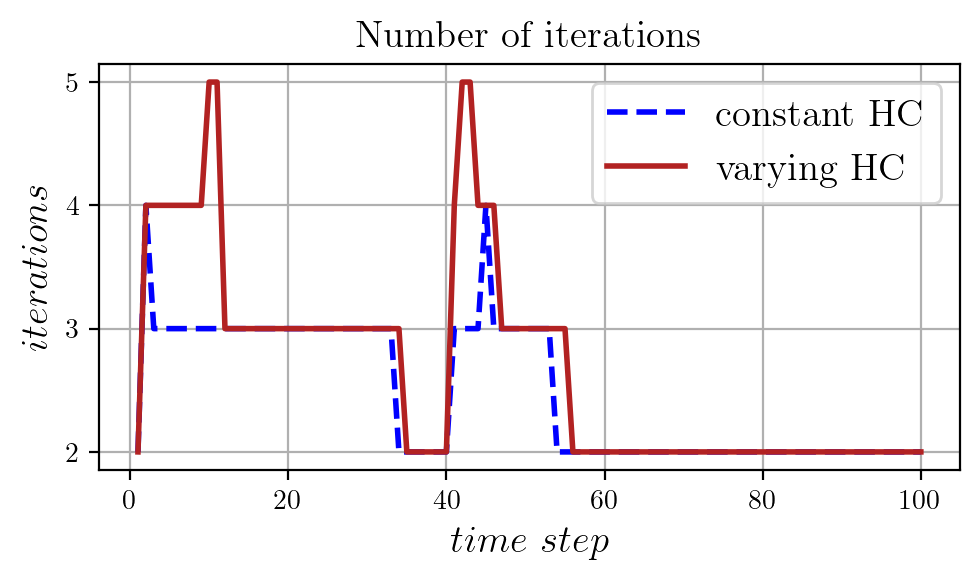}
 \caption{Number of the Newton-Raphson iterations in computing time steps (load by the inflation test according to Fig.~\ref{fig-ss0}, left.}
\label{fig-conv}
\end{figure}

As can be seen in Fig.~\ref{fig-ss1}, the channel pressure $p_f$ evolves almost synchronously with the time-dependent inlet pressure $p_f(t,L)$, whereas pressure in the inclusions asymptotically (in time) attains the channel pressure. The steady state characterized by $p_f(x,t) = p_c(x,t)$ seems to be achieved for $t\rightarrow +\infty$. The ejection valves are closed until $\hat t_2$, since $p_f$ monotonously increases everywhere in the specimen. Then they open (for some positions $x$ only, depending on the local pressure drop), releasing the fluid from the inclusions until a steady state is achieved. Note that, if the pressure drop between times $\hat t_2$ and $\hat t_3$ is smaller than $\Delta P$, the ejection valves remain close all the time, such that another steady state is achieved instantaneously for $t > \hat t_3$. In Fig.~\ref{fig-conv}, the number iterations within each time step is reported showing a good performance of the Newton-Raphson iterations with the consistent incremental homogenized coefficients.

\begin{figure}[t]
    \centering
  
\begin{tabular}{cc}
    \includegraphics[width=0.49\linewidth]{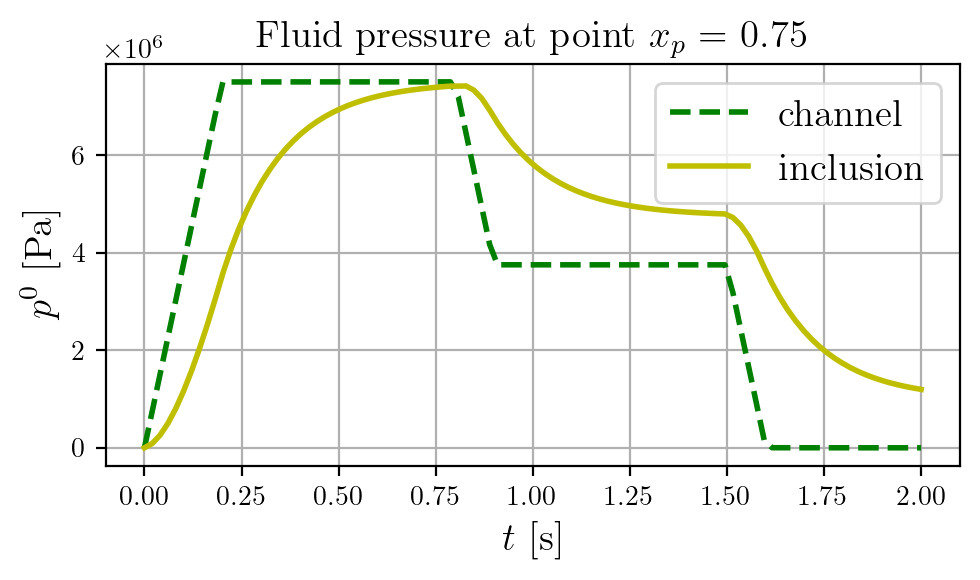}
    &
    \includegraphics[width=0.49\linewidth]{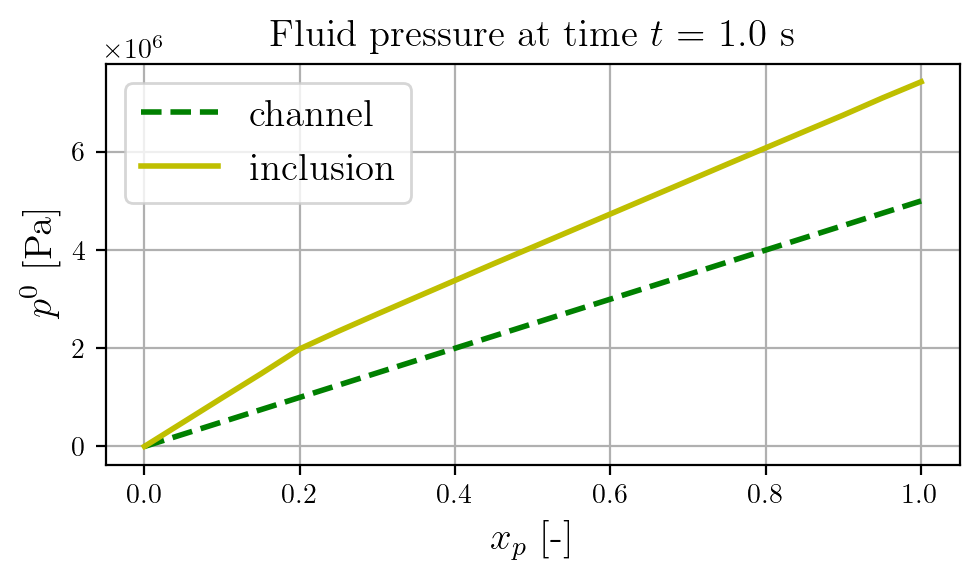}\\
    \includegraphics[width=0.49\linewidth]{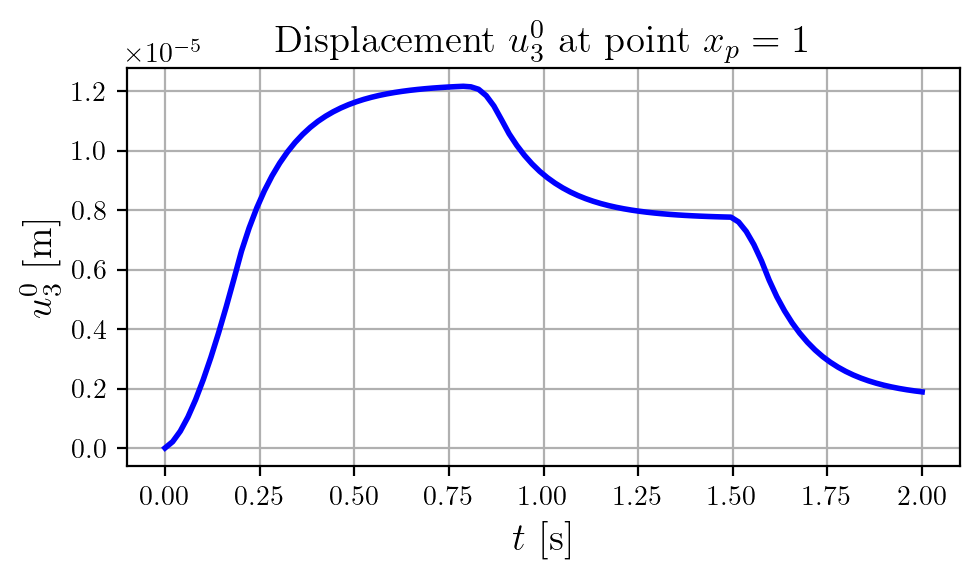}
    &
    \includegraphics[width=0.49\linewidth]{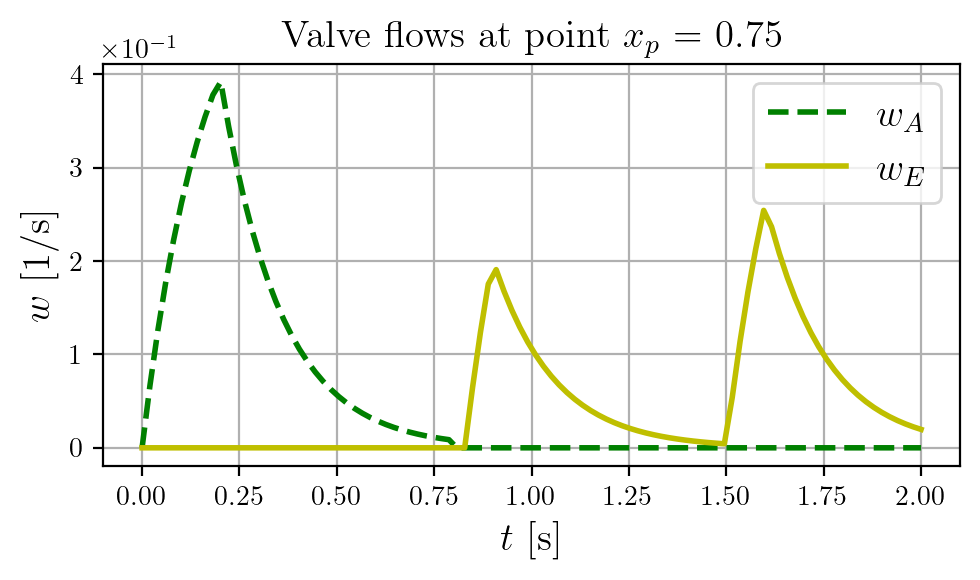}
\end{tabular}    
    \caption{The step-wise increasing -- decreasing the inflation pressure $\bar p$, for loading in Fig.~\ref{fig-ss0} (right): time evolution of $p_f(t,\hat x)$ and $p_c(t,\hat x)$ at $\hat x = 0.75$ (top-left); spatial distribution of $p_f(\hat t,x)$ for $\hat t = 1$ s (top-right);
    time evolution of displacement $u_3(t,L)$ at the sample end (bottom-left); flow through the valves $w_A(t,\hat x)$ and $w_E(t,\hat x)$ (bottom-right).}    
    \label{fig-ss2}
\end{figure}

In Fig.~\ref{fig-ss2}, a more complex loading is considered, see \ref{fig-ss0}, right. Note that the decrease in the inflation pressure $\bar p(t)$ is almost instantaneously reflected by $p_f(t,x)$ which is linear function of $x$. The inclusion are drained through the E-valve (depending on the location $x$) with a limited speed, which leads to the slow adjustment on the displacement. This shows a kind of viscoelasticity provided by the porous metamaterial, namely due to the valve limited permeability and the skeleton elasticity.


%% file: aux-Appx-inflatable-R.tex

\begin{appendices}

  \section{\newSecE{Direct numerical simulations (DNS) -- incremental formulation}}\label{sec-DNS}
  We consider the weak formulation of the porous medium inflation, as introduced in Section~\ref{sec-weak-form}, see \eq{eq-5a}-\eq{eq-wf3}, whereby, in wat follows, the dependence on the heterogenity scale $\veps$ is dropped from the notation for the sake of brevity. The solutions are computed using an incremental formulation obtained by virtue of the linearization of the valve-related nonlinearity, while the reference configurations is fixed. For a given load and BCs at time $t$, the response represented by the triplet $\ssb(t,\cdot) := (\ub(t,\cdot), \wb(t,\cdot),p_f(t,\cdot), p_c(t,\cdot)) \in \Ucalbf(\Om_{s}) \times \Wcalbf_0(\Om_f,\Gamma_\fsi)\times\Qcal_f(\Om_f)\times\Qcal_c(\Om_c)$ is computed by virtue of increments $\dlt\ssb := (\dlt\ub,\dlt\wb,\dlt p_f, \dlt_c)$, 
  \begin{equation}\label{eq-fsi3s}
\ssb^{i+1}(t,\cdot) := \ssb^i(t,\cdot) + \dlt\ssb(t,\cdot)\;,
  \end{equation}
  where $i$ labels the Newton-Raphson iterations.
  The nonlinear effects associated with equilibrium attained in the deformed configuration are not considered in the DNS based validation study reported in Section~\ref{sec:numex_dns}. However, this issue was treated in \cite{Rohan-Lukes-CaS2024}, in the context of the homogenized peristaltic-driven flow in poro-piezoelectric medium. Therein,  the FSI boundary term which appears also in the \lhs expression in \eq{eq-5a}, has been treated to obtain a ``symmetrized'' formulation of the equilibrium in the solid and in the fluid subdomains. The same approach is adopted here, leading to \eq{eq-fsi3}.

  Using the perturbations \eq{eq-fsi3s} (increments within one time level) and taking into account a finite number $N_c$ of inclusions, \ie $\Om_c = \sum_{k=1}^{N_c}$, over which $p_c$ and $q_c$ are constant, 
problem \eq{eq-5a}-\eq{eq-wf3} yields system of equations for $\dlt\ssb(t,\cdot)$,
\begin{equation}\label{eq-fsi3}
  \begin{split}
    \int_{\Om_s^\Rnil}[\Aop\eeb{\dlt\ub}]:\eb_\Rnil(\tilde \vb) \dV + \sum_{k=1}^{N_c}\int_{\pd \Om_{c,k}^\Rnil} (p_c + \dlt p_{c,k})\nb^\sx\cdot\vb \dS  &     \\
    + \int_{\Om_f^\Rnil}\left(\mu \eeb{\wb + \dlt\wb + \dot\ub + \dlt\dot\ub} - (p_f +\dlt p_f) \Ib\right):\eb_\Rnil(\tilde \vb) \dV & \\+ \int_{\Om^\Rnil_s}\sigmabf^\Rnil_s(\ub,\vphi):\eb_\Rnil(\tilde \vb)  \dV
   & =  \int_{\pd_\sigma \Om^\Rnil_f} (\bar p + \dlt \bar p)\nb\cdot \tilde \vb \dS\;,\\
    \int_{\Om^\Rnil_f}\left(\mu \eeb{\wb + \dlt\wb + \dot\ub + \dlt\dot\ub} - (p +\dlt p) \Ib\right):\eb_\Rnil(\vtheta)  \dV& =  \int_{\pd_\sigma \Om^\Rnil_f} (\bar p + \dlt \bar p) \nb \cdot \vtheta\dS\;,\\
  \end{split}
\end{equation}
\begin{equation}\label{eq-fsi3a}
\begin{split}
& \int_{\Om_f^\Rnil}\gamma (\dot p + \dlt \dot p) q   + \int_{\Om_f^\Rnil}[q \nabla\cdot(\dot{\wtilde{\ub}} + \dlt \dot{\wtilde{\ub}}) - (\wb + \dlt \wb)\cdot\nabla q ] \\
  & = \int_{\Gamma_\fsi^\Rnil} q \kappa_A \left(\posPart{p_f - p_{c,k} +\dlt p_f - \dlt p_{c,k}}\right) \hat\delta_A \dS \\
  & - \int_{\Gamma_\fsi^\Rnil} q \kappa_E \left(\posPart{p_{c,k} - p_f  - \Dlt P_E +\dlt p_{c,k} - \dlt p_f}\right)\hat\delta_E \dS  + \int_{\Gamma_\inout} \wb\cdot\nb^\fx q \dS\;,
\end{split}
\end{equation}
\begin{equation}\label{eq-fsi3b}
\begin{split}  
  &\int_{\pd\Om_{c,k}^\Rnil} (\dot{\wtilde{\ub}} + \dlt\dot{\wtilde{\ub}})\cdot \nb \dS  + |\Om_{c,k}|\gamma (\dot p_{c,k} + \dlt \dot p_{c,k}) 
  - \int_{\pd\Om_{c,k}^\Rnil} \kappa_A\bar\dlt_{A,k} \posPart{\hat p_f + \dlt \hat p_f - p_{c,k} -  \dlt p_{c,k}} \dS\\
  & - \int_{\pd\Om_{c,k}^\Rnil}  \kappa_E\bar\dlt_{E,k} \posPart{p_{c,k} + \dlt p_{c,k} - \hat p_f - \dlt \hat p_f - \Dlt P_E} \dS = 0\;, \quad k = 1,\dots,N_c\;,
\end{split}
\end{equation}
for any $\tilde \vb \in \Ucalbf_0(\Om_s^\Rnil)$, $\vtheta \in \Wcalbf_0(\Om_f^\Rnil)$, $ q \in \Qcal_{0f}(\Om_f^\Rnil)$. 
  
Discretization in time is explained in Section~\ref{sec-iter-nonlin} for the homogenized model. In analogy, it is applied to discretize in time equations \eq{eq-fsi3}-\eq{eq-fsi3b}.
We consider time levels $t_0 < t_1 < t_2,\dots$, such that backward approximations in time are used to define $\dot \ssb(t_k,x) \approx (\ssb^k - \ssb^{k-1})/\Dlt t$. Within each time level, the ``loads'' are presented by $\dlt \bar p$ 

\section{Sensitivity of the homogenized coefficients (HC)}\label{sec-SA}
  

\subsection{Auxiliary results on the shape sensitivity}\label{sec-SA-shape}
To differentiate the integral forms employed in the homogenized model, we shall use some particular results of the domain method of the shape sensitivity analysis which relies on the convection velocity field $\vec\Vcal$ defined in entire domain $Y$; see \eg. \cite{Rohan-Lukes-CaS2024} where a similar sensitivity analysis was reported for Biot coefficients of a piezoelectric porous material. For the context of the shape sensitivity analysis see \eg \cite{haslinger_makinen_2003:_introduction_shape_optimization,haug_choi_komkov_1986:_design_sensitivity_analysis}. Here by $\dltsh$ the ``partial'' shape derivative is indicated. The following auxiliary results are employed,
\begin{equation}\label{eq-sa5}
  \begin{split}
    \dltsh|Y_d| & = \int_{Y_d} \nabla_y\cdot\vec\Vcal = \int_{\pd Y_d} \vec\Vcal\cdot\nb\;,\\
    \dltsh \phi_c & = |Y|^{-1} (\dltsh|Y_c| - \phi_c\dltsh|Y|)\;,\\
  \dltsh\intY_{\Gamma_c} \vb\cdot\nb^\cx & =  - \dltsh\intY_{Y_m}\nabla_y\cdot\vb \\
  & = \intY_{Y_m} \pd_i^y\Vcal_k \pd_k^y v_i
  - \intY_{Y_m}\nabla_y\cdot\vb \nabla_y\cdot\vec\Vcal + |Y|^{-1}\dltsh|Y|\intY_{Y_m}\nabla_y\cdot\vb
  \;,\\
  \dltsh \Pi_k^{ij} & = \dlt_{ik}\Vcal_j\;.
\end{split}
\end{equation}
Note that $\dltsh|Y| = 0$ when $\vec\Vcal$ is $Y$-periodic.

For the bilinear forms involved in the local problems, the following sensitivity expression holds:
 \begin{equation}\label{eq-sa6}
\begin{split}
\dltsh \aYs{\ub}{\vb} &= \intY_{Y_m\cup Y_*} 
D_{irks}\left(\delta_{rj}\delta_{sl} \nabla_y\cdot \vec\Vcal - \delta_{jr}\pd_s^y \Vcal_l - \delta_{ls}\pd_r^y \Vcal_j\right)
e_{kl}^y(\ub) e_{ij}^y(\vb)\\
& \quad - \aYs{\ub}{\vb}\intY_{Y} \nabla_y\cdot \vec\Vcal\;.
\end{split}
 \end{equation}

\subsection{Sensitivity formulae for the poroelastic coefficients}\label{sec-SA-PE}
For the inflatable structures with the local characteristic problems featured
by the two porosities separated by the valves, the sensitivity analysis of the
HC \eq{eq-HC1} leads to the following expressions. We employ the abbreviation
$\Xibf^{ij} = \omegabf^{ij}+\Pibf^{ij}$ and $P, Q = f, c$.
\paragraph{Drained elasticity}
\begin{equation}\label{eq-HC4a}
\begin{split}
 \dlt A_{klij}^H &  = 
 \dltsh\aYs{\Xibf^{ij}}{\Xibf^{kl}} + \aYs{\Xibf^{ij}}{\dltsh\Pibf^{kl}} + \aYs{\dltsh\Pibf^{ij}}{\Xibf^{kl}}
\end{split}
\end{equation}

\paragraph{Biot coupling coefficients}
\begin{equation}\label{eq-HC1a}
  \begin{split}
 B_{ij}^P & = C_{ij}^P + \phi_P\dlt_{ij}\;,\quad \mbox{ with } C_{ij}^P = \intY_{\Gamma_P} \omegabf^{ij} \cdot \nb^\px \dSy\;,\\
 \dlt C_{ij}^P & = \dlt_\tau \left( \intY_{Y_P} \nabla_y \cdot \wtilde\omegabf^{ij}\right) +
 \dlt_\tau \aYs{\omegabf^P}{\omegabf^{ij}+\Pibf^{ij}} + \aYs{\omegabf^P}{\dlt\Pibf^{ij}}\;,\\
\dlt B_{ij}^P & = \dlt C_{ij}^P + \dlt_\tau \phi_P\dlt_{ij}\;.
   \end{split}
\end{equation} 

\paragraph{Biot modulus}
\begin{equation}\label{eq-HC2}
  \begin{split}
\dlt M^{PQ} & = -\dlt_\tau  \aYs{\omegabf^P}{\omegabf^Q} -  \dlt_\tau\intY_{\Gamma_P} \omegabf^Q \cdot \nb^\px \dSy -  \dlt_\tau\intY_{\Gamma_Q} \omegabf^P \cdot \nb^\qx \dSy  + \dlt_{PQ}\gamma \dlt_\tau\phi_P
   \end{split}
\end{equation} 
where
\begin{equation}\label{eq-HC3}
  \begin{split}
    \dlt_\tau\intY_{\Gamma_P} \omegabf \cdot \nb^\px \dSy = -|Y|^{-1}\intY_{\Gamma_P} \omegabf \cdot \nb^\px \dSy +  \intY_{Y_P} \left(\nabla_y\cdot \vec\Vcal\nabla_y\cdot \omegabf - \Ib:\nabla_y\vec\Vcal \nabla_y\omegabf \right)\;,\\
   \end{split}
\end{equation}

  \subsection{Sensitivity formula for the effective permeability}\label{sec-SA-K}
We are interested in the influence of variation of the shape of the
interface $\Gamma$ on the homogenized permeability $K_{ij}$ defined in
\eq{eq-S7a}, involving the influence of the semipermeable interface. The differentiation of $K_{ij}$ yields
\begin{equation}\label{eq-S7}
\begin{split}
\delta K_{ij} = \delta_\tau\intY_{Y_f}\psi_j^i + \intY_{Y_f}\delta\psi_j^i\;.
\end{split}
\end{equation}
To eliminate the dependence on $\delta\psi_j^i$ in the last integral
of the above expression, we differentiate \eq{eq-S3}$_1$ which yields
\begin{equation}\label{eq-S8a}
\begin{split}
\aYf{\dlt{\psibf}^i}{\vb}
+\intY_{Y_f}\nabla_y\dlt\pi^i \cdot\vb + \delta_\tau\left(\aYf{{\psibf}^i}{\vb}+
\intY_{Y_f}\nabla_y\pi^i \cdot\vb\right) = \delta_\tau\intY_{Y_f}v_i\;.
\end{split}
\end{equation}
Then we need \eq{eq-S3}$_1$ evaluated for $k= j$ with substituted $\vb = \delta \psibf^i$, 
and differentiate \eq{eq-S3}$_2$, then substitute there $q = \pi^j$, so that we obtain
\begin{equation}\label{eq-S8b}
\begin{split}
& \intY_{Y_f} \delta\psi_j^i = \intY_{Y_f}\left(\nabla_y\pi^j\cdot \dlt {\psibf}^i -\nabla_y\dlt\pi^i\cdot {\psibf}^j \right) \\
& \quad +\delta_\tau\intY_{Y_f}\psi_i^j - \dlt_\tau \aYf{{\psibf}^i}{{\psibf}^j} - \dlt_\tau\intY_{Y_f} \nabla_y \pi^i \cdot \psibf^j
\;, \\
& \dlt_\tau\intY_{Y_f} \nabla_y \pi^j \cdot \psibf^i - \dlt_\tau\ipGG{\hat\kappa\jumpG{\dlt\pi^i}}{\jumpG{\pi^j}} =  -\intY_{Y_f}\nabla_y\pi^j\cdot \dlt {\psibf}^i + \intY_{Y_f}\nabla_y\dlt\pi^i\cdot {\psibf}^j\;,
\end{split}
\end{equation}
where also \eq{eq-S3}$_2$ was employed with $q = \dlt \pi^i$ to get the 2nd equality in \eq{eq-S8b}.
We proceed by substitution $\vb = \psibf^j$ in \eq{eq-S8a} and use both the equalities in \eq{eq-S8b}, thus, we get
\begin{equation}\label{eq-S9}
\begin{split}
\intY_{Y_f} \delta\psi_j^i 
& = \delta_\tau \left(
\intY_{Y_f} \psi_i^j - \aYf{{\psibf}^j}{{\psibf}^i}
-  \intY_{Y_f} \nabla_y\pi^i  \cdot \psibf^j
\right) - \intY_{Y_f} \nabla_y\dlt\pi^i  \cdot \psibf^j +  \intY_{Y_f} \nabla_y\pi^j  \cdot \dlt\psibf^i
\\
& = \delta_\tau \left(
\intY_{Y_f} \psi_i^j - \aYf{{\psibf}^j}{{\psibf}^i} -\intY_{Y_f}  \nabla_y\pi^i \cdot \psibf^j - \intY_{Y_f} \nabla_y\pi^j  \cdot \psibf^i\right)
\\
&  +\delta_\tau  \ipGG{\hat\kappa\jumpG{\pi^i}}{\jumpG{\pi^j}}\;.
\end{split}
\end{equation}
Now we can substitute into \eq{eq-S7} and integrate by parts, 
which yields
\begin{equation}\label{eq-S10}
\begin{split}
\delta K_{ij} & = \delta_\tau  \intY_{Y_f} \left(  \psi_i^j + \psi_j^i
- \aYf{{\psibf}^j}{{\psibf}^i} + \ipGG{\hat\kappa\jumpG{\pi^i}}{\jumpG{\pi^j}}\right) \\
& + \delta_\tau  \intY_{Y_f}\left(
\pi^i \nabla_y \cdot \psibf^j +\pi^j \nabla_y \cdot \psibf^i\right)\;.
\end{split}
\end{equation}
To compute the shape derivative $\delta_\tau$, we shall employ the following expressions:
\begin{equation}\label{eq-S11a}
\begin{split}
\delta_\tau \intY_{Y_f}\psi_i^j & =
\intY_{Y_f}\psi_i^j \left(\nabla\cdot\vec\Vcal- \intY_Y\nabla_y\cdot \vec\Vcal\right)\;,\\
\delta_\tau \intY_{Y_f} \nabla_y\cdot \psibf^j
& = \intY_{Y_f} \nabla_y\cdot \psibf^j(\nabla_y \cdot \Vcal - \intY_Y\nabla_y\cdot \vec\Vcal)- \intY_{Y_f} \pd_k^y\Vcal_r \pd_r^y\psi_k^j\;,\\
\dltsh \aYf{\ub}{\vb} &= \intY_{Y_f} 
D_{irks}^f\left(\delta_{rj}\delta_{sl} \nabla_y\cdot \vec\Vcal - \delta_{jr}\pd_s^y \Vcal_l - \delta_{ls}\pd_r^y \Vcal_j\right)
e_{kl}^y(\ub) e_{ij}^y(\vb)\\
& \quad - \aYf{\ub}{\vb}\intY_{Y} \nabla_y\cdot \vec\Vcal +
\dltsh\intY_{\Gamma_f} (\hat\kappa_f)^{-1} \nb\otimes\nb : \hat\thetabf \otimes\hat\wb\;.
\end{split}
\end{equation}
The surface divergence operator $\nabla_\gamma \cdot \Vcal$ is involved in the shape derivative of the interface integrals,\chE{
\begin{equation}\label{eq-S11b}
\begin{split}
  \dltsh\intY_{\Gamma_f} (\hat\kappa_f)^{-1} \nb\otimes\nb : \hat\thetabf \otimes\hat\wb  & =
    \intY_{\Gamma_f} (\hat\kappa_f)^{-1} \nb\otimes\nb : \hat\thetabf \otimes\hat\wb (\nabla_\gamma \cdot \Vcal  -  \intY_Y\nabla_y\cdot \vec\Vcal)\\
    & \quad - \intY_{\Gamma_f} (\hat\kappa_f)^{-2} \delta_\tau \hat\kappa (\nb\otimes\nb : \hat\thetabf \otimes\hat\wb)
\;,\\
\delta_\tau \ipGG{\hat\kappa\jumpG{\pi^i}}{\jumpG{\pi^j}} & = \intY_{\Gamma_f}\hat\kappa\jumpG{\pi^i}\jumpG{\pi^i} (\nabla_\gamma \cdot \Vcal -  \intY_Y\nabla_y\cdot \vec\Vcal) \\
& \quad + \ipGG{\delta_\tau \hat\kappa\jumpG{\pi^i}}{\jumpG{\pi^j}}\;,
\end{split}
\end{equation}
}
whereby $\nabla_\gamma \cdot \Vcal = \nabla \cdot (\Vcal - \nb^\Gamma (\nb^\Gamma\cdot\Vcal))$ with $\nb^\Gamma$ being a constant unit normal vector of a flat surface $\Gamma_f$.

Finally, using \eq{eq-S10} and \eq{eq-S11a} we get the following expression
\begin{equation}\label{eq-S11}
\begin{split}
\delta K_{ij} & = \intY_{Y_f}\left( \psi_i^j + \psi_j^i - (\Dop^f\eeby{\psibf^j}) : \eeby{\psibf^i} 
+ \pi^i \nabla_y \cdot \psibf^j +\pi^j \nabla_y \cdot \psibf^i\right) \times \\
& \quad\times \left(\nabla_y\cdot \vec\Vcal - \intY_Y\nabla_y\cdot \vec\Vcal\right) \\
& - \intY_{\Gamma_f}\hat\kappa\jumpG{\pi^i}\jumpG{\pi^i} (\nabla_\gamma \cdot \Vcal -  \intY_Y\nabla_y\cdot \vec\Vcal) + \ipGG{\delta_\tau \hat\kappa\jumpG{\pi^i}}{\jumpG{\pi^j}}
 \\
& \quad + \intY_{Y_f} \left( D_{prks}^f(\delta_{qr}\pd_s^y \Vcal_l + \delta_{ls}\pd_r^y \Vcal_q)e_{kl}^y( \psibf^j) e_{pq}^y(\psibf^i)
 - \pi^i \pd_k^y\Vcal_r \pd_r^y\psi_k^j -  \pi^j \pd_k^y\Vcal_r \pd_r^y\psi_k^i 
\right)\;.
\end{split}
\end{equation}
We shall disregard $\delta_\tau \hat\kappa = 0$.
%
Although $\delta K_{ij}$ depends formally on $\vec\Vcal$, this field can
be constructed arbitrarily in the interior of $Y_f$ without any
influence on $\delta K_{ij}$.


\section{Discretized subproblems}\label{app-discr}
The homogenized model leads to problem  \eq{eq-S25-fc} which is nonlinear due to the incorporated valves. Having applied the FE discretization in space, as explained in Section~\ref{sec:numex}, in the matrix form, the problem for displacement and pressure vectors, $\ubm$ and $\ul{\pbm} = [\pbm_f,\pbm_c]$, respectively, is represented by the following equation, 
\begin{equation}\label{eq-nlFE1}
  \begin{split}
    \Abm \ubm - \ul{\Bbm}^T \ul{\pbm} & = \fbm \;,\\
    \ul{\Bbm}(\ubm - \ubm^\tp) + \Dlt t \ull{\Kbm}\ul{\pbm} + \ull{\Mbm}(\ul{\pbm} -  \ul{\pbm}^\tp) +\Dlt t { \ull{\Nbm}\Dltb\circ\ul{\pbm}} & = \ul{\gbm}\;,
    \end{split}
\end{equation}
where
\begin{equation}\label{eq-nlFE2}
  \begin{split}
    \ull{\Nbm}\Dltb\circ\ul{\pbm} =
    \left [
  \begin{array}{ll}
    \Rbm & -\Sbm \\ -\Rbm & \Sbm
  \end{array} \right]
     \left [
  \begin{array}{l}
    { \Dlt_0^+ \ul{\pbm} }\\  {\Dlt_{-P}^+ (-\ul{\pbm})}
  \end{array} \right] 
 \end{split}
\end{equation}
and the time derivatives are approximated by the backward finite differences in time
   ${\dot\ubm} \approx (\ubm - \ubm^\tp)/\Dlt t$, and ${\dot{\ul{\pbm}}} \approx (\ul{\pbm} - \ul{\pbm}^\tp)/ \Dlt t$. The valve closing/opening depending on the pressure jump is represented by $\Dlt_Q^+\ul{\pbm} = \posPart{\pbm^f - \pbm^c - Q\onebm}$, involving the ``overessure'' control $Q$.
The Newton method is applied to solve \eq{eq-nlFE1}, thus, generating increments $\dlt \ubm$ and $\dlt \ul{\pbm}$ satisfying the linear equation
\begin{equation}\label{eq-nlFE4}
  \begin{split}
       \left [
      \begin{array}{ll}
        \Abm \;,& -\ul{\Bbm}^T \\
        \ul{\Bbm}\;, &       \Dlt t \ull{\Kbm} + \ull{\Mbm} + \Dlt t \ull{\Nbm} \ull{\Jbm_+^{fc}}
      \end{array} \right] \left[\begin{array}{l}
  \dlt \ubm \\ \dlt \ul{\pbm}
       \end{array} \right] = \\
         \left[\begin{array}{l}
   \fbm - \Abm \bara\ubm + \ul{\Bbm}^T \bara{\ul{\pbm}} \\ \Dlt t  \ul{\gbm}
     - \ull{\Nbm}\Dltb\circ\ul{\pbm} +  \ul{\Bbm}(\ubm^\tp - \bara\ubm) +  \ull{\Mbm}{\ul{\pbm}^\tp} -  (\Dlt t \ull{\Kbm} + \ull{\Mbm}) \bara{\ul{\pbm}}
       \end{array} \right] \;.
     \end{split}
\end{equation}
This yields the obvious update: {$\ubm:= \bara\ubm + \dlt\ubm$} and {$\ul{\pbm}:=\bara{\ul{\pbm}} + \dlt\ul{\pbm}$}.


\section{\newSecE{Survey of used notation / symbols}}\label{sec-Notation}
For the readers convenience, the notation employed in the paper is summarized according to the context, being listed in five tables, Tab.~\ref{tab-notation1},Tab.~\ref{tab-notation2}, Tab.~\ref{tab-notation3}, Tab.~\ref{tab-notation4}, and Tab.~\ref{tab-notation5}.
\begin{table}[ht]
  \caption{Geometrical objects -- domains, boundaries, interfaces} \label{tab-notation1},
    \begin{tabular}{|c|l|}
        \hline
        $\veps$ & scaling parameter \\
        $x$, $y$ & ``macroscopic'' and ``microscopic'' coordinates\\
        $\Tuf{\cdot}$  & the unfolding operator\\
        $\Om$ & open bounded domain in $\RR^3$\\
        $\pd\Om$ & Lipschitz boundary of $\Om$\\
        $\Om_s^\veps$, $\Om_f^\veps$, $\Om_c^\veps$ & solid matrix, fluid channel, and fluid inclusion domains\\
        $\Gamma_\fsi^\veps$ & solid-fluid (channel) interface: $\Gamma_\fsi^\veps  = \ol{\Om_s^\veps} \cap \ol{\Om_f^\veps}$\\
        $\Gamma_u^\veps$, $\Gamma_p^\veps$ & boundaries with prescribed displacements and pressures\\
        $\Gamma_\sigma^\veps$, $\Gamma_w^\veps$ & boundaries with prescribed surface tractions and fluid fluxes\\
        $\pd_\ext\Om_s^\veps$, $\pd_\ext\Om_f^\veps$ & exterior boundaries: $\pd_\ext\Om_s^\veps = \Gamma_u^\veps \cup \Gamma_\sigma^\veps$, $\pd_\ext\Om_f^\veps = \Gamma_p^\veps \cup \Gamma_w^\veps$\\
        $\bar\Gamma_E^\veps$, $\bar\Gamma_A^\veps$ & admission and ejection channels connecting $\Om^\veps_f$ and $\Om^\veps_c$\\
        $\Gamma_f^{k,\veps}$ & semi-permeable interfaces in fluid channels $\Om^\veps_f$\\
        $\Gamma_\inout$  & inlet and outlet surfaces of the connected porosity \\
        $Y$ & reference periodic cell\\
        $Y_s$, $Y_f$, $Y_c$ & non-overlapping solid, fluid, and inclusion subdomains of $Y$\\
        $\Gamma_\fsi$ & solid-fluid interface: $\Gamma_\fsi  = \ol{Y_s} \cap \ol{Y_f}$\\
        ${\Gamma_f^{k}}$ & semi-permeable interfaces in fluid channels $Y_f$\\
        $\Ical_c^\veps$, $\Ical_f^\veps$ & index sets, such that: $\Om_c^\veps = \bigcup_{\kkk\in\Ical_c^\veps}{\Om_c^{\kkk,\veps}}$, $\Om_f^\veps = \bigcup_{\kkk\in\Ical_f^\veps}{\Om_f^{\kkk,\veps}}$\\
        \hline
    \end{tabular}
    \color{black}
\end{table}

\begin{table}[ht]
  \caption{Derivatives, jumps}\label{tab-notation2}
    \begin{tabular}{|c|l|}
        \hline
        $\posPart{a}$ & positive part of $a$\\
        $\jump{q}^{k,\veps}$ & pressure jump on interface $\Gamma_f^{k,\veps}$: $\jump{q}^{k,\veps} = q|_{k+1} - q|_{k}$\\
        $\nabla_x = (\pd_i^x)$, $\nabla_y = (\pd_i^y)$ & differentiation \wrt coordinate $x$ and $y$\\
        $\eeb{\ub} = \frac{(\nabla\ub)^T + \nabla\ub}{2}$ & strain of vector function $\ub$\\
        $\dot\ub$ & material derivative of $\ub$\\
        $\dltsh$ &  shape derivative, see Section~\ref{sec-SA-shape}\\
        \hline
    \end{tabular}
    \color{black}
\end{table}

\begin{table}[ht]
  \caption{Variables, model coefficients and related parameters}\label{tab-notation3}
    \begin{tabular}{|c|l|}
        \hline
        $\ub^\veps$, ${\wtilde\ub}^\veps$ & displacement field in $\Om_s^\veps$, its smooth extension to $\Om_f^\veps$\\
        $p_c^\veps$, $p_f^\veps$ & fluid pressure fields in $\Om_c^\veps$ and $\Om_f^\veps$\\
        $\wb^\veps$ & seepage velocity in $\Om_f^\veps$\\
        $\nb^{[S]}$ & unit normal vector on a boundary, outward to domain ``S''\\
        $\nb^\Gamma$ & unit normal vector on the fluid-fluid interface $\Gamma_f^\veps$\\
        $\sigmabf^\veps$ & second-order stress tensor\\
        $\Ib$ & identity tensor: $\Ib = (I_{ij}) = \delta_{ij}$\\
        $\Aop^\veps=(A_{ijkl}^\veps)$ & elasticity fourth-order symmetric positive definite tensor\\
        $\Dop_f^\veps=(D_{ijkl}^{f,\veps})$ & fluid viscosity fourth-order symmetric positive definite tensor\\
        $\mu^\veps$, $\bar\mu$ & fluid dynamic viscosity: $\mu^\veps = \veps^2\bar\mu$\\
        $\gamma$ & fluid compressibility\\
        $\hat\delta_V^\veps$, $\bar\delta_V^\veps$, $V=A,E$ & distribution function to define flow through valves \\
                $\kappa_A^\veps$, $\kappa_E^\veps$ & permeabilities of channels $\bar\Gamma_E^\veps$, $\bar\Gamma_A^\veps$\\
        $\bar\kappa_A$, $\bar\kappa_E$ & rescaled permeabilities: $\kappa_A^\veps = \veps \bar\kappa_A$, $\kappa_E^\veps = \veps \bar\kappa_E$\\
        $\gamma_A$, $\gamma_E$ & valve opening/closing switch, see \eq{eq-sts6i-V1}\\
        $\Dlt P_E$ & ejection pressure threshold \\
        $\bar w_A^\veps$, $\bar w_A^\veps$ & flows in channels $\bar\Gamma_E^\veps$, $\bar\Gamma_A^\veps$\\
        $\hat\kappa_f^\veps$, $\hat\kappa_f$ & permeability of the semi-permeable fluid-fluid interface $\Gamma_f^{k, \veps}$\\   
        $\fb^{s,\veps}$, $\fb^{s,\veps}$ & volume forces in solid and fluid domains\\
        $\bb^{s,\veps}$ & surface tractions applied on solid surface $\pd_\sigma\Om^\veps_s$\\
        $\ub_\pd$, $ p_\pd$ & prescribed displacements on $\Gamma^\veps_u$ and pressures on $\Gamma^\veps_p$\\
\hline
    \end{tabular}
    \color{black}
\end{table}

\begin{table}[ht]
  \caption{Specific notation of the homogenized two-scale model}\label{tab-notation4}
    \begin{tabular}{|c|l|}
      \hline
        $\ub^0$, $p^0_f$, $p^0_c$, $\hat\wb$ & macroscopic displacement, pressure, and velocity fields\\
        $\ub^1$, $p_f^1$ & displacement and fluid pressure microscopic fluctuations\\
        $\vb^\veps$, $q_f^\veps$, $q_c^\veps$, $\thetabf^\veps$ & displacement, pressure, and velocity test functions\\
        $\vb^0$, $q_f^0$, $q_c^0$ & macroscopic displacement and pressure test functions\\
      $\vb^1$, $q_f^1$, $\hat\thetabf$& microscopic test functions\\
      $\phi_f$, $\phi_c$ & porosity: $\phi_f = \vert Y_f \vert / \vert Y \vert$, $\phi_c = \vert Y_c \vert / \vert Y \vert$\\       
        $\aYs{\cdot}{\cdot}$ & ``elasticity'' bilinear form\\
        $\aYf{\cdot}{\cdot}$ & ``viscosity'' bilinear form\\
        $\psibf^k$, $\pi^k$ & velocity and pressure correctors in $Y_f$ \\
        $\omegabf^{ij}$, $\omegabf^c$, $\omegabf^f$ & strain and pressure correctors in $Y_s$\\
        $\Hop$ &  homogenized coefficients (tensor, or scalar), generic form\\       
        $\Kb$ & homogenized permeability tensor\\ 
        $\Cop$ & homogenized elasticity tensor\\
        $\Bb^f$, $\Bb^c$ & homogenized Biot coupling coefficients\\
        $M^{ff}$, $M^{cc}$, $M^{fc}$  & homogenized Biot modulus\\
        $\delta\ub$, $\delta p_f$, $\delta p_c$ & displacement and pressure increments\\
        $\sbb$ & state variable at time $t$: $\ssb = (\ub,\ul{p})$, where $\ul{p} = (p_f,p_c)$\\
        $\sbb\prevstep$ & state variable at time $t - \Delta t$\\
        $\Psi^t$ & residual function associated with time $t$\\
        \hline
    \end{tabular}
    \color{black}
\end{table}

\begin{table}[ht]
  \caption{Employed function spaces, admissibility sets}\label{tab-notation5}
    \begin{tabular}{|c|l|}
      \hline
        $L^2(\Om)$ & space of square integrable functions\\
        $\Hdb(\Om)$ & Sobolev space of differentiable square integrable vector functions\\
        $\Ucalbf^\veps(\Om_{s}^\veps)$ & , set of admissible displacement field, see \eq{eq-5a}\\
        $\Wcalbf_0^\veps(\Om_f^\veps,\Gamma_\fsi^\veps)$ & space of fluid velocity field, see \eq{eq-5a}\\
        $\Qcal_c^\veps(\Om_c^\veps)$ & space of piecewise constant inclusion pressures, see \eq{eq-5a}\\
        $\Qcal_f^\veps(\Om_f^\veps)$  & set of admissible fluid pressure field, see \eq{eq-5a}\\
        $\Ucalbf(\Om)$ & set of admissible displacement field, see \eq{eq-UPspaces}\\
        $\Qcal(\Om)$ & set of admissible fluid pressure field, see \eq{eq-UPspaces}\\
        \hline
    \end{tabular}
    \color{black}
\end{table}

\color{black}

\end{appendices}